\documentclass{birkjour_t2}

\usepackage[utf8]{inputenc}
\usepackage{amsmath,amssymb,paralist,babel}
\usepackage{latexsym,amsfonts,amscd}
\usepackage{graphicx}
\usepackage{dsfont,stix}
\usepackage{comment,doi}
\usepackage{tikz}
\usepackage{caption}
\usepackage{subcaption}
\usepackage{comment}
\usepackage{mathtools}
\usetikzlibrary{arrows, backgrounds}
\usepackage{enumitem}
\usepackage{stix}
\usepackage{tikz,pgfplots}
\pgfplotsset{compat=1.18}

\usepackage{soul}
\usepackage{placeins}

\usepackage{color}
\usepackage[dvipsnames]{xcolor}
\definecolor{LinkColor}{rgb}{0,0,1}
\definecolor{LinkColor2}{rgb}{0,0.5,0}
\definecolor{lbcolor}{rgb}{0.85,0.85,0.85}
\definecolor{FrameColor}{rgb}{0.85,0.85,0.85}
\definecolor{rosso}{rgb}{0.8,0,0}
\definecolor{lightgray}{rgb}{0.5,0.5,0.5}
\definecolor{violet}{rgb}{0.65,0,0.65}
\definecolor{darkgreen}{rgb}{0,0.5,0}

\usepackage{hyperref} 
\hypersetup{%
	colorlinks	=true,% Definition der Links im PDF File
	linkcolor	=LinkColor,%
	anchorcolor	=LinkColor,%
	citecolor	=LinkColor2,%
	filecolor	=LinkColor,%
	menucolor	=LinkColor,%
	urlcolor	=LinkColor,%
}

\makeatletter
\renewenvironment{proof}[1][\proofname]{%
	\par\pushQED{\qed}\normalfont%
	\topsep6\p@\@plus6\p@\relax
	\trivlist\item[\hskip\labelsep\bfseries#1\@addpunct{.}]%
	\ignorespaces
}{%
	\popQED\endtrivlist\@endpefalse
}
\makeatother

\makeatletter
\renewcommand\paragraph{\@startsection{paragraph}{4}{\z@}%
	{1ex \@plus1ex \@minus.2ex}%
	{-1em}%
	{\indent\normalfont\normalsize\bfseries}
    }
    
\renewcommand\subparagraph{\@startsection{paragraph}{4}{\z@}%
	{1ex \@plus1ex \@minus.2ex}%
	{-1em}%
	{\indent\normalfont\normalsize\itshape}
    }
\makeatother

\newtheorem{theo}{Theorem}[section]
\newtheorem{example}[theo]{Example}
\newtheorem{remark}[theo]{Remark}

\newtheorem{corrol}[theo]{Corollary}
\newtheorem{lemma}[theo]{Lemma}
\newtheorem{definition}[theo]{Definition}
\numberwithin{equation}{section}

\def\beq{\begin{equation}}
\def\beqn{\begin{eqnarray*}}
\def\eeq{\end{equation}}
\def\eeqn{\end{eqnarray*}}
\def\beqa{\begin{eqnarray}}
\def\eeqa{\end{eqnarray}}

\def\R{\mathbb{R}}
\def\eps{\varepsilon}

\def\Z{\mathbb{Z}}

\def\L{\mathcal{L}}
\def\N{\mathbb{N}}

\usetikzlibrary{arrows, backgrounds}

\newcommand{\abs}[1]{| #1 |}

\newcommand{\norm}[1]{\| #1 \|}
\newcommand{\bignorm}[1]{\big\| #1 \big\|}

\newcommand{\ang}[2]{ \langle #1 , #2  \rangle}

\newcommand{\scp}[2]{ \left( #1 , #2  \right)}

\newcommand{\mean}[1]{\langle #1 \rangle}

\newcommand{\del}{\partial}
\newcommand{\delt}{\partial_{t}}

\newcommand{\intO}{\int_\Omega}

\newcommand{\dx}{\;\mathrm{d}x}
\newcommand{\dt}{\;\mathrm dt}
\newcommand{\ds}{\;\mathrm ds}

\newcommand{\ddt}{\frac{\mathrm d}{\mathrm dt}}

\newcommand{\dds}{\frac{\mathrm d}{\mathrm ds}}

\newcommand{\Grad}{\nabla}
\newcommand{\Lap}{\Delta}
\newcommand{\Div}{\textnormal{div}}

\newcommand{\emb}{\hookrightarrow}

\newcommand{\doublehookrightarrow}{%
    \mathrel{\mathrlap{{\mspace{4mu}\lhook}}{\hookrightarrow}}
}

\newcommand{\skp}[2]{\left({#1\, ,\, #2}\right)}
\usepackage{array}  %% for \newcolumntype{...}[...]{...}
\newcolumntype{M}[1]{>{\centering\arraybackslash}m{#1}}
\newcolumntype{L}[1]{>{\centering\arraybackslash}m{#1}} % for labels
\newcolumntype{P}[1]{>{\centering\arraybackslash}p{#1}} % for images

\begin{document}

\title[Pattern formation in biomembranes]{Curvature-driven pattern formation in biomembranes:\\ A gradient flow approach}
\author{Patrik Knopf}
\address{Faculty of Mathematics, University of Regensburg, 93053 Regensburg, Germany \\
and Department of Mathematics, Karlsruhe Institute of Technology (KIT), 76128 Karlsruhe, Germany.}
\email{\href{mailto:patrik.knopf@ur.de}{patrik.knopf@ur.de}, \href{mailto:patrik.knopf@kit.edu}{patrik.knopf@kit.edu}}

\author{Anastasija Pe\v{s}i\'{c}}
\address{Weierstrass Institute for Applied Analysis and Stochastics, 10117 Berlin, Germany.}
\email{\href{mailto:anastasija.pesic@wias-berlin.de}{anastasija.pesic@wias-berlin.de}}

\author{Dennis Trautwein}
\address{Faculty of Mathematics, University of Regensburg, 93053 Regensburg, Germany.}
\email{\href{mailto:dennis.trautwein@ur.de}{dennis.trautwein@ur.de}}
\date{\today}



\begin{abstract}
In this work, we study a phase-field model for curvature-driven pattern formation in biomembranes. The model is derived as a gradient flow of an energy functional that approximates the two-phase Canham--Helfrich energy. This leads to a Cahn--Hilliard-type equation with cross diffusion for the relative chemical concentration of one lipid phase, coupled to a fourth-order reaction-diffusion equation describing the height profile of the membrane.
We first prove the existence of weak solutions for the case of regular double-well potentials, using a minimizing movement scheme to construct approximate solutions. The analysis is then extended to singular potentials, e.g., the Flory--Huggins potential, by approximating them with a Moreau--Yosida regularization. 
For both cases, we establish higher regularity, continuous dependence on the initial data, and consequently the uniqueness of weak solutions.
Finally, we propose a well-posed finite element discretization of the model and present numerical experiments illustrating the effect of different physical parameters on the resulting membrane patterns. Depending on the parameter regime, we observe purely striped, dotted, or snake-like structures.
\end{abstract}

\maketitle

%\begin{small}
%\setcounter{tocdepth}{3}
%\hypersetup{linkcolor=black}
%\tableofcontents
%\end{small}
%\newpage 

\setlength\parskip{1ex}

\section{Introduction}

\begin{figure}
\centering
\includegraphics[width=0.8\linewidth]{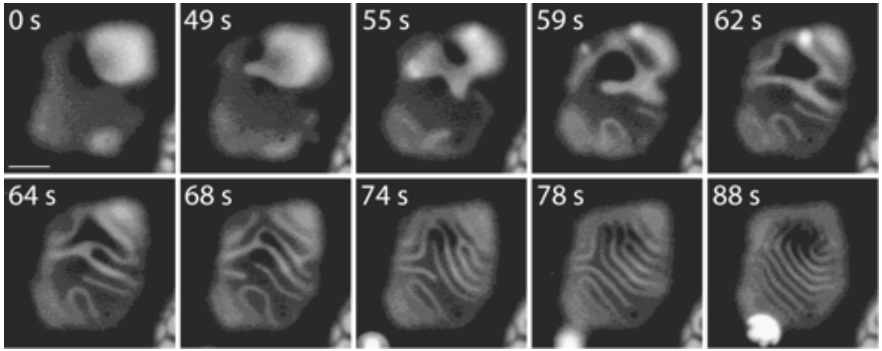}
\caption{Epifluorescence microscopy images showing a circular domain undergoing shape instability, evolving into branched arms and eventually forming a stripe-like structure. 
Reprinted with permission from \cite{rozovsky2005formation}. % J. Am. Chem. Soc. 2005, 127, 1, 36--37. 
Copyright~2005 American Chemical Society.
%
% \DT{We now have the license for arxiv. Do not forget to obtain the permission for reuse from the \href{https://pubs.acs.org/doi/full/10.1021/ja046300o}{\textit{publisher}} when submitting to a journal.}
}\vspace{-1em}
\label{fig:img_rozovsky}
\end{figure}

We investigate a gradient flow equation associated with the energy of a two-phase biomembrane, first introduced in \cite{Komura:Langmuir:2006}. In most living organisms, the cellular membrane consists of a lipid bilayer composed of various lipid types. These are believed to organize into domains known as lipid rafts \cite{simons-ikonen:1997}. Such lipid domains have been observed both in living cells \cite{hurst2020lipid, toulmay2013direct, furukawa2022signaling, bramkamp2015exploring} and in model membranes \cite{rozovsky2005formation}, where patterns such as stripes and dots have been reported.
An example of this behavior is shown in Figure~\ref{fig:img_rozovsky} (reproduced from \cite{rozovsky2005formation}), where a circular lipid domain undergoes a shape instability and develops into a stripe-like superstructure.
%\tbd{Here more references about models for pattern formation.}

To model this phenomenon, \cite{Komura:Langmuir:2006} proposed an energy functional that approximates the two-phase Canham--Helfrich functional. We consider a generalization of this model, given by \eqref{DEF:F}, where the physical parameters are matrix-valued rather than scalars. This allows us to enforce the formation of stripes in broader regimes. Specifically, for a dimension $d\in \N$ and domain $\Omega := \mathbb{T}^d$, we consider the energy functional:
\begin{equation}
    \label{DEF:F}
    F(u,h) 
    = \intO \left( \frac{1}{\eps} W(u)
        + \frac{\eps}{2} \abs{\Grad u}^2 
        + \frac 12 G \Grad h \cdot \Grad h
        + \frac\kappa2 \abs{\Delta h}^2
        - L \Grad h \cdot \Grad u \right) \dx,
\end{equation}
where $u$ represents the relative chemical concentration of one lipid phase, and $h$ denotes the height profile of the membrane.
The parameters $\eps, \kappa > 0$ are the diffuse interface width and bending rigidity, while $G, L \in \R^{d\times d}$ represent the surface tension and coupling strength, respectively.
Moreover, $W$ is a double-well potential with minima at (or near) $\pm1$.

Using a stability analysis for the isotropic version of the energy (i.e., $G=\sigma \mathrm{I}$ and $L=\Lambda \mathrm{I}$ with $\sigma,\Lambda>0$, where $\mathrm{I}\in \mathbb{R}^{d\times d}$ denotes the identity matrix) the authors of \cite{Komura:Langmuir:2006} identified parameter regimes where various patterns, such as stripes and dots, may form. They classified regions in the parameter space of coupling strength and surface tension where a homogeneous initial state becomes unstable with respect to long-wave perturbations. The remaining parameters, such as bending rigidity $\kappa$ and interface constant $\varepsilon$, were chosen based on experimental data from the literature.
One of our aims is to numerically confirm that different pattern types indeed emerge in these regimes.

Another goal of this work concerns the analytical understanding of such models. We briefly summarize the current state of the literature.
The existence of minimizers for the energy functional \eqref{DEF:F} with scalar parameters, as well as the scaling behavior of the minimal energy with respect to the problem parameters, was analyzed in \cite{ginster2024nonlinear}. That work suggested that for large coupling parameters, uniform minimizers are not energetically favorable, whereas in the weak coupling regime, uniform structures are energetically preferred. % In that work, the constructions used for proving the upper energy bound also helped to identify parameter regimes in which non-uniform structures are expected. It was suggested that for large coupling parameters, uniform minimizers are not energetically favorable, whereas in the weak coupling regime, uniform structures are energetically preferred.
The latter was rigorously confirmed for certain sub-regimes in \cite{FonHayLeoZwi, ginster2024sharp, brazke2023gamma}, where the $\Gamma$-limit of a related nonlocal family was shown to be a multiple of a perimeter functional. Furthermore, it was recently proven in \cite{daneri2024rigorous} that, in a critical parameter regime, minimizers of a related sharp-interface model are uniform periodic stripes.

% Another goal of us is concerning the analysis of such a model. Here we shortly discuss the state of the art in the literature.
% The existence of minimizers for a functional with scalar parameters, as well as the scaling behavior of the infimal energy with respect to problem parameters, was studied in \cite{ginster2024nonlinear}. There, the constructions in the proof of the upper bound were used to identify the parameter regimes in which non-uniform structures are expected to occur. It was suggested that for large coupling parameters, uniform minimizers are not energetically favorable. Conversely, in the weak coupling regime, uniform structures are expected to be energetically preferred. This has been rigorously shown in \cite{FonHayLeoZwi, ginster2024sharp, brazke2023gamma} in certain sub-regimes, where the $\Gamma$-limit of a related nonlocal family was computed to be the multiple of a perimeter functional. Moreover, recently it was shown in \cite{daneri2024rigorous}, that in a critical regime, the minimizers of a related sharp-interface model with a non-local term are uniform periodic stripes.

In the present paper, we adopt a gradient flow approach to study the existence and stability of (local) minimizers of $F$. We derive a system of partial differential equations from the energy functional \eqref{DEF:F}, relying on the principle that the long-time behavior of gradient flows is closely related to the energy minimizers (see, e.g., \cite{ball1991dynamics}).
More precisely, for $T>0$ and $\Omega_T := (0,T)\times \mathbb{T}^d$ we consider the system of equations:
\begin{subequations}
    \label{PDE}
    \begin{alignat}{2}
        \label{PDE:1}
        &\delt u = \Lap \mu
        &&\quad\text{in $\Omega_T$},
        \\
        \label{PDE:2}
        &\mu = - \eps \Lap u + \frac{1}{\eps} W'(u) + \Div(L\Grad h)
        &&\quad\text{in $\Omega_T$},
        \\
        \label{PDE:3}
        &\delt h = \Div(G\Grad h) - \kappa \Delta^2 h + \Div(L\Grad u)
        &&\quad\text{in $\Omega_T$},
    \end{alignat}
\end{subequations}
with initial conditions $u(0)=u_0$ and $h(0)=h_0$ in $\Omega$. This system is derived from \eqref{DEF:F} (see Section \ref{sec:modelling}), by taking an $H^{-1}$-gradient flow with respect to the first and an $L^2$-gradient flow with respect to the second argument in the energy $F(u,h)$. The result is a coupled system involving a Cahn--Hilliard-type equation for $u$. 
We choose the volume-preserving $H^{-1}$-flow for $u$ because the lipid molecules are assumed to diffuse along the membrane, conserving the relative lipid concentration.

% \tbd{Here goes a long list of references for previous works on both CH and AC?}
% \tbd{Discuss differences to \cite{GLNS_2024_pattern}, where a Cahn--Hilliard--Swift--Hohenberg system was used to describe pattern formation. The energy and the model equations are similar to our model.}

We note that in the case of isotropic surface tension ($G=\sigma \mathrm{I}$), equation \eqref{PDE:3} can be reformulated as:
\begin{align}
    \label{PDE:SH}
    \partial_t h + \kappa \left[ 
        \left(\Delta - \frac{\sigma}{2\kappa}\right)^2 
        - \frac{\sigma^2}{4\kappa^2} 
        \right] h
        = \Div(L \nabla u)
        \quad\text{in $\Omega_T$.}
\end{align}
If $\sigma<0$, this equation is an inhomogeneous, linear version of the \textit{Swift--Hohenberg equation}. This equation is commonly used to describe pattern formation because the resulting negative diffusion coefficient creates an instability that causes the membrane to wrinkle on its own. For instance, a Cahn--Hilliard--Swift--Hohenberg model of this type was proposed in \cite{Martinez-Agustin2022} and further investigated in \cite{GLNS_2024_pattern}.
However, since we consider $\sigma>0$, equation \eqref{PDE:SH} is not strictly a Swift--Hohenberg equation. A positive surface tension serves as a stabilizing term that tries to keep the membrane flat. Nevertheless, we still expect patterns to form. This occurs not due to the structure of equation \eqref{PDE:3}, but because the lipid concentration $u$ forces the membrane to change shape through the coupling term involving $L$.
% Models involving the Swift--Hohenberg equation are commonly used in the description of pattern formation processes. For example, a Cahn--Hilliard--Swift--Hohenberg model was proposed in \cite{Martinez-Agustin2022} and further investigated in \cite{GLNS_2024_pattern}. 
%In the present paper, as we consider $\sigma>0$, equation \eqref{PDE:SH} cannot be identified as a Swift--Hohenberg equation in the strict sense. Nevertheless, we expect the formation of patterns in our model, which is not only due to the structure of equation \eqref{PDE:3} but occurs mainly due to the interaction between the order parameter $u$ and the height function $h$ through the terms involving $L$.

\textbf{Outline of the main results.}
Our main analytical results are presented in Theorems~\ref{THM:WS:REG} and~\ref{THM:WS:SING}, where we prove the existence of weak solutions to \eqref{PDE} for regular and singular potentials~$W$, respectively. In Theorem~\ref{THM:ContDep:REG}, we further show the continuous dependence on the initial data, which in turn implies uniqueness.
The precise definitions of weak solutions for the regular and singular cases are given in Definitions~\ref{DEF:WS:REG} and~\ref{DEF:WS:SING}. In addition, we obtain higher regularity results, stated in Corollary~\ref{COR:HIGHREG} and Theorem~\ref{THM:WS:SING}.
% Our main analytical results are stated in Theorem~\ref{THM:WS:REG} and Theorem~\ref{THM:WS:SING}, in which we prove the existence of weak solutions to \eqref{PDE} for regular and singular potentials $W$, respectively. This requires an adequate definition of a weak solution in both cases, which are given by Definitions~\ref{DEF:WS:REG} and~\ref{DEF:WS:SING}. 
\\
%\PK{Why, in contrast to the literature, we also consider singular potentials:}
Originally, the model was proposed with a smooth (regular) double-well potential $W$ (cf.~\cite{Komura:Langmuir:2006}). However, this does not ensure that the order parameter $u$ stays within its physical bounds $u \in [-1,1]$. In contrast, singular potentials inherently enforce this constraint. This makes the model more physically consistent and prevents nonphysical values of the phase-field variable. For that reason, we also study singular potentials in this work.
\\
To prove the existence results for both regular and singular potentials, we follow two related approaches.
For the proof in the case of regular potentials (Theorem~\ref{THM:WS:REG}), we employ a minimizing movement scheme that reflects the gradient flow structure of \eqref{PDE}. Using the direct method of calculus of variations, we first show that each time step of the scheme admits a minimizer. By exploiting the gradient-flow formulation, we rewrite the scheme as a time-discrete PDE system. Uniform estimates and compactness arguments then allow us to extract a convergent subsequence whose limit is a weak solution in the sense of Definition~\ref{DEF:WS:REG}. The corresponding energy dissipation law follows by combining these compactness results with lower semicontinuity and a weak-strong convergence argument. Finally, Corollary~\ref{COR:HIGHREG} provides improved regularity for the weak solutions.
\\
For singular potentials (Theorem~\ref{THM:WS:SING}), we approximate the singular potential using a Moreau--Yosida regulari\-zation. The existence of a solution to the regularized system then follow directly from Theorem~\ref{THM:WS:REG}. Using uniform estimates and compactness arguments, we pass to the limit in the regularization and thereby establish the existence of a solution to the system with singular potentials in the sense of Definition~\ref{DEF:WS:SING}. The proof of the energy dissipation law proceeds in the same way as for regular potentials. Higher regularity is obtained using tools from elliptic regularity theory.
\\
For the uniqueness of weak solutions, we show continuous dependence on the initial data in both the regular and singular potential settings. In particular, if the subdifferential of the potential is single-valued (as for regular potentials) and the initial data coincide, then the continuous dependence estimate implies that the corresponding weak solutions must coincide and hence uniqueness holds. This is summarized in Theorem~\ref{THM:ContDep:REG}.

\textbf{Structure of the paper.}
Our paper is organized as follows. In Section~\ref{sec:prelim}, we introduce the notation and assumptions used throughout the work, including the assumptions on the potential, as well as the Moreau--Yosida regularization needed for singular potentials. In Section~\ref{sec:modelling}, we derive the system \eqref{PDE} as a gradient flow of the energy \eqref{DEF:F}. Our main results are stated in Section~\ref{sec:mainresult}.
Sections~\ref{SEC:proof_reg} and~\ref{sec:proof_sing} contain the proofs of the existence results for regular (Theorem~\ref{THM:WS:REG}) and singular potentials (Theorem~\ref{THM:WS:SING}), respectively.
Higher regularity for weak solutions in the regular case is proved in Section~\ref{SEC:proof_highreg} (Corollary~\ref{COR:HIGHREG}), while the corresponding result for singular potentials is included in Section~\ref{sec:proof_sing}. The continuous dependence and uniqueness result (Theorem~\ref{THM:ContDep:REG}) is shown in Section~\ref{SEC:proof_reg_ContDep}.
Finally, Section~\ref{sec:num} introduces a well-posed finite element discretization of the model and presents numerical experiments that show the different pattern types observed for various parameter choices.

% \newpage

\section{Notations, assumptions and preliminaries}\label{sec:prelim}

\subsection{General assumptions}
We start with introducing assumptions on the domain and model parameters.
\begin{enumerate}[label=\textnormal{\bfseries (A\arabic*)}, topsep=1ex, itemsep = 1ex, leftmargin=*]
    \item \label{ASS:DOM}
    We choose $\Omega$ to be the flat torus $\mathbb{T}^d \coloneqq [\R\setminus\Z]^d$ with dimension $d\in\N$ with $d\ge 2$. For an arbitrary final time $T>0$ we write
    \begin{equation*}
        \Omega_T \coloneqq (0,T) \times \Omega.
    \end{equation*}

    \item \label{ASS:COEFF}
    We assume that $\eps,\kappa > 0$ are positive constants, while $G, L \in\R^{d\times d}$ are symmetric, positive definite matrices. Hence, there exist constants $G^*,G_*,L^*,L_*>0$ such that for all $\zeta, \xi\in\R^d$, it holds
    \begin{equation*}
        \xi^T G \xi \ge G_* \abs{\xi}^2,
        \quad
        \xi^T L \xi \ge L_* \abs{\xi}^2,
        \quad
        \zeta^T G \xi \le G^* \abs{\zeta} \abs{\xi},
        \quad
        \zeta^T L \xi \le L^* \abs{\zeta} \abs{\xi}.
    \end{equation*}
\end{enumerate}

\bigskip

\noindent Next, let us introduce some notations that will be used throughout this paper.

\begin{enumerate}[label=\textnormal{\bfseries (N\arabic*)},topsep=1ex, itemsep = 1ex, leftmargin=*]

\item \textbf{Notation for general Banach spaces.} 
For any normed space $X$, we denote its norm by $\|\cdot\|_X$,
its dual space by $X^*$ and the duality pairing between $X^*$ and $X$ by $\langle\cdot,\cdot\rangle_X$.
Besides, if $X$ is a Hilbert space, we write $(\cdot,\cdot)_X$ to denote the corresponding inner product.
% Furthermore, for any vector space $X$, corresponding spaces of vector-valued or matrix-valued functions with each component belonging to $X$ are denoted by $\mathbf{X}$.

\item \textbf{Lebesgue and Sobolev spaces.} 
Suppose that \ref{ASS:DOM} holds.
For $1 \leq p \leq \infty$ and $k \in \N$, the standard Lebesgue spaces and Sobolev spaces of periodic functions on the flat torus $\Omega$ are denoted by $L^p(\Omega)$ and $W^{k,p}(\Omega)$, respectively, and their standard norms are denoted by $\|\cdot\|_{L^p}$ and $\|\cdot\|_{W^{k,p}}$, respectively.
In the case $p = 2$, we use the notation $H^k(\Omega) = W^{k,2}(\Omega)$. %We point out that $H^k(\Omega)$ is a Hilbert space, and that $H^0(\Omega)$ can be identified with $L^2(\Omega)$.
%For brevity, we simply write 
%\begin{equation*}
%    (\cdot,\cdot) \coloneqq (\cdot,\cdot)_{L^2},
%    \quad
%    \|\cdot\|\coloneqq\|\cdot\|_{L^2}
%    \quad
%    \langle\cdot,\cdot\rangle \coloneqq \langle\cdot,\cdot\rangle_{H^1(\Omega)}.
%\end{equation*}
Moreover, for any interval $I\subset\R$, any Banach space $X$, $1 \leq p \leq \infty$ and $k \in \N$, we write $L^p(I;X)$, $W^{k,p}(I;X)$ and $H^{k}(I;X) = W^{k,2}(I;X)$ to denote the Lebesgue and Sobolev spaces of functions with values in $X$. The standard norms are denoted by $\|\cdot\|_{L^p(I;X)}$, $\|\cdot\|_{W^{k,p}(I;X)}$ and $\|\cdot\|_{H^k(I;X)}$, respectively. 
For any $k\in\N$, we further write $H^{-k}(\Omega)$ to denote the dual space of $H^k(\Omega)$.  

\item \textbf{Spaces of continuous functions.}
For any interval $I\subset\R$ and any Banach space $X$, $C(I;X)$ denotes the space of continuous functions mapping from $I$ to $X$.
% and $C_b(I;X)$ denotes the space of functions in $C(I;X)$, which are additionally bounded. 
% Moreover, $C_\mathrm{w}(I;X)$ denotes the space of functions mapping from $I$ to $X$, which are continuous on $I$ with respect to the weak topology on $X$, and $BC_\mathrm{w}(I;X)$ denotes the space of functions in $C_\mathrm{w}(I;X)$, which are additionally bounded.

\item\label{N:MEAN} \textbf{Spaces of functions with prescribed mean.}
Let $m\in \R$ be arbitrary.
For any distribution $\varphi\in H^{-1}(\Omega)$, its generalized spatial mean is defined as
\begin{equation*}
    \mean{\varphi}\coloneqq \langle \varphi,1 \rangle_{H^1(\Omega)}.
\end{equation*}%
Using this definition, we introduce the affine spaces
\begin{align*}
    H^{-1}_{(m)}(\Omega) &\coloneqq \big\{ u\in H^1(\Omega)^* \,:\, \mean{u} = m \big\} \subset H^{-1}(\Omega),\\
    L^2_{(m)}(\Omega) &\coloneqq \big\{ u\in L^2(\Omega) \,:\, \mean{u} = m \big\} \subset L^2(\Omega),\\
    H^k_{(m)}(\Omega) &\coloneqq \big\{ u\in H^k(\Omega) \,:\, \mean{u} = m \big\} \subset H^k(\Omega)
\end{align*}
for any $k\in\N$.
If $m=0$, these spaces are closed linear subspaces of the respective Hilbert space. Hence, they are also Hilbert spaces. On $H^{-1}_{(0)}(\Omega)$, we introduce the bilinear form
\begin{equation*}
    \scp{f}{g}_{-1} \coloneqq \intO \Grad(-\Lap)^{-1} f \cdot \Grad(-\Lap)^{-1} g \dx
    \quad\text{for all $f,g \in H^{-1}_{(0)}(\Omega)$},
\end{equation*}
where 
\begin{equation*}
    (-\Lap)^{-1}:H^{-1}_{(0)}(\Omega) \to H^1_{(0)}(\Omega),
    \quad
    f \mapsto (-\Lap)^{-1}f = u_f 
\end{equation*}
denotes the solution operator of Poisson's equation
\begin{equation*}
    -\Lap u_f = f \quad\text{a.e. in $\Omega$}
\end{equation*}
with periodic boundary conditions. 
The bilinear form $\scp{\cdot}{\cdot}_{-1}$ defines an inner product on $H^{-1}_{(0)}(\Omega)$. The induced norm reads as
\begin{equation*}
    \norm{\cdot}_{-1} = \scp{\cdot}{\cdot}_{-1}^{1/2},
\end{equation*}
and the space $$\big(H^{-1}_{(0)}(\Omega), \, \scp{\cdot}{\cdot}_{-1}, \, \norm{\cdot}_{-1}\big)$$ is a Hilbert space. We point out that $\scp{\cdot}{\cdot}_{-1}$ also defines an inner product on the spaces $L^2_{(0)}(\Omega)$ and $H^1_{(0)}(\Omega)$, but these spaces are not complete with respect to the induced norm.
\end{enumerate}  

\subsection{Assumptions on the potential}
In the following, we address both regular and singular potentials. By regular potentials we mean functions that satisfy the assumption stated below.
\begin{enumerate}[label=\textnormal{\bfseries (R)}, topsep=1ex, leftmargin=*]
    \item\label{ASS:POT:REG} The potential $W\colon\, \R\to [0,+\infty)$ can be decomposed as $W=W_1+W_2$, where functions $W_1\colon\, \R\to [0,+\infty)$ and $W_2\colon\, \R\to \R$ have the following properties.
    \begin{enumerate}[label=\textnormal{\bfseries (R\arabic*)}, topsep=1ex, itemsep=1ex, leftmargin=*]
        \item \label{ASS:POT:REG:1}
        $W_1$ is continuously differentiable and convex. Without loss of generality, $W_1(0)=0$. 
        Moreover, there exists $p\in\R$ with $2 \le p < \infty$ if $d=2$ and $2 \le p \leq \tfrac{2d}{d-2}$ if $d\geq3$, as well as
        a constant $C_{1}>0$ such that
        \begin{align}
            \label{COND:W'}
            \abs{W_1'(s)} &\le C_{1} (1+\abs{s}^{p-1}), \quad\text{ for all } s\in\R.
        \end{align}
        As a consequence, there exists a constant $C_{0}>0$ such that
        \begin{align}
            \label{COND:W}
            \abs{W_1(s)} &\le C_{0} (1+\abs{s}^{p}), \quad\text{ for all } s\in\R.
        \end{align}
    
        \item \label{ASS:POT:REG:2}
        $W_2$ is continuously differentiable and its derivative $W_2'$ is Lipschitz continuous with Lipschitz constant $\gamma\geq0$. This implies that there exist constants $A>0$ and $B\in\R$ such that
        \begin{equation*}
            \abs{W_2(s)} \le A s^2 + B \quad\text{for all $s\in\R$.}
        \end{equation*}

        \item \label{ASS:POT:REG:3}
        For $A$ as in \ref{ASS:POT:REG:2} and 
        \begin{align*}
            D\coloneqq 1 + A + \frac{\eps(L^*)^2}{\kappa}
        \end{align*}
        there exists $B_D\in\R$ such that 
        \begin{equation}
            \label{COND:W:2}
            W_1(s) \ge D s^2 + B_D
            \quad\text{for all $s\in\R$.}
        \end{equation}
    \end{enumerate}
\end{enumerate}

On the other hand, singular potentials are assumed to satisfy the following condition.
\begin{enumerate}[label=\textnormal{\bfseries (S)}, topsep=1ex, leftmargin=*]
    \item\label{ASS:POT:SING} The potential $W\colon\, \R\to [0,+\infty]$ can be decomposed as $W=W_1+W_2$, where functions 
    $W_1\colon\, \R\to [0,+\infty]$ and $W_2\colon\, \R\to \R$ have the following properties.
    \begin{enumerate}[label=\textnormal{\bfseries (S\arabic*)}, topsep=1ex, itemsep=1ex, leftmargin=*]
        \item \label{ASS:POT:SING:1}
            $W_1$ is a proper, convex and lower semicontinuous function, which satisfies $W_1(0) = 0$ and
            \begin{equation}\label{COND:COERC}
                \underset{\abs{s}\to\infty}{\lim}\; \frac{W_1(s)}{\abs{s}^2} = +\infty.
            \end{equation}
        \item \label{ASS:POT:SING:2}
            $W_2$ is continuously differentiable and its derivative $W_2'$ is Lipschitz continuous with Lipschitz constant $\gamma\geq0$. This implies that there exist constants $A>0$ and $B\in\R$ such that
            \begin{equation*}
                \abs{W_2(s)} \le As^2 + B \quad\text{for all $s\in\R$.}
            \end{equation*}
    \end{enumerate}
\end{enumerate}
Note that the condition \ref{ASS:POT:SING:1} already entails that $W_1$ is bounded from below. The effective domain of $W_1$ reads as
\begin{equation*}
    D(W_1) = \big\{ s\in\R \,\vert\, W_1(s) < +\infty \big\}.
\end{equation*}
We write 
\begin{equation*}
    \mathrm{w}_1 \coloneqq \del W_1
\end{equation*}
to denote the subdifferential of $W_1$. It is well known that $\mathrm{w}_1$ is a maximal monotone graph in $\R\times\R$ (see, e.g., \cite[Theorem~2.8]{Barbu2010}). Its effective domain is given by
\begin{equation*}
    D(\mathrm{w}_1) = \big\{ s\in D(W_1) \,\vert\, \mathrm{w}_1(s) \neq \emptyset \big\}.
\end{equation*}
%$\phantom{x}$
%and we assume that $0\in D(\mathrm{w}_1)$.

\medskip

\begin{figure}[h]
\centering
\begin{tikzpicture}[scale=0.65, transform shape]
\begin{axis}[
axis x line=middle, 
axis y line=middle,
xlabel=$s$, ylabel=$W_\mathrm{reg}(s)$,
xlabel style={at={(axis description cs:0.95,-0.05)}, anchor=north}, 
ylabel style={at={(axis description cs:0.48,0.85)}, anchor=east},
xtick={-1, 1},
xticklabels={-1, 1},
ymajorticks=false,
domain=-1.6:1.6
]
\addplot[samples=100, smooth, sharp plot] {10*(x-1)*(x-1)*(x+1)*(x+1)};
\end{axis}
\end{tikzpicture}
\quad 
%%%%%
\begin{tikzpicture}[scale=0.65, transform shape,
declare function={ 
F(\x) = 0.5*((1-\x)*ln(1-\x) + (1+\x)*ln(1+\x));    
}
]
\begin{axis}[
axis x line=middle, 
axis y line=middle,
xlabel=$s$, ylabel=$W_\mathrm{log}(s)$,
xlabel style={at={(axis description cs:0.95,-0.05)}, anchor=north}, 
ylabel style={at={(axis description cs:0.48,0.85)}, anchor=east},
xtick={-1,1},
xticklabels={-1, 1},
xmin=-1.3,xmax=1.3,
ymin=0.55, ymax=0.7,
ymajorticks=false,
domain=-1:1
]
\addplot[samples=100, smooth, sharp plot] {0.5*((1-x)*ln(1-x)+(1+x)*ln(1+x)) + 0.65*(1-x*x) - 0.05};
\addplot[dashed, sharp plot] coordinates {(1,0.5) (1,0.63)};
\addplot[dashed, sharp plot] coordinates {(-1,0.5) (-1,0.63)};
\end{axis}
\end{tikzpicture}
\quad
%%%%%
\begin{tikzpicture}[scale=0.65, transform shape]
\begin{axis}[
axis x line=middle, 
axis y line=middle,
xlabel=$s$, ylabel=$W_\mathrm{obst}(s)$,
xlabel style={at={(axis description cs:0.95,-0.05)}, anchor=north}, 
ylabel style={at={(axis description cs:0.48,0.85)}, anchor=east},
xtick={-1, 1},
xticklabels={-1, 1},
ymajorticks=false,
ymin=0, ymax=2.5,
xmin=-1.5, xmax=1.5,
domain=-1:1
]
\addplot[samples=100, smooth, sharp plot] {(1-x*x)};
\addplot[samples=100, smooth, sharp plot,domain=-1.5:-1] {2.4};
\addplot[samples=100, smooth, sharp plot,domain=1:1.5] {2.4};
\addplot[dashed, sharp plot] coordinates {(1,0) (1,2.4)}
node [pos=1, below right] {$+\infty$};
\addplot[dashed, sharp plot] coordinates {(-1,0) (-1,2.4)}
node [pos=1, below left] {$+\infty$};
\end{axis}
\end{tikzpicture}
\caption{Sketch of the polynomial double-well potential $W_\mathrm{reg}$, the logarithmic Flory--Huggins potential $W_\mathrm{log}$, and the double-obstacle potential $W_\mathrm{obst}$, from left to right.}
\label{fig:W_examples}
\end{figure}
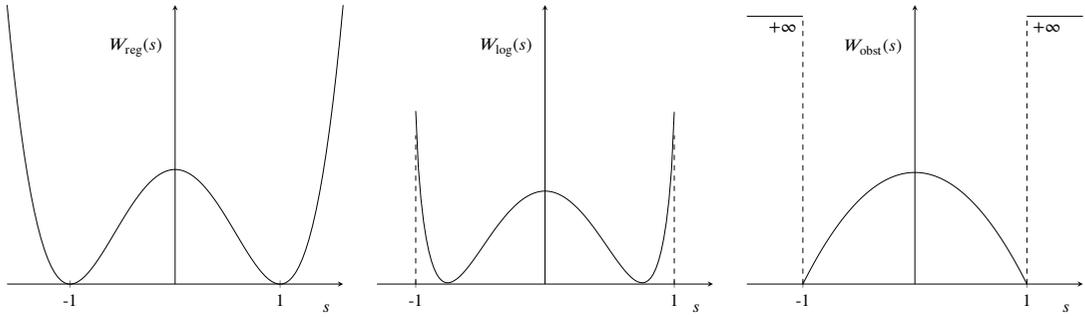

\begin{example} \normalfont We list some typical potentials and those relevant for the physical applications. For illustration, see Figure~\ref{fig:W_examples}.
\begin{enumerate}[topsep=1ex,itemsep=1ex]
\item[(i)] Polynomial potential:
\begin{equation}\label{DEF:WREG}
W_\mathrm{reg}\colon\, \R\to\R, \quad
W_\mathrm{reg}(s) = a_4 s^4 - a_2 s^2 - \mu_0 s + a_0
\end{equation}
for some constants $a_4,a_2>0$ and $a_0,\mu_0\in\R$.\\
Note that $W_\mathrm{reg}$ satisfies Assumption~\ref{ASS:POT:REG} with $p=4$ if $d\le 4$.

\item[(ii)] Flory--Huggins potential:
\begin{align}\label{DEF:WLOG}
    W_\mathrm{log}\colon\, [-1,1]\to\R,\quad
    W_\mathrm{log}(s) = \frac{\theta}{2} \big( (1+s)\ln(1+s) + (1-s)\ln(1-s) \big) - \frac{\theta_c}{2} s^2 - \mu_0 s,
\end{align}
with the convention $0 \ln0 \coloneq 0$, for some constants $0<\theta<\theta_c$, and $\mu_0\in\R$.\\
The Flory--Huggins potential satisfies \ref{ASS:POT:SING} with
\begin{align*}
    &W_{\mathrm{log},1}(s) \coloneqq \frac{\theta}{2} \big( (1+s)\ln(1+s) + (1-s)\ln(1-s) \big),
    &&\text{for all $s\in D(W_{\mathrm{log},1}) = [-1,1]$,}
    \\
    &W_{\mathrm{log},1}'(s) \phantom{:}= \frac{\theta}{2} \big( \ln(1+s) - \ln(1-s) \big),
    &&\text{for all $s\in D(\del W_{\mathrm{log},1}) = (-1,1)$}.
\end{align*}

\item[(iii)] Double-obstacle potential:
\begin{equation}\label{DEF:WOBST}
    W_\mathrm{obst}\colon\,\R\to\R, \quad
    W_\mathrm{obst}(s) = \mathcal{I}_{[-1,1]} - a_2 s^2 - \mu_0 s ,
\end{equation}
for some constants $a_4,a_2>0$, $\mu_0\in\R$, where
\begin{equation}
    \mathcal{I}_{[-1,1]}(s) = 
    \begin{cases}
        0, &\text{if $s\in [-1,1]$};\\
        +\infty, &\text{if $s \not\in [-1,1]$}.
    \end{cases}
\end{equation}
The double-obstacle potential $W_\mathrm{obst}$ satisfies \ref{ASS:POT:SING} with
\begin{align*}
    &W_{\mathrm{obst},1}(s) \coloneqq \mathcal{I}_{[-1,1]}(s)
    \qquad\text{for all $s\in D(W_{\mathrm{obst},1}) = [-1,1]$.}
    % \\
    % &W_{\mathrm{obst},1}'(s) \phantom{:}= 0
    % &&\text{for all $s\in (-1,1)$.}
\end{align*}
Moreover, we have $D(\del W_{\mathrm{obst},1}) = [-1,1]$ with $\del W_{\mathrm{obst},1}(-1) = (-\infty,0]$ and $\del W_{\mathrm{obst},1}(1) = [0,\infty)$.
\end{enumerate}
\end{example}

\FloatBarrier
% --------------- FIGURE WITH POTENTIALS: --------------------
%Here are the tikz pictures, see Figure~\ref{fig:W_examples}.
%------------------------------------------------------------------
\FloatBarrier

\subsection{Moreau-Yosida regularization}\label{sec:MorYosReg}

% \tbd{Sketch?}
In the proof of the existence result for the case of singular potentials, we will need the so-called \emph{Moreau-Yosida approximation}. For a (singular) potential $W_1: \R \to [0, +\infty]$ satisfying property \ref{ASS:POT:SING:1}, we define the \emph{Moreau-Yosida regularization} of $W_1$ as
\begin{equation*}
	W_{1,\lambda}:\R\to\R, \quad
	W_{1,\lambda}(r) = \inf_{s\in\R}\left\{\frac{1}{2\lambda}\abs{r-s}^2 + W_1(s)\right\}
\end{equation*}
for every $\lambda>0$. The function $W_{1,\lambda}$ is differentiable and its derivative can be expressed as
\begin{equation*}
	\mathrm{w}_{1,\lambda}:\R\to\R, \quad
	\mathrm{w}_{1,\lambda}(r) = \frac{1}{\lambda}\Big( r - (\mathrm{id}_\R + \lambda \mathrm{w}_1)^{-1}(r)\Big).
\end{equation*}
It is referred to as the \textit{Yosida regularization} of the subdifferential $\mathrm{w}_1$. We point out that even though $\mathrm{w}_1$ might be multi-valued, its Yosida approximation $\mathrm{w}_{1,\lambda}$ is single-valued for every $\lambda>0$.  It is well-known that the Moreau--Yosida approximation has the following properties:
\begin{enumerate}[label=\textnormal{\bfseries (M\arabic*)}, itemsep=1ex, topsep=1ex, leftmargin=*]
	\item \label{M1}
    For every $\lambda>0$, it holds
    $W_{1,\lambda}\in C^{1,1}(\R)$ and $\lambda^{-1}$ is a Lipschitz constant for $W_{1,\lambda}^\prime$. 
	Since $W_1(0) = 0$, it further holds $W_{1,\lambda}(0) = 0$. 
	Consequently, for all $r\in \R$, we have
	\begin{equation*}
		\abs{W_{1,\lambda}^\prime(r)} \le \frac{1}{\lambda}\abs{r} + c
		\qquad\text{and}\qquad
		\abs{W_{1,\lambda}(r)} \le \frac{1}{\lambda}\abs{r}^2 + c\abs{r} \,,
	\end{equation*}
	where $c\coloneqq W_{1,\lambda}^\prime(0)$. 
    
	\item \label{M2}
    For all $r\in\R$, it holds
	\begin{equation*}
		0 \le W_{1,\lambda}(r) \le W_1(r)
		\qquad\text{and}\qquad 
		W_{1,\lambda}(r) \to W_1(r)
		\quad\text{as $\lambda\to 0$.}
	\end{equation*}
    
    \item \label{M3}
    For any $M>0$, there exist $B_M\in\R$ and $\lambda_M>0$ such that for all $\lambda\in (0,\lambda_M)$, it holds
    \begin{equation*}
		W_{1,\lambda}(r) \ge M r^2 + B_M
        \quad\text{for all $r\in\R$.}
	\end{equation*}
    This property follows since $W_1$ is assumed to have super-quadratic growth (cf.~\eqref{COND:COERC} in Assumption~\ref{ASS:POT:SING:1}).
    
\end{enumerate}
For more details, we refer, for example, to \cite[Section~3]{Colli2023} and \cite[Section~5]{Garcke2017}.
Comprehensive overviews of the theory of monotone operators and the Moreau--Yosida approximation can be found, e.g., in \cite[Chapter~2]{Barbu2010}, \cite[Chapter~2]{Brezis} or \cite[Chapter~IV]{Showalter}.

\section{The gradient flow interpretation of our model}\label{sec:modelling}

We start by recalling the energy functional $F$ introduced in \eqref{DEF:F} in the setting of \ref{ASS:DOM} and \ref{ASS:COEFF}: 
\begin{equation}
    F(u,h) = \int_\Omega \left( \frac{1}{\eps} W(u) 
        + \frac{1}{2} \eps |\nabla u|^2 
        + \frac{1}{2} G \nabla h \cdot \nabla h 
        + \frac{1}{2} \kappa |\Delta h|^2 
        - L \nabla h \cdot \nabla u \right) \dx.
\end{equation}
As we are interested in detecting local minimizers of $F$, we consider the gradient flow  
\begin{equation}
    \label{EQ:GFE:1}
    \scp{\delt u}{v}_{-1} + \scp{\delt h}{g}_{L^2} = - \del F(u,h)[v,g]
    \quad\text{on $[0,T]$}
\end{equation}
of the energy functional $F$. Formally speaking, the solution to the gradient flow equation decreases the energy $F$ in the fastest possible way with respect to a certain metric. Therefore, a solution to the gradient flow equation \eqref{EQ:GFE:1} is expected to approach a local minimizer of the functional $F$ as time tends to infinity.
In \eqref{EQ:GFE:1}, the metric which determines the steepest descent direction originates from the inner products $\scp{\cdot}{\cdot}_{-1}$ and $\scp{\cdot}{\cdot}_{L^2}$. For the order parameter $u$, we consider an $H^{-1}$ metric to ensure the conservation of mass in the sense that $\mean{u(t)} = \mean{u(0)}$ for all $t\in [0,T]$. Concerning the height-function $h$, simply choose the $L^2$ metric. This does not enforce the conservation of the mean. Nevertheless, since the energy functional $F$ depends on $h$ only through its gradient and Laplacian, the spatial mean of 
$h$ remains conserved in model~\eqref{PDE}.
For more information about gradient flows in metric spaces, we refer to the book \cite{AGSbook}.

Using integration by parts over $\Omega$ and exploiting periodicity of functions on the domain $\Omega$, the first variation of~$F$ along an arbitrary direction $(v,g) \in H^1(\Omega) \times H^2(\Omega)$ can be expressed as
\begin{equation}
    \label{EQ:firstvar}
    \begin{split}
        \partial F(u,h)[v,g] 
        %\gry = \lim_{s\to 0} \frac{F(u+sv, h+sg) - F(u,h)}{s} 
        &= \int_\Omega \left( \frac{1}{\eps}W'(u)v 
            + \eps \nabla u \cdot \nabla v 
            + G \nabla h \cdot \nabla g 
            + \kappa \Delta h\Delta g 
            - L \nabla u \cdot \nabla g 
            - L \nabla h \cdot \nabla v \right) \dx 
        \\
        %&\gry \overset{(\ast)}{=} \int_\Omega \left( \frac{1}{\eps} W'(u) 
         %       - \eps \Delta u 
          %      - \sum_{i=1}^n \partial_{i}\Big( \sum_{j=1}^n L_{ij} \, \partial_{j}h\Big) \right)v 
        %\\
        %&\quad \gry \
         %   + \left( - \sum_{i=1}^n \partial_{i}\Big( \sum_{j=1}^n G_{ij} \, \partial_{j}h\Big)
          %      + \kappa \Delta^2 h 
          %      - \sum_{i=1}^n \partial_{i} \Big( \sum_{j=1}^n L_{ij} \, \partial_{j}u \Big) \right) g \dx 
        %\\
        &= \int_\Omega \left( \frac{1}{\eps} W'(u) 
                - \eps \Delta u + \Div (L \nabla h) \right)v 
            + \left( - \Div( G \nabla h) 
                + \kappa \Delta^2 h 
                + \Div(L \nabla u) \right) g \dx. 
    \end{split}
\end{equation}

To motivate our model \eqref{PDE}, we show that every sufficiently regular solution of the gradient flow equation is a solution to system \eqref{PDE}. Note that the gradient flow equation implicates the condition $\mean{\delt u}\equiv 0$ on $[0,T]$. 
Integrating by parts, the first summand on the left-hand side of \eqref{EQ:GFE:1} can be reformulated as
\begin{equation}\label{EQ:deltu}
    \scp{\delt u}{v}_{-1} 
    = \int_\Omega \nabla (-\Delta)^{-1} \delt u \cdot \nabla (-\Delta)^{-1} v \dx
    = \int_\Omega (-\Delta)^{-1} \delt u \, v \dx,
    \quad\text{on $[0,T]$.}
\end{equation}
We further introduce the function
\begin{equation}
    \label{PDE:GF:2}
    \mu \coloneqq - \eps \Delta u + \frac{1}{\eps} W'(u) - \Div (L \nabla h)
    \quad\text{in $\Omega_T$.}
\end{equation}
Combining \eqref{EQ:firstvar}--\eqref{PDE:GF:2}, we infer
\begin{align}
    \label{EQ:GFE:2}
    \begin{split}
        &\intO (-\Delta)^{-1} \delt u \, v \dx
            + \intO \delt h \, g \dx
        \\ &\quad
        = - \intO \mu v \dx 
           \, - \intO \left( - \Div( G \nabla h) 
                + \kappa \Delta^2 h 
                + \Div(L \nabla u) \right) g \dx.
    \end{split}
\end{align}
for almost all $t\in(0,T)$ and for all $(v,g) \in H^1(\Omega) \times H^2(\Omega)$. Choosing $g\equiv 0$ in \eqref{EQ:GFE:2} leads to
\begin{equation*}
    (-\Delta)^{-1} \delt u = - \mu \quad\text{in $\Omega_T$.}
\end{equation*}
By the definition of the operator $(-\Delta)^{-1}$ (see~\ref{N:MEAN}), this means that %$\mu$ is a solution of
\begin{equation}
    \label{PDE:GF:1}
    \delt u = \Lap \mu 
    \quad\text{in $\Omega_T$.}
\end{equation}
Moreover, choosing $v\equiv 0$ in \eqref{EQ:GFE:2}, we obtain 
\begin{equation}
    \label{PDE:GF:3}
    \delt h = \Div( G \nabla h) 
                - \kappa \Delta^2 h 
                - \Div(L \nabla u) 
    \quad\text{in $\Omega_T$.}
\end{equation}
Collecting \eqref{PDE:GF:1}, \eqref{PDE:GF:2} and \eqref{PDE:GF:3}, we conclude that the functions $u$, $\mu$ and $h$ satisfy the equations
\begin{subequations}
    \label{PDE*}
    \begin{alignat}{2}
        \label{PDE:1*}
        &\delt u = \Lap \mu
        &&\quad\text{in $\Omega_T$},
        \\
        \label{PDE:2*}
        &\mu = - \eps \Lap u + \eps^{-1} W'(u) + \Div(L\Grad h)
        &&\quad\text{in $\Omega_T$},
        \\
        \label{PDE:3*}
        &\delt h = \Div(G\Grad h) - \kappa \Delta^2 h - \Div(L\Grad u)
        &&\quad\text{in $\Omega_T$},
    \end{alignat}
\end{subequations}
which is exactly the system \eqref{PDE}.

Formally, because of its gradient flow structure, the system \eqref{PDE*} is thermodynamically consistent, meaning that its solutions satisfy the following energy inequality on $[0,T]$:
\begin{align}\label{PDE:EnIneq*}
    \ddt F(u,h) 
    = - \|\partial_t u\|_{-1}^2 - \|\partial_t h\|_{L^2(\Omega)}^2
    = - \|\nabla \mu\|_{L^2(\Omega)}^2 - \|\partial_t h\|_{L^2(\Omega)}^2
    \leq 0.
\end{align}
This relation is obtained by testing \eqref{EQ:GFE:1} with $v=\partial_t u$ and $g=\partial_t h$.
An energy dissipation law similar to \eqref{PDE:EnIneq*} will play a key role in the mathematical analysis presented in this work.

% \subsection{Reduced models (or some other title)}
% \begin{itemize}
% \item Perform a non-dimensionalization in order to remove the physical units from the variables and the parameters, see Table \ref{tab:physical_parameters}.

% \item With the choice $L = \mathrm{diag}\,(\Lambda, \delta_1)$, $G = \mathrm{diag}\,(\sigma, \delta_2)$, where $\delta_1$, $\delta_2$ are small, Patrik solved the ODEs. Then, for the following choice of parameters, one can expect striped patterns:
% \begin{align*}
%     0 < \Lambda < \frac{\sigma^2}{4 \kappa}.
% \end{align*}
% This was verified in our simulations.
% \end{itemize}

% \begin{table}
% \centering
% \begin{tabular}{c||c|c|c|c|c|c|c}
%      Parameter & vesicle size & period & $\Lambda$ & $\sigma$ & $a_4/\varepsilon$ & $\kappa$ & $\varepsilon$
%      \\[1ex] \hline 
%      Physical unit &
%      [m] & [m] & [J/m] & [J/$\mathrm{m}^2$] & [J/$\mathrm{m}^2$] & [J] & [J]
%      \\[1ex] \hline 
%      Estimated values & 
%      $10^{-5}$ m & $10^{-6}$ m & $10^{-12}$ J/m & $10^{-5}$ J/$\mathrm{m}^2$ & $10^{-5}$ J/$\mathrm{m}^2$ & $10^{-19}$ J & $5\cdot 10^{-19}$ J
% \end{tabular}
% \caption{Physical parameters from \cite{Komura:Langmuir:2006,rozovsky2005formation}.}
% \label{tab:physical_parameters}
% \end{table}

\section{Mathematical analysis}
\subsection{Main results} \label{sec:mainresult}
\subsubsection{Existence of weak solutions for regular potentials} 
$\phantom{.}$\\ 
%Also the command ``$\backslash$hfill$\backslash\backslash$'' enforces a line break after the sub(sub)section command.\hfill\\
%
Before stating the result, we first introduce the notion of a weak solution in case $W$ is a regular potential.

\begin{definition}\label{DEF:WS:REG}
Suppose that the assumptions \ref{ASS:DOM}, \ref{ASS:COEFF} and \ref{ASS:POT:REG} hold, and let $u_0\in H^1(\Omega)$ and $h_0\in H^2(\Omega)$
be prescribed initial data. We set $m\coloneqq\mean{u_0}$ and $\rho\coloneqq\mean{h_0}$.
A triplet $(u,\mu,h)$ is called a weak solution to system \eqref{PDE} on $[0,T]$ if the following conditions are fulfilled:
\begin{enumerate}[label=\textnormal{(\roman*)}, topsep=1ex, itemsep=1ex, leftmargin=*]
\item\label{DEF:WS:REG:1} 
The functions $u$, $\mu$ and $h$ have the regularities
\begin{subequations}
\label{REG:WS:REG}
\begin{align}
    u &\in H^1(0,T;H^{-1}(\Omega)) 
        \cap L^\infty(0,T;H^1(\Omega)) 
        \cap C([0,T];L^2_{(m)}(\Omega)), 
    \\%[0.5ex]
    \mu &\in L^2(0,T;H^1(\Omega)), 
    \\%[0.5ex]
    h &\in H^1(0,T;L^2(\Omega)) 
        \cap L^\infty(0,T;H^2(\Omega))
        \cap C([0,T];H^1_{(\rho)}(\Omega)).
\end{align}
\end{subequations}

\item\label{DEF:WS:REG:2} 
The weak formulation
\begin{subequations}
\label{WF:REG}
\begin{align}
    \label{WF:REG:1}
    &\ang{\delt u}{w}_{H^1(\Omega)} 
    = - \intO \Grad \mu \cdot \Grad w \dx,
    \\
    \label{WF:REG:2}
    &\intO \mu\, \vartheta \dx
    = \intO \eps \Grad u \cdot \Grad \vartheta 
    + \frac 1\eps W'(u)\, \vartheta 
    - L \Grad h \cdot \Grad \vartheta \dx,
    \\
    \label{WF:REG:3}
    &   \int_\Omega \delt h\, \eta \dx
    %\ang{\delt h}{\eta}_{H^2(\Omega)}
    + \intO G \Grad h \cdot \Grad \eta \dx
    + \kappa \intO \Lap h \, \Lap \eta \dx
    = \intO L \Grad u \cdot \Grad \eta \dx,
\end{align}
\end{subequations}
holds almost everywhere on $[0,T]$ and for all $w,\vartheta \in H^1(\Omega)$ and $\eta\in H^2(\Omega)$. Moreover, it holds
\begin{equation}
    \label{WS:REG:INI}
    u\vert_{t=0} = u_0, 
    \quad
    h\vert_{t=0} = h_0
    \quad\text{a.e. in $\Omega$.}
\end{equation}

\item\label{DEF:WS:REG:3} The weak energy dissipation law
\begin{align}
    F(u(t),h(t)) 
    + \frac12 \int_0^t \left( \norm{\Grad\mu(s)}_{L^2(\Omega)}^2 
    + \norm{\partial_t h (s)}_{L^2(\Omega)}^2 \right) \ds
    \le F(u_0,h_0)
\end{align}
holds for  all $t\in[0,T]$. 
\end{enumerate}
\end{definition}

\medskip

\noindent The following theorem ensures the existence of a weak solution in the sense of Definition~\ref{DEF:WS:REG}.

\begin{theo}\label{THM:WS:REG}
Suppose that the assumptions \ref{ASS:DOM}, \ref{ASS:COEFF} and \ref{ASS:POT:REG} hold, and let $u_0 \in H^1(\Omega)$ and $h_0 \in L^2(\Omega)$ be prescribed initial data. Then there exists a weak solution $(u,\mu,h)$ of system~\eqref{PDE} on $[0,T]$ in the sense of Definition~\ref{DEF:WS:REG}. Moreover, we have
\begin{align}
    \label{REG:WS:REG:HLD}
    u \in C^{0,\frac 14}([0,T];L^2(\Omega))
    \quad\text{and}\quad
    h \in C^{0,\frac 12}([0,T];L^2(\Omega)).     
\end{align}
\end{theo}

\medskip

The proof of this theorem is carried out in Section \ref{SEC:proof_reg}. There, a weak solution to system \eqref{PDE} is constructed using a \textit{minimizing movement scheme}. This technique has already been employed to construct weak solutions to Cahn--Hilliard models, for example, in \cite{GarKno2020, garcke2023diffuse, Garcke2003}.
For more information on the minimizing movement scheme in more general frameworks, we refer to the book \cite{AGSbook} and the recent review article \cite{mielke2023introduction}.

\medskip

\begin{remark}
    We point out that the prefactor $\frac12$ in the weak energy dissipation law %\ref{DEF:WS:REG:3} 
    in Definition~\ref{DEF:WS:REG}\ref{DEF:WS:REG:3} 
    arises due to our construction of a weak solution via 
    a minimizing movement scheme. There, the quadratic distance term in the functional that is to be minimized includes a factor $\frac{1}{2}$, which persists on the continuous level after passing to the limit in the discretization. In principle, the exact energy dissipation law \eqref{PDE:EnIneq*} can be recovered when enough regularity of the solution is known.
\end{remark}

\medskip

\begin{corrol} \label{COR:HIGHREG}
    Suppose that the assumptions \ref{ASS:DOM}, \ref{ASS:COEFF} and \ref{ASS:POT:REG} hold with $d\in\{2,3\}$ and $p\in [2,4]$. 
    We further assume that $W_1$ and $W_2$ are twice continuously differentiable and that there exists a constant $C_2>0$ such that
    \begin{align*}
        |W_1''(s)| \le C_2(1+|s|^{p-2})
        \quad\text{for all $s\in\R$.}
    \end{align*}
    Moreover, let $u_0 \in H^1(\Omega)$ and $h_0 \in L^2(\Omega)$ be prescribed initial data, and let $(u,\mu,h)$ be a corresponding weak solution in the sense of Definition~\ref{DEF:WS:REG}.
    Then $u$ and $h$ have the additional regularities
    \begin{align}
    \label{REG:WS:REG:HIGH}
    u\in L^2(0,T;H^3(\Omega))
    \quad\text{and}\quad
    h\in L^2(0,T;H^4(\Omega)).
\end{align}
\end{corrol}

\subsubsection{Existence of weak solutions for singular potentials} 
$\phantom{.}$\\ 
As in the case of regular potentials, we first introduce the notion of a weak solution in case $W$ is a singular potential.

\begin{definition}\label{DEF:WS:SING}
Suppose that the assumptions \ref{ASS:DOM}, \ref{ASS:COEFF} and \ref{ASS:POT:SING} hold, and let $u_0\in H^1(\Omega)$ and $h_0\in H^2(\Omega)$ be prescribed initial data with
\begin{equation}
    \label{ASS:INI}
    W_1(u_0) \in L^1(\Omega)
    \qquad\text{and}\qquad
    m\coloneqq\mean{u_0} \in \mathrm{int} D(\mathrm{w}_1).
\end{equation}
We further set $\rho\coloneqq\mean{h_0}$.
A quadruplet $(u,\mu,\xi,h)$ is called a weak solution to system \eqref{PDE} on $[0,T]$ if the following conditions are fulfilled:
\begin{enumerate}[label=\textnormal{(\roman*)}, topsep=1ex, itemsep=1ex, leftmargin=*]
\item\label{DEF:WS:SING:1} 
The functions $u$, $\mu$, $h$ and $\xi$ have the regularities
\begin{subequations}
\label{REG:WS:SING:S}
\begin{align}
    u &\in H^1(0,T;H^{-1}(\Omega)) 
        \cap L^\infty(0,T;H^1(\Omega)) 
        \cap C([0,T];L^2_{(m)}(\Omega)), 
    \\%[0.5ex]
    \mu &\in L^2(0,T;H^1(\Omega)), 
    \quad
    \xi \in L^2(0,T;L^2(\Omega)),
    \\[0.5ex]
    h &\in H^1(0,T;L^2(\Omega)) 
        \cap L^\infty(0,T;H^2(\Omega))
        \cap C([0,T];H^1_{(\rho)}(\Omega)),
\end{align}
\end{subequations}
and it holds         
\begin{equation}
    \label{DIFF:INC:1}
    u(t,x) \in D(\mathrm{w}_1) \quad\text{for almost all $(t,x) \in \Omega_T$.}
\end{equation} 

\item\label{DEF:WS:SING:2} 
The weak formulation
\begin{subequations}
\label{WF:SING}
\begin{align}
    \label{WF:SING:1} 
    &\ang{\delt u}{w}_{H^1(\Omega)} 
    = - \intO \Grad \mu \cdot \Grad w \dx,
    \\
    \label{WF:SING:2}
    &\intO \mu\, \vartheta \dx
    = \intO \eps \Grad u \cdot \Grad \vartheta 
        + \frac 1\eps \xi \, \vartheta 
        + \frac 1\eps W_2'(u)\, \vartheta 
        - L \Grad h \cdot \Grad \vartheta \dx,
    \\
    \label{WF:SING:3}
    &\intO \delt h \, \eta \dx
        + \intO G \Grad h \cdot \Grad \eta \dx
        + \kappa \intO \Lap h \, \Lap \eta \dx
    = \intO L \Grad u \cdot \Grad \eta \dx,
\end{align}
\end{subequations}
holds almost everywhere on $[0,T]$ for all $w,\vartheta \in H^1(\Omega)$ and $\eta\in H^2(\Omega)$, and
\begin{equation}
    \label{DIFF:INC:2}
    \xi(t,x) \in \mathrm{w}_1\big(u(t,x)\big) \quad\text{for almost all $(t,x) \in \Omega_T$.}
\end{equation}
Moreover, it holds
\begin{equation}
    \label{WS:SING:INI}
    u\vert_{t=0} = u_0, 
    \quad
    h\vert_{t=0} = h_0
    \quad\text{a.e. in $\Omega$.}
\end{equation}

\item\label{DEF:WS:SING:3} The weak energy dissipation law
\begin{align}
    F(u(t),h(t)) 
    + \frac12 \int_0^t \left( \norm{\Grad\mu(s)}_{L^2(\Omega)}^2 
    + \norm{\partial_t h (s)}_{L^2(\Omega)}^2 \right) \ds
    \le F(u_0,h_0)
\end{align}
holds for all $t\in[0,T]$.
\end{enumerate}
\end{definition}

\medskip

\noindent The existence of a weak solution in the sense of Definition~\ref{DEF:WS:SING} is ensured by the following theorem. 

\begin{theo}\label{THM:WS:SING}
    Let $\Omega_T$ be as in \ref{ASS:DOM}, $\eps,\kappa, G, L$ as in \ref{ASS:COEFF} and assume that $W$ satisfies \ref{ASS:POT:SING}. Let $u_0 \in H^1(\Omega)$ and $h_0 \in L^2(\Omega)$ such that $u_0$ satisfies \eqref{ASS:INI}.
    Then there exists a weak solution $(u,\mu,h)$ of system~\eqref{PDE} on $[0,T]$ in the sense of Definition~\ref{DEF:WS:SING}. Moreover, we additionally have
    \begin{align}
        \label{REG:WS:SING:U}
        u &\in C^{0,\frac 14}([0,T];L^2(\Omega))
            \cap L^2([0,T];H^2(\Omega))
        \\
        \label{REG:WS:SING:H}
        h &\in C^{0,\frac 12}([0,T];L^2(\Omega))
            \cap L^2([0,T];H^4(\Omega)).     
    \end{align}
\end{theo}

\medskip

\noindent The proof of this theorem can be found in Subsection~\ref{sec:proof_sing}.

\subsection{Uniqueness of the weak solution and continuous dependence on the initial data}

\noindent We point out that any potential $W$, which satisfies \ref{ASS:POT:REG} also fulfills \ref{ASS:POT:SING} with $D(W_1)=\R$ and $D(\mathrm{w}_1) =\R$. Therefore, it suffices to formulate the continuous dependence and uniqueness result simply under the assumption that \ref{ASS:POT:SING} is fulfilled.

\begin{theo}\label{THM:ContDep:REG}
Let $\Omega_T$ be as in \ref{ASS:DOM}, $\eps,\kappa, G, L$ as in \ref{ASS:COEFF} and assume that $W$ satisfies \ref{ASS:POT:SING}. 
Let $u_{0,1},u_{0,2}\in H^1(\Omega)$ and $h_{0,1},h_{0,2}\in H^2(\Omega)$ be prescribed initial data and assume that the mean values coincide, i.e.,
\begin{equation*}
    m \coloneqq \mean{u_{0,1}} = \mean{u_{0,2}}
    \quad\text{and}\quad
    \rho \coloneqq \mean{h_{0,1}} = \mean{h_{0,2}}.
\end{equation*}
We additionally assume that both $u_{0,1}$ and $u_{0,2}$ satisfy \eqref{ASS:INI}.
For $i\in\{1,2\}$, let $(u_i,\mu_i,h_i)$
denote corresponding weak solutions of \eqref{PDE} on $[0,T]$ in the sense of Definition~\ref{DEF:WS:SING}. 
Then the stability estimate
\begin{align}
    \label{THM:ContDep:REG:estimate}
    \begin{split}
    & \| u_2(t) - u_1(t) \|_{-1}^2
    + \| h_2(t) - h_1(t) \|_{L^2(\Omega)}^2 
    + \varepsilon \int_0^{t} \norm{\nabla u_2 - \nabla u_1}_{L^2(\Omega)}^2 \ds
    \\
    &\quad
    + 2G_* \int_0^{t} \norm{\nabla h_2 - \nabla h_1}_{L^2(\Omega)}^2\ds
    + \kappa \int_0^{t} \norm{\Delta h_2 - \Delta h_1}_{L^2(\Omega)}^2 \ds
    \\[1ex]
    &\leq C \big( \| u_{0,2} - u_{0,1} \|_{-1}^2
    + \| h_{0,2} - h_{0,1} \|_{L^2(\Omega)}^2 \big)
    \end{split}
\end{align}
holds for all $t\in[0,T]$ and a positive constant $C$ that depends only on $T$, $\varepsilon$, $\kappa$, ${L^*}$, and $\gamma$ (cf.~\eqref{EXACT:C}).

In particular, if $\mathrm{w}_1 = \partial W_1$ is single-valued and the initial data coincide, i.e., $u_{0,1}=u_{0,2}$ and $h_{0,1}=h_{0,2}$, then \eqref{THM:ContDep:REG:estimate} implies the uniqueness of the corresponding weak solution. 
\end{theo}

The proof of this theorem is carried out in Subsection~\ref{SEC:proof_reg_ContDep}.

\subsection{Proofs}

\subsubsection{Proof of Theorem~\ref{THM:WS:REG}} \label{SEC:proof_reg}
We construct a weak solution to system \eqref{PDE} via a \emph{minimizing movement scheme}. The proof is split into several steps.

\paragraph{Step 0: Formulation of the minimizing movement scheme.} 
    Let $N \in \N$ and let $\tau \coloneqq \frac{T}{N}$ be the time step size. 
    We inductively define a time-discrete approximate solution $\{(u^n, h^n)\}_{n=0,...,N} \subseteq H^1_{(m)}(\Omega) \times H^2(\Omega)$ as follows. First, we set $(u^0,h^0) \coloneqq (u_0, h_0) \in H^1_{(m)}(\Omega) \times H^2_{(\rho)}(\Omega)$.
    Assuming that $(u^n,h^n)$ for some $n\in\{0,...,N-1\}$ is already constructed, we define the auxiliary functional
    \begin{equation}
        \begin{split}
            &J_n : H^1_{(m)}(\Omega) \times H^2(\Omega) \to \R\\
            &J_n(u,h) \coloneqq \frac{1}{2\tau}\|u - u^n\|_{-1}^2 + \frac{1}{2\tau} \|h - h^n \|_{L^2(\Omega)}^2 + F(u,h). 
        \end{split}
    \end{equation}
    The next element of the sequence $(u^{n+1}, h^{n+1})$ is then defined as a minimizer of $J_n$, i.e,
    \begin{equation}\label{EQ:minimizing_scheme_inductive_def}
        J_n(u^{n+1}, h^{n+1}) = \min_{H^1_{(m)}(\Omega) \times H^2(\Omega)} J_n.
    \end{equation}
    To make sense of this definition, we have to ensure that the functional $J_n$ has at least one minimizer. If we can prove that $J_n$ is bounded from below, coercive and lower semicontinuous, the existence of a minimizer follows by the direct method of the calculus of variations.
    
    Let $n\in\{0,...,N-1\}$ be arbitrary. We first show that $J_n$ is bounded from below and (weakly) coercive.
    Recall that, according to \ref{ASS:POT:REG:2}, there exist constants $A>0$ and $B\in\R$ such that
    \begin{equation}
        \label{EST:W:2}
        \abs{W_2(s)} \le As^2 + B \quad\text{for all $s\in\R$.}
    \end{equation}
    As in \ref{ASS:POT:REG:3}, we set
    \begin{align*}
        D \coloneqq 1 + A + \frac{\eps(L^*)^2}{\kappa}.
    \end{align*}
    Then, according to \ref{ASS:POT:REG:3}, there exists a constant $B_D\in\R$ such that 
    \begin{equation}
        \label{EST:W:1}
        W_1(s) \ge D s^2 + B_D
        \quad\text{for all $s\in\R$.}
    \end{equation}
    Now, let $(u,h) \in H^1_{(m)}(\Omega) \times H^2(\Omega)$ be arbitrary. Then, using \eqref{EST:W:1} and \eqref{EST:W:2} as well as Hölder's and Young's inequalities, we derive the estimate
    \begin{align}
    \label{EST:COERC:F}
    \begin{aligned}
        F(u,h)
        &\ge \frac{1}{\eps} \big[(D-A)\norm{u}_{L^2(\Omega)}^2 + (B_D - B)\abs{\Omega}\big]
            + \frac{\eps}{2}
            \norm{\Grad u}_{L^2(\Omega)}^2
        \\
        &\quad
            + \frac{G_*}{2} 
            \norm{\Grad h}_{L^2(\Omega)}^2
            + \frac{\kappa}{2}
            \norm{\Lap h}_{L^2(\Omega)}^2
            - L^* \norm{\Delta h}_{L^2(\Omega)}
            \norm{u}_{L^2(\Omega)}
        \\
        &\ge \frac{1}{\eps} (D-A)\norm{u}_{L^2(\Omega)}^2
            + \frac{\eps}{2}
            \norm{\Grad u}_{L^2(\Omega)}^2
            + \frac{G_*}{2} 
            \norm{\Grad h}_{L^2(\Omega)}^2
            + \frac{\kappa}{4}
            \norm{\Lap h}_{L^2(\Omega)}^2
        \\
        &\quad
            - \frac{(L^*)^2}{\kappa} \norm{u}_{L^2(\Omega)}^2
            + \frac{1}{\eps} (B_D - B)\abs{\Omega}
        \\
        &= \frac{1}{\eps} \Bigg[ D - A - \frac{\eps (L^*)^2}{\kappa} \Bigg]
            \norm{u}_{L^2(\Omega)}^2
            + \frac{\eps}{2}
            \norm{\Grad u}_{L^2(\Omega)}^2
        \\
        &\quad
            + \frac{G_*}{2} 
            \norm{\Grad h}_{L^2(\Omega)}^2
            + \frac{\kappa}{4}
            \norm{\Lap h}_{L^2(\Omega)}^2 
            + \frac{1}{\eps} (B_D - B)\abs{\Omega}
        \\
        &= \frac{1}{\eps} 
            \norm{u}_{L^2(\Omega)}^2
            + \frac{\eps}{2}
            \norm{\Grad u}_{L^2(\Omega)}^2
            + \frac{G_*}{2} 
            \norm{\Grad h}_{L^2(\Omega)}^2
            + \frac{\kappa}{4}
            \norm{\Lap h}_{L^2(\Omega)}^2
        \\
        &\quad 
            + \frac{1}{\eps} (B_D - B)\abs{\Omega}.
    \end{aligned}
    \end{align}
    According to the Poincar\'e--Wirtinger inequality, there exists a constant $c_p>0$ that depends only on $\Omega$ such that
    \begin{align*}
        \|h\|_{L^2(\Omega)} \le c_P  \| \nabla h \|_{L^2(\Omega)} + |\Omega|^{\frac 12 } \mean{h} .
    \end{align*}
    Hence, we find constants $c_1,c_2,c_3>0$ depending only on $\Omega$ and the system parameters, and independent of $N$, $n$ and $\tau$, such that
    \begin{align}
    \label{EST:COERC:F:3}
    \begin{aligned}
        F(u,h)
        &\ge c_1 \|u\|_{H^1(\Omega)}^2 
            + c_2 \|h\|_{H^2(\Omega)}^2
            - c_3\mean{h}^2
            - c_3
    \end{aligned}
    \end{align}
    for all $(u,h) \in H^1_{(m)}(\Omega) \times H^2(\Omega)$. 
    
    Moreover, for any $h\in H^2(\Omega)$, we use Hölder's inequality and the triangle inequality in $L^2(\Omega)$ to obtain 
    \begin{align}
        \mean{h}^2
        \le \frac{1}{|\Omega|} \|h\|_{L^2(\Omega)}^2
        \le \frac{2}{|\Omega|} \|h-h^n\|_{L^2(\Omega)}^2
            + \frac{2}{|\Omega|} \|h^n\|_{L^2(\Omega)}^2\,.
    \end{align}
    Hence, demanding
    \begin{align*}
        \tau \le \frac{|\Omega|}{4c_3} \eqqcolon \tau_* \, ,
    \end{align*}
    the functional $J_n$ can be bounded from below by
    \begin{align}
        \begin{split}
        J_n(u,h)
        &\ge \frac{1}{2\tau} \|h - h^n \|_{L^2(\Omega)}^2 + F(u,h)
        \\
        &\ge \Big( \frac{1}{2\tau} - \frac{2c_3}{|\Omega|} \Big)  
            \|h - h^n \|_{L^2(\Omega)}^2
            + c_1 \|u\|_{H^1(\Omega)}^2 
            + c_2 \|h\|_{H^2(\Omega)}^2
            - \frac{2c_3}{|\Omega|} \|h^n \|_{L^2(\Omega)}^2
            - c_3
        \\
        &\ge c_1 \|u\|_{H^1(\Omega)}^2 
            + c_2 \|h\|_{H^2(\Omega)}^2
            - \frac{2c_3}{|\Omega|} \|h^n \|_{L^2(\Omega)}^2
            - c_3
        \end{split}
    \end{align}
    for all $(u,h) \in H^1_{(m)}(\Omega) \times H^2(\Omega)$. 
    This shows that $J_n$ is bounded from below and coercive. As a consequence of the Banach--Alaoglu theorem, we further infer that $J_n$ is also coercive with respect to the weak topology of $H^1(\Omega) \times H^2(\Omega)$.
    This means that the level sets of $J_n$, which are subsets of $H^1_{(m)}(\Omega) \times H^2(\Omega)$, are precompact in $H^1(\Omega) \times H^2(\Omega)$ with respect to the weak topology.
    
    It remains to show the weak lower semicontinuity of $J_n$. To this end, let $(u,h) \in H^1_{(m)}(\Omega) \times H^2(\Omega)$ be arbitrary and let $\{u_k,h_k\}_{k\in \N}$ be a sequence in $H^1_{(m)}(\Omega) \times H^2(\Omega)$ with
    \begin{align*}
        u_k &\to u \quad\text{weakly in $H^1(\Omega)$,}
        \\
        h_k &\to u \quad\text{weakly in $H^2(\Omega)$}
    \end{align*}
    as $k\to \infty$. Using the compact embedding $H^1(\Omega) \doublehookrightarrow L^2(\Omega)$, we infer that, up to subsequence extraction,
    \begin{align*}
        u_k &\to u \quad\text{strongly in $L^p(\Omega)$ and a.e.~~in $\Omega_T$,}
    \end{align*}
    where $p$ is the exponent from \ref{ASS:POT:REG}.
    Using Fatou's lemma and the continuity of $W$, this implies that
    \begin{equation}
        \label{LSC:W}
        \liminf_{k \to \infty} \int_\Omega W(u_k) \dx 
        \geq \int_\Omega W(u) \dx.
    \end{equation}
    The weak lower semicontinuity of $J_n$ in $H^1(\Omega) \times H^2(\Omega)$ then follows directly from \eqref{LSC:W} and the lower semicontinuity of the involved norms.
    
    Therefore, the direct method of the calculus of variations implies the existence of at least one minimizer $(u^{n+1}, h^{n+1})\in H^1_{(m)}(\Omega) \times H^2(\Omega)$. 
    
\paragraph{Step 1: Discretized weak formulation.}
    Let now $n\in\{0,...,N-1\}$ be arbitrary.
    Since $(u^{n+1}, h^{n+1})$ is a minimizer of $J_n$, it solves the corresponding Euler-Lagrange equation:
    \begin{equation*}
        \delta J_n(u^{n+1}, h^{n+1})[\xi, \eta] = 0 \quad\text{for all $(\xi, \eta) \in H^1_{(0)}(\Omega)\times H^2(\Omega)$}.
    \end{equation*}
    More precisely, for all $\xi \in H^1_{(0)}(\Omega)$ and all $\eta \in H^2(\Omega)$ it holds 
    \begin{equation}\label{EQ:EL_for_discrete}
        \begin{split}
           0&= \delta J_n(u^{n+1}, h^{n+1})[\xi, \eta] 
           \\
           % &= \frac{1}{\tau} \left(u^{n+1} -u^n, \xi \right)_{-1} + \frac{1}{\tau}\left(h^{n+1} - h^n, \eta \right)_{L^2(\Omega)} + \int_\Omega \Big( \frac{1}{\eps}W'(u^{n+1})\xi + \eps \nabla u^{n+1} \cdot \nabla \xi \\
           % &\qquad + G \nabla h^{n+1} \cdot \nabla \eta + \kappa \Delta h^{n+1} \Delta \eta - L \nabla h^{n+1} \cdot \nabla \xi - L \nabla \eta \cdot \nabla u^{n+1}\Big) \dx\\
           &= \frac{1}{\tau} \left(u^{n+1} -u^n, \xi \right)_{-1} + \int_\Omega \left( \frac{1}{\eps}W'(u^{n+1}) \xi + \eps \nabla u^{n+1} \cdot \nabla \xi - L \nabla h^{n+1} \cdot \nabla \xi \right)\dx\\
           &\quad + \frac{1}{\tau}\left(h^{n+1} - h^n, \eta \right)_{L^2(\Omega)} + \int_\Omega \left( G \nabla h^{n+1} \cdot \nabla \eta + \kappa \Delta h^{n+1} \Delta \eta - L \nabla \eta \cdot \nabla u^{n+1}\right) \dx.
        \end{split}
    \end{equation}
    Let $\mathring{\mu}^{n+1} \in H^1_{(0)}(\Omega)$ be defined as 
    \begin{equation*}
        \mathring{\mu}^{n+1} \coloneqq (-\Delta)^{-1} \left( -\frac{u^{n+1} - u^n}{\tau}\right).
    \end{equation*}
    This means that $\mathring{\mu}^{n+1}$ satisfies the weak formulation
    \begin{equation*}
        \int_\Omega \nabla \mathring{\mu}^{n+1} \cdot \nabla v \dx = - \left(\frac{u^{n+1} - u^n}{\tau}, v \right)_{L^2(\Omega)} \quad\text{for all  $v \in H^1(\Omega)$}.
    \end{equation*}
    Choosing $v=(-\Delta)^{-1} \xi$ for any function $\xi \in H^1_{(0)}(\Omega)$, we deduce that
    \begin{equation*}
        \begin{split}
            \int_\Omega \mathring{\mu}^{n+1} \xi \dx
            &= \int_\Omega \mathring{\mu}^{n+1} (-\Delta)(-\Delta)^{-1}\xi \dx = \int_\Omega \nabla \mathring{\mu}^{n+1} \cdot \nabla (-\Delta)^{-1} \xi \dx = - \int_{\Omega} \left(\frac{u^{n+1} - u^n}{\tau} \right)(-\Delta)^{-1} \xi \dx\\
            &= -\int_\Omega \nabla (-\Delta)^{-1} \left(\frac{u^{n+1} - u^n}{\tau} \right) \cdot  \nabla (-\Delta)^{-1} \xi \dx = - \left(\frac{u^{n+1} - u^n}{\tau} , \xi \right)_{-1}.
        \end{split}
    \end{equation*}
    Hence, choosing $\eta = 0$ in \eqref{EQ:EL_for_discrete}, we get
    \begin{equation}
        \label{EQ:MUCIRC}
        \int_\Omega \mathring{\mu}^{n+1} \xi \dx = \int_\Omega \frac{1}{\eps} W'(u^{n+1}) \xi  + \eps \nabla u^{n+1}\cdot \nabla \xi  - L\nabla h^{n+1} \cdot \nabla \xi \dx \quad \text{for all $\xi \in H^1_{(0)}(\Omega)$}.
    \end{equation}
    Now, we want to replace $\mathring{\mu}^{n+1}$ by a function $\mu^{n+1}$, which satisfies \eqref{EQ:MUCIRC} even for all $\xi \in H^1(\Omega).$
    To this end, we make the ansatz
    \begin{equation*}
        \mu^{n+1} \coloneqq \mathring{\mu}^{n+1} + C^{n+1},
    \end{equation*}
    where $C^{n+1}$ is a positive constant that may depend on $u^{n+1}$.
    Let now $\xi \in H^1(\Omega)$ be arbitrary. Then it obviously holds $\mathring\xi\coloneqq \xi - \mean{\xi} \in H^1_{(0)}(\Omega)$. Thus, using $\mathring \xi$ as the test function in \eqref{EQ:MUCIRC}, we deduce that
    \begin{equation*}
        \begin{split}
            &\int_\Omega \mu^{n+1} \xi \dx 
            =
            \int_\Omega \big(\mathring{\mu}^{n+1} + C^{n+1}\big) \big(\xi - \mean{\xi}\big) \dx 
            \;+\; \mean{\xi} \left(\int_\Omega (\mathring{\mu}^{n+1} + C^{n+1} )\dx \right)\\
            &\quad= \int_\Omega \mathring{\mu}^{n+1} \left(\xi - \mean{\xi} \right)\dx 
            \;+\; C^{n+1}\mean{\xi} \\
            &\quad= \int_\Omega \left(\frac{1}{\eps} W'(u^{n+1}) \left(\xi - \mean{\xi}\right) + \eps \nabla u^{n+1} \cdot \nabla \xi - L\nabla h^{n+1} \cdot \nabla \xi\right) \dx 
            \;+\; C^{n+1}\mean{\xi}\\
            &\quad= \int_\Omega \left(\frac{1}{\eps} W'(u^{n+1})\xi + \eps \nabla u^{n+1}\cdot \nabla \xi - L\nabla h^{n+1} \cdot \nabla \xi \right)\dx 
            \;+\; \left(C^{n+1} - \int_\Omega \frac{1}{\eps} W'(u^{n+1}) \dx \right) \mean{\xi}.  
        \end{split}
    \end{equation*}
    Therefore, choosing $C^{n+1} \coloneqq \frac{1}{\eps} \int_\Omega W'(u^{n+1}) \dx$ and recalling the definition of $\mathring{\mu}^{n+1}$, we obtain 
    \begin{equation*}
        \mu^{n+1} = (-\Delta)^{-1} \left(- \frac{u^{n+1} - u^n}{\tau } \right) + \frac{1}{\eps} \int_\Omega W'(u^{n+1})\dx.
    \end{equation*}
    Consequently, $\mu^{n+1}$ satisfies
    \begin{equation*}
        \int_\Omega \mu^{n+1} \xi \dx = \int_\Omega \frac{1}{\eps} W'(u^{n+1})\xi + \eps \nabla u^{n+1} \cdot \nabla \xi - L \nabla h^{n+1} \cdot \nabla \xi \dx \quad\text{for all $\xi \in H^1(\Omega)$}.
    \end{equation*}
    Finally, testing \eqref{EQ:EL_for_discrete} with $\xi = 0$ and an arbitrary function $\eta\in H^2(\Omega)$, we conclude that
    \begin{equation*}
        \frac{1}{\tau}\left(h^{n+1} - h^n, \eta \right)_{L^2(\Omega)} = - \int_\Omega \left(G\nabla h^{n+1} \cdot \nabla \eta + \kappa \Delta h^{n+1} \Delta \eta - L\nabla u^{n+1} \cdot \nabla \eta \right)\dx.
    \end{equation*}
    Thus, we arrive at 
    \begin{subequations}\label{EQ:disc_weak_form}
        \begin{align}
            \label{EQ:disc_weak_form:1}
            &\left(\frac{u^{n+1}-u^n}{\tau}, v \right)_{L^2(\Omega)} 
            = - \int_\Omega \nabla \mu^{n+1} \cdot \nabla v\dx, 
            \\
            \label{EQ:disc_weak_form:2}
            &\int_\Omega \mu^{n+1} \xi \, \dx
                = \int_\Omega \left( \frac{1}{\eps} W'(u^{n+1}) 
                    + \eps \nabla u^{n+1} \cdot \nabla \xi 
                    - L\nabla h^{n+1} \cdot \nabla \xi \right) \dx ,
            \\
            \label{EQ:disc_weak_form:3}
            &\left(\frac{h^{n+1} - h^n}{\tau} , \eta \right)_{L^2(\Omega)}
                + \int_\Omega \left( G \nabla h^{n+1} \cdot \nabla \eta + \kappa \Lap h^{n+1} \, \Lap \eta \right) \, \dx
            = \intO L \Grad u^{n+1} \cdot \Grad \eta \dx,
        \end{align}
    \end{subequations}
    for all $v, \xi \in H^1(\Omega)$ and all $\eta \in H^2(\Omega)$, which is a time-discrete version of the weak formulation \eqref{WF:REG}.

\paragraph{Step 2: Extensions and a priori bounds}
The next step is to extend the time-discrete approximate solutions onto the interval $[0,T]$ and to derive uniform bounds on these extensions.

% \subparagraph{Step 2a: Extensions and a priori bounds}
To this end, let $n\in\{1,...,N\}$ be arbitrary. Since $(u^n, h^n)$ is a minimizer of $J_{n-1}$, we notice that
\begin{equation}\label{EQ:un_hn_minimizers}
    \begin{split}
    &\frac{1}{2\tau} \|u^n - u^{n-1}\|_{-1}^2 + \frac{1}{2\tau} \|h^n - h^{n-1}\|_{L^2(\Omega)}^2 + F(u^n, h^n)  
    \\
    &\quad = J_{n-1}(u^n, h^n)
    \leq J_{n-1}(u^{n-1}, h^{n-1}) 
    = F(u^{n-1}, h^{n-1}).
    \end{split}
\end{equation}
Thus, by finite induction, we infer the a priori estimate
\begin{equation}\label{EQ:apriori_bound}
    \inf_{H^1_{(m)}(\Omega) \times H^2(\Omega)} F \leq F(u^n, h^n) \leq F(u^{n-1}, h^{n-1}) \leq F(u_0,h_0).
\end{equation}
Now, using the time-discrete approximate solution $\{(u^n, \mu^n, h^n)\}_{n=0,...,N}$, we introduce the piecewise constant extension
\begin{equation*}
    \bar{u}_\tau \in L^\infty(0,T; H^1_{(m)}(\Omega)), \qquad \bar{\mu}_\tau \in L^\infty(0,T; H^1(\Omega)), \qquad \bar{h}_\tau \in L^\infty(0,T; H^2(\Omega))
\end{equation*}
by defining
\begin{equation*}
    \left(\Bar{u}_\tau (t) , \Bar{\mu}_\tau (t), \Bar{h}_\tau(t) \right) 
    \coloneqq (u^n, \mu^n, h^n) 
    \qquad \text{if $\; t\in ((n-1)\tau, n\tau ]\;$ with $\;n \in \{1, \dots, N\}\,$.}
\end{equation*}
Moreover, we introduce the piecewise affine extension
\begin{equation*}
    \hat{u}_\tau \in C([0,T]; H^1_{(m)}(\Omega)), \qquad \hat{h}_\tau \in C([0,T]; H^2(\Omega)), %\qquad \hat{\mu}_\tau \in C([0,T]; H^1(\Omega)),
\end{equation*}
% as the s of $\{u^n\}, \{h^n\}$ and $\{\mu^n\}$, i.e., for $n = 1,\dots,N$ and $\theta \in[0,1]$ we set
as follows: For any $t\in (0,T]$, there exist some unique $n\in \{1,...,N\}$ and $\theta\in (0,1]$ such that 
$t= \tau(n-1)(\theta-1) + \tau n\theta$. In this case, we set
\begin{equation*}
\begin{split}
        \hat{u}_\tau (t) &\coloneqq (1-\theta)u^{n-1} + \theta u^n,
        % \\
        % \hat{\mu}_\tau (t) &\coloneqq (1-\theta)\mu^{n-1} + \theta \mu^n,
        \\
        \hat{h}_\tau (t) &\coloneqq (1-\theta) h^{n-1} + \theta h^n.
\end{split}
\end{equation*}
Moreover, we set $(\hat{u}_\tau,\hat{h}_\tau)(0) = (u_0, h_0)$. 

Our next goal is to derive uniform bounds on both extensions, which are independent of $N$ and thus also of the step size $\tau$. Therefore, in the remainder of this proof, the letter $C$ will denote generic positive constants that may depend on the initial data and the system parameters, but not on $N$, $n$ or $\tau$. The exact value of $C$ may vary throughout this proof.

\subparagraph{Step 2a: Uniform bounds on the piecewise constant extension.}

First, we derive uniform bounds on $\|\Bar{u}_\tau \|_{L^\infty (0,T; H^1(\Omega))}$ and $\|\Bar{h}_\tau \|_{L^\infty (0,T; H^2(\Omega))}$. Since
\begin{equation}
    \label{EST:CONTEX:1}
    \|\Bar{u}_\tau \|_{L^\infty (0,T; H^1(\Omega))} 
    = \max_{n \in \{0,\dots, N\}}\|u^n\|_{H^1(\Omega)} 
    \quad\text{ and }\quad 
    \|\Bar{h}_\tau \|_{L^\infty (0,T; H^2(\Omega)} 
    = \max_{n \in \{0,\dots, N\}}\|h^n\|_{H^2(\Omega)},
\end{equation}
it suffices to show that $\{u^n\}_{n \in \{0,\dots, N\}}$ and $\{h^n\}_{n \in \{0,\dots, N\}}$ are uniformly bounded in $H^1(\Omega)$ and $H^2(\Omega)$, respectively. 

To this end, let $n\in\{1,...,N\}$ be arbitrary. 
Testing the discrete weak formulation \eqref{EQ:disc_weak_form:3} with $\eta\equiv 1$, we first notice that 
\begin{align*}
    h^{k} \in H^2_{(\rho)}(\Omega)
    \quad\text{for all $k\in\{0,...,N\}$.}
\end{align*}
Hence, invoking \eqref{EST:COERC:F:3} and the a priori bound \eqref{EQ:apriori_bound}, we find that
\begin{align}
    \label{EST:COERC:F:4}
    \begin{aligned}
        c_1 \|u^n\|_{H^1(\Omega)}^2 
            + c_2 \|h^n\|_{H^2(\Omega)}^2
        &\le F(u^n,h^n)
            + c_3 (1+\rho^2)
        \\
        &\le F(u_0,h_0)
            + c_3 (1+\rho^2).
    \end{aligned}
\end{align}
From the assumptions on the initial data $u_0$ and $h_0$ and on the potential $W$ (see \ref{ASS:POT:REG}), we infer that
\begin{align}
    \label{EST:F:INI}
    F(u_0,h_0) \le C.
\end{align}
Thus, we conclude the uniform bound
\begin{align}
    \label{EST:UNI:2}
    \begin{split}
    \|u^n\|_{H^1(\Omega)}^2
        + \|h^n\|_{H^2(\Omega)}^2
    &\le C
    \quad\text{for all $n\in\{0,...,N\}$.}
    \end{split}
\end{align}
Recalling \eqref{EST:CONTEX:1}, this directly implies that
\begin{align}
    \label{EST:CONSTEX:UNI:UH}
    \|\Bar{u}_\tau \|_{L^\infty (0,T; H^1(\Omega))} 
    \leq C
    \quad\text{and}\quad
    \|\Bar{h}_\tau \|_{L^\infty (0,T; H^2(\Omega)} 
    \leq C.
\end{align}

The next goal is to derive a uniform bound on $\|\bar\mu_{\tau}\|_{H^1(\Omega)}$. 
To this end, let $n\in\{0,...,N-1\}$ be arbitrary.
According to \ref{ASS:POT:REG}, we have 
\begin{equation*}
    |W'(s)| \leq C_1(1 + |s|^{p-1}) + \gamma|s|
    \leq C(1  + |s|^{p-1})
    \quad \text{for all } s\in\mathbb{R},
\end{equation*} 
where $p$ is the exponent from \ref{ASS:POT:REG:1} and $\gamma$ is the Lipschitz constant from \ref{ASS:POT:REG:2}. 
Let now $\xi \in C_c^\infty(\Omega; [0,1])$ be arbitrary. Invoking the continuous Sobolev embedding $H^1(\Omega)\hookrightarrow L^p(\Omega)$ and the uniform estimate \eqref{EST:UNI:2}, we deduce that
\begin{equation}
    \label{EST:DW:UN}
    \begin{split}
        \left| \intO \frac{1}{\eps} W'(u^{n+1}) \xi \dx \right| 
        &\leq \frac{1}{\eps} \intO C\left(1  + |u^{n+1}|^{p-1} 
        \right) |\xi| \dx
        \\
        &\leq C \left(1 + \|u^{n+1}\|_{L^p(\Omega)}^{p-1} \right)
        \\
        &\leq C \left(1 + \|u^{n+1}\|_{H^1(\Omega)}^{p-1} \right)
        \leq C.
    \end{split}
\end{equation}
Hence, testing the discrete weak formulation \eqref{EQ:disc_weak_form:2} with $\xi$, we obtain
\begin{equation*}
    \left| \intO \mu^{n+1} \xi \right| \le \left| \intO \frac{1}{\eps} W'(u^{n+1})\xi \dx\right| + \eps \left| \intO \nabla u^{n+1} \cdot \nabla \xi \dx\right| + \left|\intO L\nabla h^{n+1} \cdot \nabla \xi \dx \right| \leq C\big( 1 + \|\xi\|_{H^1(\Omega)}\big).
\end{equation*}
Based on this estimate, we proceed as in \cite[Proof of Lemma~4.5]{GarKno2020} and use a generalized Poincaré inequality (see, e.g., \cite[Sect.~8.16]{Alt2016} or also \cite[Lemma~A.1]{GarKno2020}) to derive the bound
\begin{equation}
    \label{EST:UNI:POINC}
    \|\mu^{n+1}\|_{L^2(\Omega)} \leq C (\|\nabla \mu^{n+1}\|_{L^2(\Omega)} + 1)
    \quad\text{for all $n\in\{0,...,N-1\}$.}
\end{equation}

Therefore, it remains to establish a uniform bound on $\|\nabla \mu^{n+1}\|_{L^2(\Omega)}$.
Let $t\in (0,T]$ be arbitrary. Then there exists an $n\in\{1,...,N\}$ such that $t\in ( (n-1)\tau, n\tau ]$. Thus, by the definition of the piecewise constant extension, we have
\begin{align*}
    \big(\Bar{u}_\tau(t),\Bar{\mu}_\tau(t),\Bar{h}_\tau(t)\big)
    = (u^n,\mu^n,h^n).
\end{align*}
Then, recalling \eqref{EQ:un_hn_minimizers} and the definitions of $\mu^n$ and of the piecewise affine extension $\hat{h}_\tau$, we obtain
\begin{equation}
    \label{EST:DISC:EN:1}
    \begin{split}
        F(\Bar{u}_\tau(t), \Bar{h}_\tau (t)) 
        &+ \frac{1}{2} \int_{(n-1)\tau}^t \left(\|\nabla \Bar{\mu}_\tau (s)\|_{L^2(\Omega)}^2 + \|\partial_s \hat{h}_\tau(s) \|_{L^2(\Omega)}^2 \right) \ds 
        \\
        &\le F(\Bar{u}_\tau(t), \Bar{h}_\tau (t)) 
            + \frac{1}{2}\int_{(n-1)\tau}^{n\tau} \left( \|\nabla \mu^n\|_{L^2(\Omega)}^2 + \|\partial_s \hat{h}_\tau(s) \|_{L^2(\Omega)}^2 \right) \ds
        \\
        &= F(\Bar{u}_\tau(t), \Bar{h}_\tau (t)) 
            + \frac{1}{2}\int_{(n-1)\tau}^{n\tau} \left( \Big\|\frac{u^{n} - u^{n-1}}{\tau}\Big\|_{-1}^2 + \Big\|\frac{h^n - h^{n-1}}{\tau}\Big\|_{L^2(\Omega)}^2 \right) \ds
        \\
        &= F(u^n, h^n) 
            + \frac{1}{2\tau}\|u^{n} - u^{n-1}\|_{-1}^2 + \frac{1}{2\tau} \|h^{n} - h^{n-1}\|_{L^2(\Omega)}^2 
        \\
        &= J_{n-1}(u^n,h^n) \leq J_{n-1}(u^{n-1}, h^{n-1}) = F(u^{n-1}, h^{n-1}) .
    \end{split}
\end{equation}
In particular, choosing $t=k\tau$ for any $k\in\{1, \dots, N\}$, we infer that
\begin{equation}
    \label{EST:DISC:EN:2}
    F(u^k, h^k) + \frac{1}{2}\int_{(k-1)\tau}^{k\tau} \|\nabla \Bar{\mu}_{\tau}(s)\|_{L^2(\Omega)}^2 \ds 
    + \frac{1}{2} \int_{(k-1)\tau}^{k\tau} \|\partial_t \hat{h}_\tau(s) \|_{L^2(\Omega)}^2 \ds  \leq F(u^{k-1}, h^{k-1}).
\end{equation}
for all $k\in\{1, \dots, N\}$.
Combining \eqref{EST:DISC:EN:1} with \eqref{EST:DISC:EN:2} written for $k=1,...,n-1$, and using the a priori estimate \eqref{EQ:apriori_bound}, we conclude that
\begin{equation}\label{EQ:approx_dissipation}
    F(\bar{u}_\tau(t),\bar{h}_\tau(t)) + \frac{1}{2}\int_0^t \|\nabla \bar{\mu}_{\tau}(s)\|_{L^2(\Omega)}^2\ds + \frac{1}{2}\int_0^t \|\partial_t \hat{h}_\tau (s)\|_{L^2(\Omega)}^2 \ds \leq F(u_0, h_0).
\end{equation}
for all $t\in [0,T]$.
In particular, choosing $t=T$, this entails that
\begin{align}
    \int_0^T \|\nabla \bar{\mu}_{\tau}(s)\|_{L^2(\Omega)}^2\ds
    \le 2 F(u_0, h_0) \le C.
\end{align}
In combination with \eqref{EST:UNI:POINC}, this finally shows that
\begin{equation}
\label{EST:CONSTEX:UNI:MU}
\begin{split}
    \|\Bar{\mu}_\tau\|_{L^2(0,T; H^1(\Omega))}^2 
    = \int_0^T \|\Bar{\mu}_\tau \|_{H^1(\Omega)}^2 \ds 
    \leq C \left( 1 + \int_0^T \|\nabla \Bar{\mu}_\tau\|_{L^2(\Omega)}^2 \ds \right)  
    \leq C.
\end{split}
\end{equation}

\subparagraph{Step 2b: Uniform bounds and Hölder estimates for the piecewise affine extension.}
Recalling the construction of the piecewise affine extension, we immediately obtain the uniform bounds
\begin{align}
    \label{EST:AFFEX:UNI:U}
    \|\hat{u}_\tau \|_{L^\infty (0,T; H^1(\Omega))}
    &\le \|\Bar{u}_\tau \|_{L^\infty (0,T; H^1(\Omega))}
    \le C,
    \\
    \label{EST:AFFEX:UNI:H}
    \|\hat{h}_\tau \|_{L^\infty (0,T; H^2(\Omega))}
    &\le \|\Bar{h}_\tau \|_{L^\infty (0,T; H^2(\Omega))}
    \le C.
\end{align}
Let now $t\in(0,T]$ be arbitrary. Then, there exists $n\in\{0,...,N-1\}$ such that $t\in (n\tau,(n+1)\tau]$.
In view of \eqref{EQ:disc_weak_form:1}, it holds
\begin{equation}
    \label{EQ:AFFEX:1}
    \int_\Omega \partial_t \hat{u}_\tau (t)\,  \zeta \dx 
    = \int_\Omega \frac{u^{n+1} - u^n}{\tau}\, \zeta \dx 
    = - \intO \nabla \mu^{n+1} \cdot \nabla \zeta \dx 
    = - \intO \nabla \Bar{\mu}_\tau(t) \cdot \nabla\zeta \dx
\end{equation}
for all $\zeta \in H^1(\Omega)$. 
In particular, for any $\zeta \in H^1(\Omega)$, it follows that
\begin{equation}
    \big| \langle \partial_t \hat{u}_\tau (t) , \zeta \rangle_{H^1(\Omega)} \big| 
    = \left| \int_\Omega \partial_t \hat{u}_\tau (t)\,  \zeta \dx \right| 
    \le \norm{\Bar{\mu}_\tau(t)}_{H^1(\Omega)}\, \norm{\zeta}_{H^1(\Omega)}
\end{equation}
for all $t\in [0,T]$. Hence, after taking the supremum over all $\zeta \in H^1(\Omega)$ with $\|\zeta\|_{H^1(\Omega)} \le 1$ and integrating over $[0,T]$, we use the uniform bound \eqref{EST:CONSTEX:UNI:MU} to conclude that 
\begin{align}
    \label{EST:AFFEX:UNI:DTU}
    \|\hat{u}_\tau \|_{L^2(0,T;H^{-1}(\Omega))}
    \le \norm{\Bar{\mu}_\tau}_{L^2(0,T;H^1(\Omega))}
    \le C.
\end{align}
Furthermore, it follows directly from \eqref{EQ:approx_dissipation} (written for $t=T$) that
\begin{align}
    \label{EST:AFFEX:UNI:DTH}
    \|\partial_t \hat{h}_\tau \|_{L^2(0,T;L^2(\Omega))}
    \le C.
\end{align}

Now, we fix arbitrary times $0 \leq t_1 < t_2 \leq T$. Choosing $\zeta \coloneqq \hat{u}_\tau(t_2) - \hat{u}_\tau(t_1)$ in \eqref{EQ:AFFEX:1} and using Hölder's inequality, we obtain
\begin{equation*}
\begin{split}
    \|\hat{u}_\tau(t_2) - \hat{u}_\tau(t_1)\|_{L^2(\Omega)}^2 
    &= \intO \left( \hat{u}_\tau(t_2) - \hat{u}_\tau(t_1)\right) \left(\int_{t_1}^{t_2} \partial_\tau \hat{u}_\tau(t)\dt \right) \dx\\ 
    &=\int_{t_1}^{t_2} \intO \partial_t \hat{u}_\tau(t) \left( \hat{u}_\tau(t_2) - \hat{u}_\tau(t_1) \right) \dx \dt\\
    &= - \int_{t_1}^{t_2} \intO \nabla \Bar{\mu}_\tau (t) \cdot \left( \nabla\hat{u}_\tau(t_2) - \nabla \hat{u}_\tau(t_1)\right)\dx\dt\\
    &\leq C \|\Bar{u}_\tau\|_{L^\infty(0,T; H^1(\Omega))} \int_{t_1}^{t_2} \|\nabla \Bar{\mu}_\tau(t)\|_{L^2(\Omega)} \dt\\
    &\leq C\, \|\Bar{u}_\tau\|_{L^\infty(0,T; H^1(\Omega))} \, \|\Bar{\mu}_\tau\|_{L^2(0,T; H^1(\Omega))}\,|t_2 - t_1|^{1/2}.
\end{split}
\end{equation*}
In combination with \eqref{EST:CONSTEX:UNI:UH}, we obtain the Hölder estimate
\begin{equation}
    \label{EST:AFFEX:HLD:U}
    \|\hat{u}_\tau(t_2)- \hat{u}_\tau(t_1)\|_{L^2(\Omega)} \leq C\, |t_2-t_1|^{1/4} 
    \quad\text{for all $0\leq t_1<t_2\leq T$.}
\end{equation}
Similarly, by Hölder's inequality, we have
\begin{equation*}
\begin{split}
    \|\hat{h}_\tau(t_2) - \hat{h}_\tau(t_1)\|_{L^2(\Omega)}^2 
    &=\intO \left| \int_{t_1}^{t_2} \partial_t \hat{h}_\tau(s) \ds \right|^2 \dx 
    \\
    &\leq |t_2-t_1|\, |\Omega|\, \int_0^T \intO |\partial_t \hat{h}_\tau (s)|^2 \dx \ds 
    \\[0.5ex]
    &\leq C\, |t_2-t_1|\, \|\partial_t \hat{h}_\tau\|_{L^2(0,T;L^2(\Omega))}.
\end{split}
\end{equation*}
Hence, using the uniform bound \eqref{EST:AFFEX:UNI:DTH}, we conclude the Hölder estimate
\begin{align}
    \label{EST:AFFEX:HLD:H}
    \|\hat{h}_\tau(t_2) - \hat{h}_\tau(t_1)\|_{L^2(\Omega)} \le C\, |t_2-t_1|^{1/2}.
\end{align}

\paragraph{Step 3: Convergence to a weak solution.} 
So far, we have constructed approximate solutions of piecewise constant functions $\{(\bar{u}_\tau, \bar{\mu}_\tau, \bar{h}_\tau)\}_{\tau \in (0,\tau_*)}$  
as well as piecewise affine functions $\{(\hat{u}_\tau, \hat{\mu}_\tau, \hat{h}_\tau)\} _{\tau \in (0,\tau_*)}$. The goal of this step is to send $N\to\infty$ or, equivalently, $\tau\to 0$.

In view of the uniform bounds \eqref{EST:CONSTEX:UNI:UH} and \eqref{EST:CONSTEX:UNI:MU}, the Banach--Alaoglu theorem implies that there exist functions $u$, $\mu$ and $h$ such that
\begin{alignat}{2}
    \label{CONV:CONSTEX:U}
    \Bar{u}_\tau &\to u 
    &&\quad\text{weakly-star in $L^\infty(0,T;H^1(\Omega))$}, 
    \\
    \label{CONV:CONSTEX:MU}
    \Bar{\mu}_\tau &\to \mu 
    &&\quad\text{weakly in $L^2(0,T;H^1(\Omega))$}, 
    \\
    \label{CONV:CONSTEX:H}
    \Bar{h}_\tau &\to h 
    &&\quad\text{weakly-star in $L^\infty(0,T;H^2(\Omega))$},
\end{alignat}
along a non-relabeled subsequence, as $\tau\to 0$. Furthermore, based on the uniform bounds \eqref{EST:AFFEX:UNI:U}, \eqref{EST:AFFEX:UNI:H}, \eqref{EST:AFFEX:UNI:DTU} and \eqref{EST:AFFEX:UNI:DTH}, the Banach--Alaoglu theorem and the Aubin--Lions--Simon lemma imply that there exist functions $\hat u$, $\hat \mu$ and $\hat h$ such that, up to subsequence extraction, it holds for $\tau\to 0$
\begin{alignat}{2}
    \nonumber
    \label{CONV:AFFEX:U} 
    \hat{u}_\tau &\to \hat u 
    &&\quad\text{weakly-star in $L^\infty(0,T;H^1(\Omega))$, 
        weakly in $H^1(0,T;H^{-1}(\Omega))$}, 
    \\
    &{}&& \qquad\text{strongly in $C(0,T;L^2(\Omega))$},
    \\
    \label{CONV:AFFEX:H}
    \nonumber
    \hat{h}_\tau &\to \hat h 
    &&\quad\text{weakly-star in $L^\infty(0,T;H^2(\Omega))$,
        weakly in $H^1(0,T;L^2(\Omega))$},
    \\
    &{}&& \qquad\text{strongly in $C(0,T;H^1(\Omega))$}.
\end{alignat}

Next, let $t\in (0,T]$ be arbitrary. Then there exists a unique $n\in \{1,...,N\}$ with $t\in ((n-1)\tau, n\tau]$, as well as a unique $\theta\in (0,1]$ such that $t= \tau(n-1)(\theta-1) + \tau n\theta$. Recalling the definitions of the piecewise constant and piecewise affine extensions, we use the Hölder estimate \eqref{EST:AFFEX:HLD:U} to deduce that
\begin{align}
    \label{COMP:EXT}
    \|\hat{u}_\tau(t) - \Bar{u}_\tau(t)\|_{L^2(\Omega)}
    = \|\hat{u}_\tau(t) - \hat{u}_\tau(\tau n)\|_{L^2(\Omega)}
    \le C |t-n\tau|^\frac{1}{4} 
    \le C \tau^\frac{1}{4} .
\end{align}
Similarly, we obtain
\begin{align}
    \|\hat{h}_\tau(t) - \Bar{h}_\tau(t)\|_{L^2(\Omega)}
    \le C \tau^\frac{1}{2} .
\end{align}
Sending $\tau\to 0$, the above convergences, i.e., \eqref{CONV:CONSTEX:U}, \eqref{CONV:CONSTEX:H}, \eqref{CONV:AFFEX:U} and \eqref{CONV:AFFEX:H}, imply that
\begin{align}
    \hat u = u
    \quad\text{and}\quad
    \hat h = h
    \quad\text{a.e.~in $\Omega_T$.}
\end{align}
Moreover, from \eqref{COMP:EXT} we deduce that
\begin{align*}
    \|\Bar{u}_\tau - u\|_{L^2(\Omega_T)}
    &\le \|\hat{u}_\tau - \Bar{u}_\tau\|_{L^2(\Omega_T)}
        + \|\hat{u}_\tau - u\|_{L^2(\Omega_T)}
    \\
    &\le CT\tau^\frac{1}{2} + \|\hat{u}_\tau - u\|_{L^2(\Omega_T)} .
\end{align*}
Consequently, up to subsequence extraction, \eqref{CONV:AFFEX:U} further yields
\begin{align}
    \label{CONV:CONTEX:PTW:U}
    \Bar{u}_\tau \to u
    \quad\text{strongly in $L^2(\Omega_T)$ and a.e.~in $\Omega_T$} ,
\end{align}
as $\tau\to 0$.
Since
\begin{align*}
    \mean{u^n}=m 
    \quad\text{and}\quad
    \mean{h^n}=\rho 
    \quad\text{for all $n\in\{0,\ldots,N\}$,}
\end{align*}
it follows from the definition of the piecewise affine extension that
\begin{align*}
    \mean{\hat{u}_\tau(t)}=m 
    \quad\text{and}\quad
    \mean{\hat{h}_\tau(t)}=\rho 
    \quad\text{for all $t\in [0,T]$.}
\end{align*}
Hence, the convergences \eqref{CONV:AFFEX:U} and \eqref{CONV:AFFEX:H} imply that
\begin{align*}
    u(t) \in H^1_{(m)}(\Omega)
    \quad\text{and}\quad
    h(t) \in H^2_{(\rho)}(\Omega)
    \quad\text{for all $t\in [0,T]$.}
\end{align*}
Altogether, we conclude that the functions $(u,\mu,h)$ have the regularities demanded in Definition~\ref{DEF:WS:REG}\ref{DEF:WS:REG:1}.

Next, we show that the triplet $(u,\mu, h)$ satisfies the weak formulation. 
To this end, let $\varpi \in C_c^\infty([0,T])$ be arbitrary. Multiplying the equations in \eqref{EQ:disc_weak_form} by $\varpi$ and integrating over $[0,T]$, we infer that
\begin{subequations}
    \label{EQ:DW}
    \begin{align}
        \label{EQ:DW:1}
        &\int_0^T \langle \partial_t \hat{u}_\tau , v\rangle_{H^1(\Omega)} \varpi\dt 
        = - \iint_{\Omega_T} \nabla \bar{\mu}_\tau \cdot \nabla v \,\varpi\dx\dt,
        \\
        \label{EQ:DW:2}
        &\iint_{\Omega_T} \bar{\mu}_\tau \xi \,\varpi \dx\dt 
        = \iint_{\Omega_T} \frac{1}{\eps} W'(\bar{u}_\tau) \xi \,\varpi
            + \eps \nabla \bar{u}_\tau \cdot \nabla \xi \,\varpi
            - L \nabla \bar{h}_\tau \cdot \nabla \xi \,\varpi\dx\dt,
        \\
        \label{EQ:DW:3}
        &\int_0^T (\partial_t \hat{h}_\tau, \eta)_{L^2(\Omega)}\,\varpi \dt 
            + \iint_{\Omega_T} G \nabla \bar{h}_\tau \cdot \nabla \eta \,\varpi
                + \kappa \Delta \bar{h}_\tau \Delta \eta \,\varpi\dx \dt 
        = \iint_{\Omega_T} L \nabla \bar{u}_\tau \cdot \nabla \eta \,\varpi\dx\dt ,
    \end{align}
\end{subequations}
for all $v, \xi \in H^1(\Omega)$ and all $\eta \in H^2(\Omega)$.

In order to pass to the limit $\tau \to 0$, we still need to discuss the nonlinear term, which involves $W'$. 
For the exponent $p$ that was introduced in \ref{ASS:POT:REG:1}, we define
\begin{align*}
    q\coloneqq \frac{p}{p-1},
\end{align*}
which satisfies $\tfrac{1}{p} + \tfrac{1}{q} = 1$.
Recalling the growth assumptions on $W_1$ (see~\ref{ASS:POT:REG:1}), we infer that
\begin{align}
    \int_\Omega |W_1'(\bar{u}_\tau)|^{q} \dx
    \le C + C \int_\Omega |\bar{u}_\tau|^{q(p-1)} \dx
    = C + C \int_\Omega |\bar{u}_\tau|^{p} \dx
\end{align}
for all $t\in[0,T]$.
By Sobolev's embedding theorem, we have $H^1(\Omega) \emb L^p(\Omega)$ continuously. Hence, we deduce the uniform bound
\begin{align}
    \|W_1'(\bar{u}_\tau)\|_{L^\infty(0,T;L^q(\Omega))}^q
    \le C + C \|\bar{u}_\tau\|_{L^\infty(0,T;H^1(\Omega))}^{p}
    \le C.
\end{align}
Now, the Banach--Alaoglu theorem implies that there exists a function 
$W^* \in L^\infty(0,T; L^q(\Omega))$ such that, up to subsequence extraction,
\begin{equation*}
    W_1'(\bar{u}_\tau) %\overset{\ast}{\rightharpoonup} 
    \to W^* \qquad \text{weakly-star in } L^\infty(0,T; L^q(\Omega)).
\end{equation*} 
as $\tau\to 0$. However, since $W_1'$ is continuous, the pointwise convergence \eqref{CONV:CONTEX:PTW:U} further entails that
\begin{align*}
    W_1'(\bar{u}_\tau) \to W_1'(u)
    \quad\text{a.e.~in $\Omega_T$}
\end{align*}
However, as a consequence of Egorov's theorem (see \cite[Proposition 9.1c]{dibenedetto2002real}), we know that the weak limit and the pointwise limit coincide. This means that
\begin{equation}
    W^* = W_1'(u) \quad \text{a.e.~in } \Omega_T
\end{equation}
Consequently, we have
\begin{equation}
    \label{CONV:DW1U}
    W_1'(\bar{u}_\tau) %\overset{\ast}{\rightharpoonup} 
    \to W_1'(u) 
    \quad\text{weakly-star in } L^\infty(0,T; L^q(\Omega)).
\end{equation} 
Since $W_2'$ was assumed to be Lipschitz continuous, $W_2'$ has at most linear growth. Hence, proceeding similarly as above, we conclude that
\begin{equation}
    \label{CONV:DW2U}
    W_2'(\bar{u}_\tau) %\overset{\ast}{\rightharpoonup} 
    \to W_2'(u) 
    \quad\text{weakly-star in } L^\infty(0,T; L^q(\Omega)).
\end{equation}

Now, let $v, \xi \in H^1(\Omega)$, $\eta \in H^2(\Omega)$ and $\varpi\in C_c^\infty(\Omega)$ be arbitrary. In particular, by the Sobolev embedding theorem, we have $\xi\varpi \in L^2(0,T;L^p(\Omega))$.
Hence, using the convergences \eqref{CONV:CONSTEX:U}--\eqref{CONV:AFFEX:H}, \eqref{CONV:DW1U} and \eqref{CONV:DW2U}, we can pass to the limit $\tau\to 0$ in \eqref{EQ:DW}. In this way, we obtain
\begin{subequations}
    \label{EQ:DW:LIM}
    \begin{align}
        \label{EQ:DW:LIM:1}
        &\int_0^T \langle \partial_t u , v\rangle_{H^1(\Omega)} \varpi\dt 
        = - \iint_{\Omega_T} \nabla \mu \cdot \nabla v \,\varpi\dx\dt,
        \\
        \label{EQ:DW:LIM:2}
        &\iint_{\Omega_T} \mu \xi \,\varpi \dx\dt 
        = \iint_{\Omega_T} \frac{1}{\eps} W'(u) \xi \,\varpi
            + \eps \nabla u \cdot \nabla \xi \,\varpi
            - L \nabla h \cdot \nabla \xi \,\varpi\dx\dt,
        \\
        \label{EQ:DW:LIM:3}
        &\int_0^T (\partial_t h, \eta)_{L^2(\Omega)}\,\varpi \dt 
            + \iint_{\Omega_T} G \nabla h \cdot \nabla \eta \,\varpi
                + \kappa \Delta h \Delta \eta \,\varpi\dx \dt 
        = \iint_{\Omega_T} L \nabla u \cdot \nabla \eta \,\varpi\dx\dt.
    \end{align}
\end{subequations}
As $\varpi$ is arbitrary, this shows that the triplet $(u,\mu,h)$ satisfies the weak formulation \eqref{WF:REG} for all test functions $v, \xi \in H^1(\Omega)$, $\eta \in H^2(\Omega)$.
Finally, since
\begin{align*}
    \hat{u}_\tau(0) = u_0
    \quad\text{and}\quad
    \hat{h}_\tau(0) = h_0
\end{align*}
the strong convergences in \eqref{CONV:AFFEX:U} and \eqref{CONV:AFFEX:H} ensure that the initial condition \eqref{WS:REG:INI} is fulfilled. This means that all conditions demanded in Definition~\ref{DEF:WS:REG}\ref{DEF:WS:REG:2} are verified.
Moreover, by passing to the limit $\tau\to 0$ in the Hölder estimates \eqref{EST:AFFEX:HLD:U} and \eqref{EST:AFFEX:HLD:H}, we conclude that \eqref{REG:WS:REG:HLD} holds.

\paragraph{Step 4: Energy dissipation law.} 
It remains to verify the weak energy dissipation law in Definition~\ref{DEF:WS:REG}\ref{DEF:WS:REG:3}.
Taking the limes inferior with respect to $\tau\to 0$ in \eqref{EQ:approx_dissipation}, we obtain that
\begin{equation*}
    \begin{split}
         \liminf_{\tau \to 0} \left(F(\bar{u}_\tau(t), \bar{h}_\tau(t)) + \frac{1}{2}\int_0^t \|\nabla \bar{\mu}_\tau(s)\|_{L^2(\Omega)}^2 + \|\partial_t \hat{h}_\tau(s)\|_{L^2(\Omega)}^2 \ds  \right)
         \leq F(u_0,h_0).
    \end{split} 
\end{equation*}
for all $t\in[0,T]$.
Due to the convergences \eqref{CONV:CONSTEX:MU} and \eqref{CONV:AFFEX:H}, and the weak lower semicontinuity of the occurring norms, it follows that
\begin{equation}\label{EQ:diss_ineq_before_limit}
     \liminf_{\tau \to 0} F(\bar{u}_\tau(t), \bar{h}_\tau(t)) 
     + \frac{1}{2}\int_0^t \|\nabla \mu(s)\|_{L^2(\Omega)}^2 
     + \|\partial_t h(s)\|_{L^2(\Omega)}^2 \ds 
    \leq F(u_0,h_0).
\end{equation}
for all $t\in[0,T]$. Recalling the growth assumptions on $W_1$ and $W_2$ (see~\ref{ASS:POT:REG}) with $2 \le p < \infty$ if $d=2$ and $2 \le p \leq \tfrac{2d}{d-2}$ if $d\geq3$, the continuous Sobolev embedding $H^1(\Omega)\emb L^p(\Omega)$, and the uniform estimate \eqref{EST:CONSTEX:UNI:UH}, we deduce that
\begin{align*}
    \intO W(\bar{u}_\tau(t)) \dx 
    \le C + C \intO |\bar{u}_\tau(t)|^p \dx
    \le C + C \|\bar{u}_\tau(t))\|_{H^1(\Omega)}^p
    \le C.
\end{align*}
As $W$ is continuous, the pointwise convergence \eqref{CONV:CONTEX:PTW:U} entails that
\begin{align*}
    W(\bar{u}_\tau(t)) \to W(u(t))
\end{align*}
a.e.~in $\Omega$ for almost all $t\in [0,T]$. Hence, Fatou's lemma implies that
\begin{align}
    \intO W(u(t)) \dx \le \liminf_{\tau \to 0} \intO W(\bar{u}_\tau(t)) \dx
\end{align}
for almost all $t\in [0,T]$. Moreover, using the convergences \eqref{CONV:CONTEX:PTW:U} and \eqref{CONV:CONSTEX:H} along with the weak-strong convergence principle,
\begin{equation}\label{EQ:conv_L_term}
    \lim_{\tau \to 0} \int_\Omega L\nabla \bar{h}_\tau \cdot \nabla \bar{u}_\tau \dx = \int_\Omega L \nabla h \cdot \nabla u \dx 
\end{equation}
for almost all $t\in [0,T]$. All other summands of the energy $F$ can be handled by the convergences \eqref{CONV:CONSTEX:U} and \eqref{CONV:CONSTEX:H} along with the weak lower semicontinuity of the respective norms. In this way, we conclude that
\begin{align}
    F(u(t), h(t)) 
    \le \liminf_{\tau \to 0} F(\bar{u}_\tau(t), \bar{h}_\tau(t))
\end{align}
for almost all $t\in [0,T]$. In view of \eqref{EQ:diss_ineq_before_limit}, this yields
\begin{equation}\label{EQ:dissipation}
    F(u(t), h(t)) + \frac{1}{2}\int_0^t \|\nabla \mu(s)\|_{L^2(\Omega)}^2 + \|\partial_t h(s)\|_{L^2(\Omega)}^2 \ds
    \leq F(u_0, h_0) . 
\end{equation}
for almost all $t \in [0,T]$.

It remains to show that \eqref{EQ:dissipation} actually holds true for all $t\in [0,T]$.
We first notice that the energy inequality is clearly fulfilled at time zero. 
Furthermore, we deduce that there exists a null set $N\subset (0,T]$ such that \eqref{EQ:dissipation} and the estimates
\begin{align}
    \label{REG:WS:REG:t}
    \norm{u(t)}_{H^1(\Omega)}
    \le \norm{u}_{L^\infty(0,T;H^1(\Omega))}
    \quad\text{and}\quad
    \norm{u(t)}_{H^1(\Omega)}
    \le \norm{u}_{L^\infty(0,T;H^1(\Omega))}
\end{align}
hold for all $t\in (0,T]\setminus N$. 
Let now $t\in N$ be arbitrary. Since $N\subset (0,T]$ is a null set, there exists a sequence $(t_k)_{k\in\mathbb N}\subset (0,t)\setminus N$ with $t_k \nearrow t$ as $k\to\infty$.
In view of \eqref{REG:WS:REG:t}, we use the Banach--Alaoglu theorem to infer that there exist functions $u_*\in H^1(\Omega)$ and $h_*\in H^2(\Omega)$ such that
\begin{alignat*}{2}
    u(t_k) &\to u_* &&\quad\text{weakly in $H^1(\Omega)$ and a.e.~in $\Omega$},
    \\
    h(t_k) &\to h_* &&\quad\text{weakly in $H^2(\Omega)$ and a.e.~in $\Omega$},
\end{alignat*}
as $k\to\infty$, up to subsequence extraction. Due to the continuity properties in \eqref{REG:WS:REG}, we know that
\begin{align*}
    u_* = u(t), \quad
    h_* = h(t) \quad\text{a.e.~in $\Omega$.}
\end{align*}
Now, since $W_2$ has at most quadratic growth, we obtain
\begin{align*}
    \intO W_2(u(t)) \dx = \lim_{k\to\infty}\intO W_2(u(t_k)) \dx  
\end{align*}
by means of Lebesgue's general convergence theorem (see, e.g., \cite[Sect.~3.25]{Alt2016}).
Furthermore, since $W_1$ is continuous and convex with polynomial growth of order $p$, we deduce that the functional
\begin{align*}
    L^p(\Omega) \ni v \mapsto \intO W_1(v) \dx \in \R
\end{align*}
is lower semicontinuous. Thus, it follows that
\begin{align*}
    \intO W_1(u(t)) \dx \le \liminf_{k\to\infty}\intO W_1(u(t_k)) \dx
\end{align*}
Moreover, we notice that the integral term in \eqref{EQ:dissipation} is continuous with respect to $t$.
Hence, using the above convergences along with the weak-strong convergence principle and the weak lower semicontinuity of the respective norms, we conclude that
\begin{align*}
    &F(u(t), h(t)) + \frac{1}{2}\int_0^t \|\nabla \mu(s)\|_{L^2(\Omega)}^2 + \|\partial_t h(s)\|_{L^2(\Omega)}^2 \ds
    \\
    &\le \underset{k\to\infty}{\liminf} F(u(t_k),h(t_k))
    + \underset{k\to\infty}{\lim} \int_0^{t_k}\|\mu(s)\|_{L^2(\Omega)}^2 + \|\partial_t h(s)\|_{L^2(\Omega)}^2 \ds
    \\
    &\le \underset{k\to\infty}{\limsup} \Bigg[ F(u(t_k),h(t_k))
    + \int_0^{t_k}\|\mu(s)\|_{L^2(\Omega)}^2 + \|\partial_t h(s)\|_{L^2(\Omega)}^2 \ds \Bigg]
    \\
    &\le F(u_0,h_0)
\end{align*}
This proves that the energy inequality \eqref{EQ:dissipation} is also valid in $t$. As $t\in N$ was arbitrary, we have finally shown that \eqref{EQ:dissipation} holds for all $t\in [0,T]$. Hence, Definition~\ref{DEF:WS:REG}\ref{DEF:WS:REG:3} is verified.

In summary, we have shown that the triplet $(u,\mu,h)$ is a weak solution Definition~\ref{DEF:WS:REG}, which additionally satisfies the regularity property \eqref{REG:WS:REG:HLD}. Thus, the proof is complete.
$\phantom{.}$\hfill$\Box$

\subsubsection{Proof of Corollary~\ref{COR:HIGHREG}} \label{SEC:proof_highreg}
In view of the growth assumptions from \ref{ASS:POT:REG}, we use the continuous Sobolev embedding $H^1(\Omega)\emb L^6(\Omega)$, for $d\in\{2,3\}$, to obtain the bound
\begin{align*}
    \intO |W'(u)|^2 \dx 
    \le C + C\intO |u|^6 \dx
    \le C + C\|u\|_{L^\infty(0,T;H^1(\Omega))}^6
\end{align*}
a.e.~on $[0,T]$.
This shows that $W'(u)\in L^\infty(0,T;L^2(\Omega))$. From the weak formulation \eqref{WF:REG:2} we deduce that for almost all $t\in [0,T]$, $u(t)$ is a weak solution of 
\begin{align}
    \label{EQ:LAPU:STRG}
    -\Lap u(t) 
    = \frac{1}{\eps}\Big( \mu(t) 
        - \frac{1}{\eps}W'\big(u(t)\big) 
        + \Div\big( L \nabla h(t) \big) \Big) .
\end{align}
As the right-hand side belongs to $L^2(0,T;L^2(\Omega))$, we can apply elliptic regularity theory to infer that the order parameter satisfies $u\in L^2(0,T;H^2(\Omega))$.
Due to the Sobolev embedding $H^2(\Omega) \emb L^\infty(\Omega)$ for $d\in\{2,3\}$, this entails that $u\in L^2(0,T;L^\infty(\Omega))$. 

Furthermore, \eqref{WF:REG:3} can be reformulated as 
\begin{align*}
    \intO \Lap h \, \Lap \eta \dx
    = \frac{1}{\kappa} \intO \Big( \delt h 
        - \Div(G \Grad h) 
        - \Div(L \Grad u) \Big) \eta \dx
\end{align*}
a.e.~on $[0,T]$.
This means that for almost all $t\in [0,T]$, $h(t)$ is a weak solution of the inhomogeneous bi-Laplace equation
\begin{align}
    \label{EQ:BLAPH:STRG}
    \Delta^2 h(t)
    = \frac{1}{\kappa}\Big( \delt h(t) 
        - \Div(G \Grad h(t)) 
        - \Div(L \Grad u(t)) \Big)
\end{align}
Since $\partial_t h\in L^2(0,T;L^2(\Omega))$ and $u,h\in L^2(0,T;H^2(\Omega))$, the right-hand side of this equation belongs to $L^2(0,T;L^2(\Omega))$. Hence, using elliptic regularity theory, we deduce that $h\in L^2(0,T;H^4(\Omega))$.

Next, using the growth assumption on $W''$, the chain rule, and Hölder's inequality, we find that
\begin{align*}
    \intO |\nabla \left( W'(u)\right)|^2 \dx 
    &= \intO |W''(u)|^2 \, |\nabla u|^2 \dx
    \le C\intO ( 1 + |u|^4 ) \, |\nabla u|^2 \dx
    \\
    &\le C \big( 1 + \|u\|_{L^\infty(\Omega)}^4 \big) \|u\|_{L^\infty(0,T;H^1(\Omega))}^2.
\end{align*}
This implies $\nabla W'(u) \in L^2(0,T;L^2(\Omega))$, and thus we have $W'(u)\in L^\infty(0,T;H^1(\Omega))$. Together with the regularities $h\in L^2(0,T;H^4(\Omega))$ and $\mu\in L^2(0,T;H^1(\Omega))$, this implies that the right-hand side of \eqref{EQ:LAPU:STRG}
belongs to $L^2(0,T;H^1(\Omega))$. Finally, employing elliptic regularity theory, we conclude that $u\in L^2(0,T;H^3(\Omega))$. Thus all desired regularities are established.
$\phantom{.}$\hfill$\Box$

\subsubsection{Proof of Theorem~\ref{THM:WS:SING}} \label{sec:proof_sing}

The main idea of this proof is to regularize the singular convex part of the potential $W$ using a Moreau--Yosida approximation. Then Theorem~\ref{THM:WS:REG} ensures the existence of a weak solution to the approximated model. The proof is split into several steps.

\paragraph{Step 1: Construction of approximate solutions via Moreau--Yosida regularization.} 

We first approximate the convex part $W_1$ of the double-well potential $W$ by a sequence of regular potentials. For this purpose, we use the \textit{Moreau--Yosida regularization} $W_{1,\lambda}$, as introduced in Section \ref{sec:MorYosReg}. Note that, by choosing 
\begin{align*}
    M = D \coloneqq 1 + A + \frac{\eps(L^*)^2}{\kappa}
\end{align*}
in \ref{M3}, where $A$ is the constant from Assumption~\ref{ASS:POT:SING:2}, we have
\begin{equation}
    \label{EST:W1L:A}
    W_{1,\lambda}(s) \ge D s^2 + B_D
    \quad\text{for all $s\in\R$, $\lambda\in (0,\lambda_D)$.}
\end{equation}
For $\lambda\in (0,\lambda_D)$, we now define the family of regularized potentials
\begin{equation*}
    W_\lambda:\R\to\R, \quad
    W_\lambda(s) = W_{1,\lambda}(s) + W_2(s).
\end{equation*}
From the properties \ref{M1}--\ref{M3} and \ref{ASS:POT:SING:2}, we conclude that for every $\lambda\in (0,\lambda_D)$, $W_\lambda$ is a regular potential, which satisfies Assumption~\ref{ASS:POT:REG} with $p=2$. 
In particular, because of \ref{ASS:POT:SING:2} and \eqref{EST:W1L:A}, $W_\lambda$ is bounded from below by
\begin{equation}
    \label{BND:WL}
    W_\lambda(r) \ge B_D - B
    \quad\text{for all $r\in\R$, $\lambda\in (0,\lambda_D)$.}
\end{equation}
Hence, for every $\lambda\in (0,\lambda_D)$, Theorem~\ref{THM:WS:REG} implies that there exists a weak solution $(u_\lambda,\mu_\lambda,h_\lambda)$ of System~\eqref{PDE} with potential $W_\lambda$ in the sense of Definition~\ref{DEF:WS:REG}. We now intend to show that, as $\lambda\to 0$, the approximate solutions $(u_\lambda,\mu_\lambda,h_\lambda)$ converge to a weak solution of System~\eqref{PDE} with the singular potential $W$ in the sense of Definition~\ref{DEF:WS:SING}.

\paragraph{Step 2: Uniform bounds on the approximate solutions.} 

We now aim to establish bounds on the approximate solutions $(u_\lambda,\mu_\lambda,h_\lambda)$, which are uniform with respect to $\lambda \in (0,\lambda_D)$. 
% To this end, let $M\in (0,D)$ and $\lambda\in (0,\lambda_M)$ be arbitary. 
In the remainder of this proof, the letter $C$ will denote a generic constant that may depend on the initial data and the system parameters, but is independent of $\lambda$. The exact value of $C$ may vary throughout the computations.

As for any $\lambda\in (0,\lambda_D)$, $(u_\lambda,\mu_\lambda,h_\lambda)$ is the weak solution of System~\eqref{PDE} with potential $W_\lambda$ in the sense of Definition~\ref{DEF:WS:REG}, we know that
\begin{align}
    \label{EST:UNI:S1}
    F_\lambda\big(u_\lambda(t),h_\lambda(t)\big) 
    + \frac12 \int_0^t \left( \norm{\Grad\mu_\lambda(s)}_{L^2(\Omega)}^2 
    + \norm{\partial_t h_\lambda (s)}_{L^2(\Omega)}^2 \right) \ds
    \le F_\lambda(u_0,h_0)
\end{align}
for almost all $t\in [0,T]$. In view of \ref{M2}, we have
\begin{equation}
    \label{EST:UNI:S2}
    F_\lambda(u_0,h_0) \le F(u_0,h_0) \le C.
\end{equation}
We point out that $F(u_0,h_0)$ is bounded due to the assumption $W(u_0) \in L^1(\Omega)$ (see~\eqref{ASS:INI}). 
Furthermore, we know that 
\begin{equation}
    \label{EST:UNI:MEAN}
    \mean{u_\lambda(t)} = \mean{u_0} = m
    \quad\text{and}\quad
    \mean{h_\lambda(t)} = \mean{h_0} = \rho
\end{equation}
for almost all $t\in [0,T]$. Consequently, using the Poincar\'e--Wirtinger inequality, we find a constant $c_P>0$ depending only on $\Omega$, $m$ and $\rho$ such that
\begin{align}
    \label{EST:POINC}
    \norm{u_\lambda}_{L^2(\Omega)}^2 
    \le c_P \big(\norm{\Grad u_{\lambda}}_{L^2(\Omega)}^2 + 1 \big)
    \quad\text{and}\quad
    \norm{h_\lambda}_{L^2(\Omega)}^2 
    \le c_P \big(\norm{\Grad h_{\lambda}}_{L^2(\Omega)}^2 + 1 \big)
\end{align}
almost everywhere on $[0,T]$.
Proceeding exactly as in the derivation of \eqref{EST:COERC:F:3}, we find constants $c_1,c_2,c_3>0$ depending only on $\Omega$ and the system parameters, but not on $\lambda$, such that
\begin{align}
\label{EST:UNI:S4}
\begin{aligned}
    F(u_\lambda,h_\lambda)
    &\ge c_1 \|u_\lambda\|_{H^1(\Omega)}^2 
        + c_2 \|h_\lambda\|_{H^2(\Omega)}^2
        - c_3 (1+ \rho^2)
\end{aligned}
\end{align}
for all $\lambda\in (0,\lambda_D)$.
Combining \eqref{EST:UNI:S4} with \eqref{EST:UNI:S1} and \eqref{EST:UNI:S2}, we conclude that
\begin{equation}
    \label{EST:UNI:S5}
    \begin{aligned}
    &\norm{u_\lambda}_{L^\infty(0,T;H^1(\Omega))}
    + \norm{\Grad \mu_\lambda}_{L^2(0,T;L^2(\Omega))}
    \\
    &\quad
    + \norm{h_\lambda}_{H^1(0,T;L^2(\Omega))}
    + \norm{h_\lambda}_{L^\infty(0,T;H^2(\Omega))}
    \le C
    \end{aligned}
\end{equation}
for all $\lambda\in (0,\lambda_D)$. In view of Definition~\ref{DEF:WS:REG}\ref{DEF:WS:REG:2}, we further know that
\begin{subequations}
\label{WF:LAM}
\begin{align}
    \label{WF:LAM:1}
    &\ang{\delt u_\lambda}{w}_{H^1(\Omega)} 
    = - \intO \Grad \mu_\lambda \cdot \Grad w \dx,
    \\
    \label{WF:LAM:2}
    &\intO \mu_\lambda\, \vartheta \dx
    = \intO \eps \Grad u_\lambda \cdot \Grad \vartheta 
        + \frac 1\eps W_\lambda'(u_\lambda)\, \vartheta 
        - L \Grad h_\lambda \cdot \Grad \vartheta \dx,
    \\
    \label{WF:LAM:3}
    &   \int_\Omega \delt h_{\lambda}\, \eta \dx
    %\ang{\delt h}{\eta}_{H^2(\Omega)}
    + \intO G \Grad h_{\lambda} \cdot \Grad \eta \dx
    + \kappa \intO \Lap h_{\lambda} \, \Lap \eta \dx
    = \intO L \Grad u_{\lambda} \cdot \Grad \eta \dx,
\end{align}
\end{subequations}
holds almost everywhere on $[0,T]$ for all $\lambda\in (\lambda,\lambda_D)$ and all test functions $w,\vartheta\in H^1(\Omega)$. 
For any $w\in H^1(\Omega)$, we thus have
\begin{equation*}
    \abs{\ang{\delt u_\lambda}{w}_{H^1(\Omega)}}
    \le \norm{\Grad \mu_\lambda}_{L^2(\Omega)}
        \norm{w}_{H^1(\Omega)}
\end{equation*}
almost everywhere on $[0,T]$. Hence, taking the supremum over all $w\in H^1(\Omega)$ with $\norm{w}_{H^1(\Omega)} \le 1$, we infer
\begin{equation}
    \label{EST:UNI:S6}
    \norm{\delt u_\lambda}_{H^{-1}(\Omega)}
    \le \norm{\Grad \mu_\lambda}_{L^2(\Omega)}
\end{equation}
almost everywhere on $[0,T]$. Next, in view of the properties \ref{M1}--\ref{M3} using the Miranville--Zelik inequality (see \cite[Appendix~A.1]{Miranville2004} or \cite[p.~908]{Gilardi2009}), we find that
\begin{align}
    \label{MZ:1}
    \norm{W'_{1,\lambda}(u_\lambda)}_{L^1(\Omega)}
    \le \intO W'_{1,\lambda}(u_\lambda)\, (u_\lambda - m) \dx .
\end{align}
Choosing $\vartheta = \eps\, (u_\lambda - m)$ in \eqref{WF:LAM:2}, the right-hand side of this estimate can be expressed as
\begin{align}
    \begin{split}
    \label{MZ:2}
    \intO W'_{1,\lambda}(u_\lambda)\, (u_\lambda - m) \dx
    &= - \eps^2 \intO \nabla u_\lambda \cdot \nabla u_\lambda \dx
        + \eps \intO (\mu_\lambda - \mean{\mu_\lambda})\, (u_\lambda - m) \dx
    \\
    &\qquad - %\eps 
    \intO W'_2(u_\lambda) \, (u_\lambda - m) \dx
    + \eps \intO L \nabla h_\lambda \cdot \nabla u_\lambda \dx
    \end{split}
\end{align}
almost everywhere on $[0,T]$.
Since $W'_2$ is Lipschitz continuous (see \ref{ASS:POT:SING:2}), it has at most linear growth.
Combining \eqref{MZ:1} and \eqref{MZ:2}, we use Hölder's inequlity, the Poincar\'e--Wirtinger inequality, and the uniform estimate \eqref{EST:UNI:S5} to deduce that
\begin{align}
    \begin{split}
    \norm{W'_{1,\lambda}(u_\lambda)}_{L^1(\Omega)}
    &\le C \big( 1 + \norm{u_\lambda}_{L^2(\Omega)} \big)
        \big( 1 + \norm{u_\lambda}_{L^2(\Omega)} 
        + \norm{\nabla \mu_\lambda}_{L^2(\Omega)} 
        + \norm{\nabla h_\lambda}_{L^2(\Omega)}\big)
    \\
    &\le C \big( 1 + \norm{\nabla \mu_\lambda}_{L^2(\Omega)} \big)
    \end{split}
\end{align}
almost everywhere on $[0,T]$. Now, choosing $\vartheta\equiv 1$ in \eqref{WF:LAM:2}, we obtain
\begin{align}
    \abs{\mean{\mu_\lambda}}
    \le C \big( 1 + \norm{\nabla \mu_\lambda}_{L^2(\Omega)} \big).
\end{align}
Hence, using the Poincar\'e--Wirtinger inequality and the uniform estimate \eqref{EST:UNI:S5}, we infer that
\begin{align}
    \label{EST:UNI:S7}
    \norm{\mu_\lambda}_{L^2(0,T;H^1(\Omega))}
    \le C \big( 1 + \norm{\nabla \mu_\lambda}_{L^2(0,T;L^2(\Omega))} \big)
    \le C.
\end{align}
As $u_\lambda\in H^1(\Omega)$ and $W'_{1,\lambda}$ is Lipschitz continuous (see \ref{M1}), a chain rule for the combination of Lipschitz and Sobolev functions (see, e.g., \cite[Corollary~3.2]{Ziemer1989}) implies that $W'_{1,\lambda}(u_\lambda) \in H^1(\Omega)$ with
\begin{align*}
    \nabla\big( W'_{1,\lambda}(u_\lambda)\big)
    = W''_{1,\lambda}(u_\lambda) \nabla u_\lambda
    \quad\text{a.e.~in $\Omega_T$.}
\end{align*}
This means that we can choose $\vartheta = \eps\, W'_{1,\lambda}(u_\lambda)$ in \eqref{WF:LAM:2}. In this way, after an integration by parts, we obtain
\begin{align}
    \begin{split}
    &\norm{W'_{1,\lambda}(u_\lambda)}_{L^2(\Omega)}^2
    + \eps^2 \intO W''_{1,\lambda}(u_\lambda)\, \abs{\nabla u_\lambda}^2 \dx
    \\
    &= \intO \big(\eps \mu_\lambda 
        - W'_2(u_\lambda) 
        - \eps \Div( L \nabla h_\lambda) 
    \big) \, W'_{1,\lambda}(u_\lambda) \dx
    \end{split}
\end{align}
almost everywhere on $[0,T]$.
Now, recalling that $W''_{1,\lambda}(u_\lambda) \ge 0$ and that $W'_2$ has at most linear growth, we use Hölder's inequality to deduce that
\begin{align}
    \norm{W'_{1,\lambda}(u_\lambda)}_{L^2(\Omega)}^2
    \le C \big( 1 + \norm{u_\lambda}_{L^2(\Omega)} 
        + \norm{\mu_\lambda}_{L^2(\Omega)} 
        + \norm{h_\lambda}_{H^2(\Omega)}  
    \big) \, \norm{W'_{1,\lambda}(u_\lambda)}_{L^2(\Omega)}
\end{align}
almost everywhere on $[0,T]$.
Thus, in view of \eqref{EST:UNI:S5} and \eqref{EST:UNI:S7}, we conclude that
\begin{align}
    \label{EST:UNI:S8}
    \norm{W'_{1,\lambda}(u_\lambda)}_{L^2(0,T;L^2(\Omega))} \le C.
\end{align}
Finally, combining \eqref{EST:UNI:S5}, \eqref{EST:UNI:S6} and \eqref{EST:UNI:S8}, we obtain the uniform bound
\begin{equation}
    \label{EST:UNI:Sx}
    \begin{aligned}
    & \norm{u_\lambda}_{H^1(0,T;H^{-1}(\Omega))}
    + \norm{u_\lambda}_{L^\infty(0,T;H^1(\Omega))}
    + \norm{\mu_\lambda}_{L^2(0,T;H^1(\Omega))}
    \\
    &\quad
    + \norm{W'_{1,\lambda}(u_\lambda)}_{L^2(0,T;L^2(\Omega))}
    + \norm{h_\lambda}_{H^1(0,T;L^2(\Omega))}
    + \norm{h_\lambda}_{L^\infty(0,T;H^2(\Omega))}
    \le C
    \end{aligned}
\end{equation}
for all $\lambda\in (0,\lambda_D)$.

\paragraph{Step 3: Convergence to a weak solution.} 

In view of \eqref{EST:UNI:Sx}, the Banach--Alaoglu theorem and the Aubin--Lions lemma imply the existence of functions $u$, $\mu$, $\xi$ and $h$, which have the regularity stated in \eqref{REG:WS:SING:S}, such that
\begin{alignat}{2}
    \label{CONV:SING:U}
    u_\lambda &\to u 
    &&\quad\text{weakly-star in $L^\infty(0,T;H^1(\Omega))$, 
    weakly in $H^1(0,T;H^{-1}(\Omega))$,}\nonumber
    \\
    {}&{}
    &&\qquad\text{strongly in $C([0,T];L^2(\Omega))$ and a.e.~in $\Omega_T$,}
    \\
    \label{CONV:SING:MU}
    \mu_\lambda &%\rightharpoonup 
    \to \mu
    &&\quad\text{weakly in $L^2(0,T;H^1(\Omega))$},
    \\
    \label{CONV:SING:W}
    W'_{1,\lambda}(u_\lambda) &%\rightharpoonup 
    \to \xi
    &&\quad\text{weakly in $L^2(0,T;L^2(\Omega))$,}
    \\
    \label{CONV:SING:H}
    h_\lambda &\to h 
    &&\quad\text{weakly-star in $L^\infty(0,T;H^2(\Omega))$, 
    weakly in $H^1(0,T;L^2(\Omega))$,}\nonumber
    \\
    {}&{}
    &&\qquad\text{strongly in $C([0,T];H^1(\Omega))$ and a.e.~in $\Omega_T$,}
\end{alignat}
as $\lambda\to 0$, up to subsequence extraction. Moreover, using \eqref{EST:UNI:MEAN}, we infer that have $\mean{u}=m$ and $\mean{h}=\rho$ almost everywhere in $[0,T]$.
Since $u_\lambda$ and $h_\lambda$ satisfy the initial condition \eqref{WS:REG:INI}, the strong convergences in \eqref{CONV:SING:U} and \eqref{CONV:SING:H} directly imply that $u$ and $h$ satisfy the initial condition \eqref{WS:SING:INI}.
Moreover, using the weak-strong convergence principle, we infer that
\begin{align*}
    \underset{\lambda\to 0}{\lim}\, 
    \iint_{\Omega_T} W'_{1,\lambda}(u_\lambda) \, u_\lambda \dx\dt
    =
    \iint_{\Omega_T} \xi \, u \dx\dt.
\end{align*}
Hence, invoking the theory of maximal monotone operators 
(see \cite[Proposition~1.1, p.~42]{Barbu1976} and also \cite[Section~5.2]{Garcke2017}), we conclude that
\begin{align}
    u_\lambda(t,x) \in D(\mathrm{w}_1)
    \quad\text{and}\quad
    \xi(t,x) \in \mathrm{w}_1\big( u(t,x) \big)  
\end{align}
for almost all $(t,x)\in\Omega_T$. This verifies \eqref{DIFF:INC:1} and \eqref{DIFF:INC:2}.
Due to the convergences \eqref{CONV:SING:U}--\eqref{CONV:SING:H}, it is now straightforward to pass to the limit in the approximate weak formulation \eqref{WF:LAM}. In this way, it follows that the quadruplet $(u,\mu,\xi,h)$ satisfies the weak formulation \eqref{WF:SING}. 
This means that the quadruplet $(u,\mu,\xi,h)$ satisfies the conditions demanded in Definition~\ref{DEF:WS:SING}\ref{DEF:WS:SING:1} and \ref{DEF:WS:SING:2}.

Using the convergences \eqref{CONV:SING:U}--\eqref{CONV:SING:H}, the weak energy dissipation law can be verified by following the line of argument in Step~4 of the proof of Theorem~\ref{THM:WS:REG}. This means that Definition~\ref{DEF:WS:SING}\ref{DEF:WS:SING:3} is fulfilled. Thus, the quadruplet $(u,\mu,\xi,h)$ is a weak solution in the sense of Definition~\ref{DEF:WS:SING}.

\paragraph{Step 4: Additional regularity.}
Let $t_1,t_2 \in [0,T]$ with $t_1<t_2$ be arbitrary. 
Using the Lions--Magenes lemma, we deduce that the mapping 
\begin{align*}
    s\mapsto \frac 12 \|u(s) - u(t_1)\|_{L^2(\Omega)}^2
\end{align*}
is differentiable almost everywhere on $[0,T]$ with
\begin{align*}
    \dds \frac 12 \|u(s) - u(t_1)\|_{L^2(\Omega)}^2
    = \big\langle \partial_t u(s) , u(s) - u(t_1) \big\rangle_{H^1(\Omega)}
\end{align*}
for almost all $s\in [0,T]$. Hence, using the fundamental theorem of calculus, the weak formulation \eqref{WF:SING:1}, and Hölder's inequality, we infer that
\begin{align*}
    \frac 12 \|u(t_2) 
    - u(t_1)\|_{L^2(\Omega)}^2
    &= \int_{t_1}^{t_2} \big\langle \partial_t u(s) , u(s) - u(t_1) \big\rangle_{H^1(\Omega)} \ds
    \\
    &= \int_{t_1}^{t_2} \intO \nabla \mu(s) \cdot \nabla \big(u(s) - u(t_1)\big) \dx\ds
    \\
    &\le 2\, \|u\|_{L^\infty(0,T;H^1(\Omega))} \,\int_{t_1}^{t_2} \|\nabla \mu(s)\|_{L^2(\Omega)} \ds
    \\
    &\le 2\, \|u\|_{L^\infty(0,T;H^1(\Omega))} \, \|\nabla \mu\|_{L^2(0,T;L^2(\Omega))} \,
        |t_2 - t_1|^{\frac 12}.
\end{align*}
This proves that $u\in C^{0,\frac 14}([0,T];L^2(\Omega))$. The regularity $h\in C^{0,\frac 12}([0,T];L^2(\Omega))$ can be established similarly, see also \eqref{EST:AFFEX:HLD:H}.

Arguing similarly to the proof of Corollary~\ref{COR:HIGHREG}, we can interpret $u(t)$, for almost all $t\in [0,T]$, as a weak solution of 
\begin{align}
    \label{EQ:LAPU:STRG:S}
    -\Lap u(t) 
    = \frac{1}{\eps}\Big( \mu(t) 
        - \frac{1}{\eps} \xi
        - \frac{1}{\eps}W_2'\big(u(t)\big) 
        + \Div\big( L \nabla h(t) \big) \Big) .
\end{align}
In view of the regularities \eqref{REG:WS:SING:S}, the right-hand side of this equation belongs to $L^2(0,T;L^2(\Omega))$. Hence, using elliptic regularity theory, we find that $u\in L^2(0,T;H^2(\Omega))$. Based on this regularity, we also obtain $h\in L^2(0,T;H^4(\Omega))$ by proceeding exactly as in Corollary~\ref{COR:HIGHREG}. This completes the proof. 
$\phantom{.}$\hfill$\Box$

\subsubsection{Proof of Theorem~\ref{THM:ContDep:REG}} \label{SEC:proof_reg_ContDep}

We now prove the continuous dependence result. Suppose there exist two weak solution triples $(u_i,\mu_i,\xi_i,h_i)$, $i\in\{1,2\}$, of system \eqref{PDE} to the initial data $u_i(0)=u_{0,i}$ and $h_i(0)=h_{0,i}$ in the sense of Definition~\ref{DEF:WS:SING}.  
We denote their differences by
\begin{align*}
    (u,\mu,\xi,h) \coloneqq (u_2-u_1,\mu_2-\mu_1,\xi_2-\xi_1,h_2-\xi_1), \quad
    (u_0,h_0) \coloneqq (u_{0,2}-u_{0,1}, h_{0,2}-h_{0,1}) \,.
\end{align*}
Then, we observe that
\begin{subequations} \label{eq:system_ContDep}
\begin{align} \label{eq:system_ContDep_u}
    &\ang{\delt u}{w}_{H^1(\Omega)} 
    = - \intO \Grad \mu \cdot \Grad w \dx,
    \\ \label{eq:system_ContDep_mu}
    &\intO \mu\, \vartheta \dx
    = \intO \eps \Grad u \cdot \Grad \vartheta 
    + \frac 1\eps (\xi_1-\xi_2)) \, \vartheta
    + \frac 1\eps (W_2'(u_2)-W_2'(u_1)) \, \vartheta 
    - L \Grad h \cdot \Grad \vartheta \dx,
    \\ \label{eq:system_ContDep_h}
    &   \int_\Omega \delt h\, \eta \dx
    %\ang{\delt h}{\eta}_{H^2(\Omega)}
    + \intO G \Grad h \cdot \Grad \eta \dx
    + \kappa \intO \Lap h \, \Lap \eta \dx
    = \intO L \Grad u \cdot \Grad \eta \dx,
\end{align}
\end{subequations}
for almost every $t\in[0,T]$ and for all $w,\vartheta \in H^1(\Omega)$ and $\eta\in H^2(\Omega)$.

For any $t_0 \in [0,T]$, $\varpi\in C^\infty([0,T])$ and $v\in H^1(\Omega)$, we define
\begin{align*}
    \zeta(t,\cdot) = \begin{cases}
        \int_t^{t_0} v(\cdot)\, \varpi(s) \ds &t \leq t_0,\\
        0 &t>t_0 .
    \end{cases}
\end{align*}
Choosing $w=\zeta(t,\cdot)$ in \eqref{eq:system_ContDep_u}, integrating over $t\in(0,T)$ and by parts over $(0,T)$, we observe
\begin{align*}
    \int_0^T \ang{\partial_t u}{\zeta}_{H^1(\Omega)} \dt
    &= - \int_0^{t_0} \ang{\partial_t \zeta}{u}_{H^1(\Omega)} \dt
    - \int_\Omega u_0 \zeta(0) \dx
    = \int_0^{t_0} \int_\Omega (u-u_0)\, v \dx \, \varpi \dt,
\end{align*}
since $\zeta(t_0)=0$. 
On the other side, on noting \eqref{eq:system_ContDep_u}, we obtain
\begin{align*}
    \int_0^T \ang{\partial_t u}{\zeta}_{H^1(\Omega)} \dt
    = \int_0^T \int_\Omega \nabla\mu \cdot \nabla\zeta \dx \dt
    = \int_0^{t_0} \int_\Omega \nabla \left( \int_0^t \mu(s)\ds \right) \cdot \nabla v \dx\, \varpi \dt.
\end{align*}
Combining both identities, we infer that
\begin{align*}
    \big[(-\Delta_N)^{-1} (u-u_0)\big](t,x) = - \int_0^t \mu(s,x) \ds
    \quad \text{and}\quad
    \big[\partial_t (-\Delta_N)^{-1} (u-u_0)\big](t,x) = - \mu(t,x)
\end{align*}
for almost all $t\in(0,t_0)$ and $x\in\Omega$.
Hence, integrating \eqref{eq:system_ContDep_u} over $(0,t_0)$, using the fundamental theorem of calculus, and choosing $w= -\mu$, we deduce that
\begin{align}
    \label{ID:CONTDEP:1}
    \begin{split}
    - \int_0^{t_0} \int_\Omega u \mu \dx\dt 
    &= \int_0^{t_0} \int_\Omega \nabla (-\Delta_N)^{-1} (u-u_0) \cdot \nabla \partial_t (-\Delta_N)^{-1} (u-u_0) \dx\dt
    \\
    &= \frac12 \int_0^{t_0} \ddt \| u(t)-u_0 \|_{-1}^2 \dt
    = \frac12 \| u(t_0)-u_0 \|_{-1}^2 \,.
    \end{split}
\end{align}
Noting $\mean{u_0}=\mean{u(t_0)}=0$, the inverse triangle inequality yields
\begin{align}\label{eq:ContDep_reg_1}
    \frac12 \| u_0 \|_{-1}^2 
    \geq \int_0^{t_0} \int_\Omega u \mu \dx\dt
    + \frac12 \| u(t_0) \|_{-1}^2 .
\end{align}

Next, we choose $\vartheta(t,\cdot) = u(t,\cdot) \chi|_{[0,t_0]}(t)$ in \eqref{eq:system_ContDep_mu} and integrate over $t\in(0,T)$. After re-sorting the terms, we obtain
\begin{align*}
    &- \int_0^{t_0} \int_\Omega \mu u \dx\dt + \varepsilon \int_0^{t_0} \norm{\nabla u}_{L^2(\Omega)}^2 \dt
    \\
    &\quad
    = - \int_0^{t_0} \int_\Omega \frac{1}{\varepsilon} \big( \xi_2 - \xi_1 \big) u \dx \dt
    - \int_0^{t_0} \int_\Omega \frac{1}{\varepsilon}\big( W'_2(u_2)- W'_2(u_1)\big) u \dx \dt
    + \int_0^{t_0} \int_\Omega L\nabla h \cdot \nabla u \dx\dt.
\end{align*}
Now, we recall that 
\begin{equation}
    \label{DIFF:INC:i}
    \xi_i(t,x) \in \mathrm{w}_1\big(u_i(t,x)\big) \quad\text{for almost all $(t,x) \in \Omega_T$ and $i=1,2$.}
\end{equation}
and that $\mathrm{w}_1$ is a maximal monotone graph. Since $u=u_2-u_1$, this implies that
\begin{align}
    \int_0^{t_0} \int_\Omega \frac{1}{\varepsilon} \big( \xi_2 - \xi_1 \big) u \dx \dt
    \ge 0.
\end{align}
We further recall that, according to \ref{ASS:POT:REG:2}, $W_2'$ is Lipschitz continuous with Lipschitz constant $\gamma$. Therefore, we infer that 
\begin{align}
    -\int_0^{t_0} \int_\Omega \mu u \dx\dt
        + \varepsilon \int_0^{t_0} \norm{\nabla u}_{L^2(\Omega)}^2 \dt
    \leq \frac\gamma\varepsilon \int_0^{t_0} \| u \|_{L^2(\Omega)}^2 \dt
        + \int_0^{t_0} \int_\Omega  L\nabla h \cdot \nabla u \dx\dt.
\end{align}
After using \eqref{ID:CONTDEP:1} to replace the first term on the left-hand side, we use the estimate \eqref{eq:ContDep_reg_1} to conclude that
\begin{align}\label{eq:ContDep_reg_2}
    \frac12 \| u(t_0)\|_{-1}^2
        + \varepsilon \int_0^{t_0} \norm{\nabla u}_{L^2(\Omega)}^2 \dt
    \leq \frac12 \| u_0 \|_{-1}^2
        + \frac\gamma\varepsilon \int_0^{t_0} \| u \|_{L^2(\Omega)}^2 \dt
        + \int_0^{t_0} \int_\Omega  L\nabla h \cdot \nabla u \dx\dt.
\end{align}

Next, we choose $\eta(t,\cdot) = h(t,\cdot) \chi|_{[0,t_0]}(t)$ in \eqref{eq:system_ContDep_h} and integrate over $(0,t_0)$. Making use of the identity
\begin{align*}
    \int_0^{t_0} \int_\Omega \delt h\, h \dx\dt
    = \frac12 \int_0^{t_0} \ddt \|h(t)\|_{L^2(\Omega)}^2
    = \frac12 \|h(t_0)\|_{L^2(\Omega)}^2 - \frac12 \|h_0\|_{L^2(\Omega)}^2,
\end{align*}
and rearranging the terms in the resulting equation, we obtain
\begin{align}\label{eq:ContDep_reg_3}
    \begin{split}
    &\frac12 \|h(t_0)\|_{L^2(\Omega)}^2  
        + \int_0^{t_0} \int_\Omega G\nabla h\cdot\nabla h \dx\dt
        + \int_0^{t_0} \int_\Omega \kappa |\Delta h|^2 \dx\dt
    \\
    &\quad= \frac12 \|h_0\|_{L^2(\Omega)}^2
        + \int_0^{t_0} \int_\Omega L\nabla u \cdot \nabla h \dx\dt.
    \end{split}
\end{align}

Combining \eqref{eq:ContDep_reg_1}, \eqref{eq:ContDep_reg_2}, and \eqref{eq:ContDep_reg_3}, and recalling the assumptions on $G$ and $L$ (see~\ref{ASS:COEFF}), we deduce that
\begin{align}
    \label{EST:CONTDEP:1}
    \begin{split}
    & \frac12\| u(t_0) \|_{-1}^2
    + \frac12\|h(t_0)\|_{L^2(\Omega)}^2 
    + \varepsilon \int_0^{t_0} \norm{\nabla u}_{L^2(\Omega)}^2 \dt
    \\
    &\quad
    + G_* \int_0^{t_0} \norm{\nabla h}_{L^2(\Omega)}^2\dt
    + \kappa \int_0^{t_0} \norm{\Delta h}_{L^2(\Omega)}^2 \dt
    \\
    &\leq \frac12\| u_0 \|_{-1}^2
    + \frac12\|h_0\|_{L^2(\Omega)}^2  
    + \frac\gamma\varepsilon \int_0^{t_0} \| u \|_{L^2(\Omega)}^2 \dt
    + 2L^* \int_0^{t_0} \intO |\nabla u| \, |\nabla h| \dx\dt .
    \end{split}
\end{align}
Recalling the definition of the norm $\|\cdot\|_{-1}$, we use Young's inequality to derive the estimate 
\begin{align*}
    \frac\gamma\varepsilon \|u\|_{L^2(\Omega)}^2 
    &= \frac\gamma\varepsilon \intO u\; (-\Delta) (-\Delta)^{-1} u \dx
    = \frac\gamma\varepsilon \intO \nabla u\cdot \nabla (-\Delta)^{-1} u \dx
    \\
    &\leq \frac{\varepsilon}{4} \|\nabla u\|_{L^2(\Omega)}^2
    + \frac{\gamma^2}{\varepsilon^3} \|u\|_{-1}^2.
\end{align*}
a.e.~on $[0,T]$.
Proceeding in a similarl manner, we obtain
\begin{align*}
    \frac{4{L^*}^2}{\eps} \|\nabla h\|_{L^2(\Omega)}^2
    \le \frac{\kappa}{2} \|\Delta h\|_{L^2(\Omega)}^2
        +  \frac{8 {L^*}^4}{\kappa \eps^2} \|h\|_{L^2(\Omega)}^2.
\end{align*}
a.e.~on $[0,T]$.
By means of Young's inequality, we infer that
\begin{align*}
    2L^* \intO |\nabla u| \, |\nabla h| \dx
    &\leq
    \frac{\eps}{4} \|\nabla u\|_{L^2(\Omega)}^2 + \frac{4{L^*}^2}{\eps} \|\nabla h\|_{L^2(\Omega)}^2
    \\
    &\leq \frac{\varepsilon}{4} \|\nabla u\|_{L^2(\Omega)}^2
        + \kappa \|\Delta h\|_{L^2(\Omega)}^2
        + \frac{\kappa}{2} \|\Delta h\|_{L^2(\Omega)}^2
        + \frac{8 {L^*}^4}{\kappa \eps^2} \|h\|_{L^2(\Omega)}^2.
\end{align*}
a.e.~on $[0,T]$.
Using these inequalities to bound the right-hand side of \eqref{EST:CONTDEP:1} and multiplying both sides by $2$, we obtain
\begin{align}
    \label{EST:CONTDEP:2}
    \begin{split}
    & \| u(t_0) \|_{-1}^2
    + \|h(t_0)\|_{L^2(\Omega)}^2 
    + \eps \int_0^{t_0} \norm{\nabla u}_{L^2(\Omega)}^2 \dt
    \\
    &\quad
    + 2G_* \int_0^{t_0} \norm{\nabla h}_{L^2(\Omega)}^2\dt
    +\kappa \int_0^{t_0} \norm{\Delta h}_{L^2(\Omega)}^2 \dt
    \\
    &\leq \| u_0 \|_{-1}^2 + \|h_0\|_{L^2(\Omega)}^2 
        + \Bigg( \frac{2\gamma^2}{\eps^3} + \frac{16 {L^*}^4}{\kappa \eps^2} \Bigg)
    \Bigg( \int_0^{t_0} \| u \|_{-1}^2 \dt
        + \int_0^{t_0} \|h\|_{L^2(\Omega)}^2 \dt \Bigg).
    \end{split}
\end{align}
Finally, applying a variant of Gronwall's lemma (see, e.g.,~\cite[Lemma~3.1]{garcke_lam_2017}), we conclude that
\begin{align}
    \label{EST:CONTDEP:FINAL}
    \begin{split}
    & \| u(t) \|_{-1}^2
    + \|h(t )\|_{L^2(\Omega)}^2 
    + \eps \int_0^{t} \norm{\nabla u}_{L^2(\Omega)}^2 \ds
    \\
    &\quad
    + 2G_* \int_0^{t} \norm{\nabla h}_{L^2(\Omega)}^2\ds
    +\kappa \int_0^{t} \norm{\Delta h}_{L^2(\Omega)}^2 \ds
    \\
    &\leq C \big( \| u_0 \|_{-1}^2
    + \|h_0\|_{L^2(\Omega)}^2 \big).
    \end{split}
\end{align}
for almost all $t\in(0,T)$, since $t_0$ was arbitrary. Here, the constant $C>0$ is explicitly given by
\begin{align}
    \label{EXACT:C}
    C = \exp\Bigg( \frac{2T\gamma^2}{\eps^3} + \frac{16T {L^*}^4}{\kappa \eps^2} \Bigg).
\end{align}
This verifies the desired estimate \eqref{THM:ContDep:REG:estimate}. 

To obtain the uniqueness statement, assume that $u_{0,1}=u_{0,2}$ and $h_{0,1}=h_{0,2}$ almost everywhere in $\Omega$ and that $\mathrm{w}_1$ is single-valued. 
Then \eqref{THM:ContDep:REG:estimate} directly implies $u_1=u_2$ and $h_1=h_2$ almost everywhere in $\Omega_T$.
As $\mathrm{w}_1$ is single-valued, we use \eqref{DIFF:INC:i} to infer that $\xi_1 = \xi_2$ almost everywhere in $\Omega_T$.
This implies that the right-hand side of \eqref{eq:system_ContDep_mu} is zero, and we finally conclude that $\mu_1=\mu_2$ 
almost everywhere in $\Omega_T$.
This proves the uniqueness of the weak solution and thus, the proof of Theorem~\ref{THM:ContDep:REG} is complete.
$\phantom{.}$\hfill$\Box$

\section{Numerical results}\label{sec:num}

In this section, we present numerical results for \eqref{PDE}. We begin by introducing a finite element discretization of \eqref{PDE} and discuss its practical implementation. Next, we examine key properties of the numerical scheme, including unique solvability, mass conservation, and a discrete energy inequality. Finally, we present numerical experiments, compare our parameter choices with values from the literature, and analyze the influence of different parameters on the model behavior.

\subsection{Finite element discretization}
In the following, we restrict to two space dimensions. We represent the flat torus $\Omega=\mathbb{T}^2 \subset\mathbb{R}^2$ by the the unit square $(0,1)^2$ with periodic boundary identification. 
We divide the time interval $[0,T)$ into equidistant subintervals $[n\tau,(n+1)\tau)$ with $\tau = \frac{T}{N}$, $N\in\mathbb{N}$ and $n \in \{0, \ldots, N-1\}$.
We will throughout use the index $k\in\mathbb{N}$ to indicate the spatial level of refinement. In fact, we require $\{\mathcal{T}_k\}_{k\in\mathbb{N}}$ to be a family of conforming partitionings of $\Omega$ into disjoint open simplices such that $\overline\Omega = \bigcup_{\mathcal{K}\in\mathcal{T}_k} \overline{\mathcal{K}}$. 
In addition, we assume that the family of triangulations is uniformly shape regular (see~\cite[Definition~3.8]{bartels_2016}).
% \footnote{\DT{Yes, this is meant by shape regular or uniformly shape regular or non-degenerate. There exists also the term `quasi-uniform' (see~\cite[Definition~3.12]{bartels_2016}), where in addition $h_K \geq C h$ for all elements $K$ has to hold (e.g., for inverse inequalities), but it is not needed here.}}
We introduce the following notation for the finite element spaces of continuous and piecewise linear functions
\begin{align*}
     \mathcal{S}_k &\coloneqq \left\{ q \in C\big(\overline\Omega\big) \;\middle|\; q|_{\mathcal{K}} \in \mathcal{P}_1(\mathcal{K}) \; \text{for all}\;  \mathcal{K}\in\mathcal{T}_k \right\} 
     \ \subset \ W^{1,\infty}(\Omega) \, ,
\end{align*}
where $\mathcal{P}_1(\mathcal{K})$ denotes the set of affine linear polynomials defined on $\mathcal{K}\in\mathcal{T}_k$. 
% We recall that $\Omega=\mathbb{T}^d$, $d\in\mathbb{R}$ is the flat torus, represented as the unit cube $(0,1)^d$ with periodic boundary identification. 

We assume that the potential $W$ satisfies \ref{ASS:POT:REG}.
This means that $W = W_1 + W_2$, where $W_1$ is convex. In addition, we assume that $W_2$ is concave.
To devise a numerical approximation of the weak formulation from Definition~\ref{DEF:WS:REG} based on the finite dimensional subspace $\mathcal{S}_k$, we need to get rid of the Laplace operators in the weak formulation
\eqref{WF:REG:3}.
To this end, we introduce the auxiliary variable $g\coloneqq -\Delta h$ and reformulate \eqref{WF:REG:3} as
\begin{align}
    \label{WF:REG:3:ALT}
       \int_\Omega \delt h\, \eta \dx
    %\ang{\delt h}{\eta}_{H^2(\Omega)}
    + \intO G \Grad h \cdot \Grad \eta \dx
    + \kappa \intO \nabla g \cdot \nabla \eta \dx
    &= \intO L \Grad u \cdot \Grad \eta \dx
    \\
    \label{WF:REG:4:ALT}
    \intO \nabla h \cdot \nabla q \dx 
    &= \intO g \, q \dx 
\end{align}
for all $\eta, q\in H^1(\Omega)$. 

Now, suppose that discrete initial data $u^0_k, h^0_k \in \mathcal{S}_k$ are given.
Then, for any $n\in\{0,\ldots,N-1\}$, the goal is to find a solution quadruplet 
\begin{align*}
    (u^{n+1}_k, \mu^{n+1}_k, h^{n+1}_k, g^{n+1}_k) \in (\mathcal{S}_k)^4
\end{align*}
such that the discretized weak formulation 
\begin{subequations}
\label{NUM*}
\begin{alignat}{2}
    \label{NUM:1*}
    0 &= \frac{1}{\tau} 
    \skp{u^{n+1}_k - u^n_k}{w_k}_{L^2(\Omega)} 
    + \skp{\nabla \mu^{n+1}_k}{\nabla w_k}_{L^2(\Omega)} \,,
    \\
    \label{NUM:2*} \nonumber
    0 &= - \skp{\mu^{n+1}_k}{\vartheta_k}_{L^2(\Omega)}
    + \eps \skp{\nabla u^{n+1}_k}{\nabla \vartheta_k}_{L^2(\Omega)}
    + \eps^{-1} \skp{W_1'(u^{n+1}_k) + W_2'(u^n_k)}{\vartheta_k}_{L^2(\Omega)}
    \\
    &\quad
    - \skp{L \nabla h^{n+1}_k}{\nabla \vartheta_k}_{L^2(\Omega)} \,,
    \\
    \label{NUM:3*}\nonumber
    0 &= \frac{1}{\tau} 
    \skp{h^{n+1}_k - h^n_k}{\eta_k}_{L^2(\Omega)} 
    + \skp{G \nabla h^{n+1}_k}{\nabla \eta_k}_{L^2(\Omega)}
    + \kappa \skp{\nabla g^{n+1}_k}{\nabla \eta_k}_{L^2(\Omega)}
    \\
    &\quad
    - \skp{L \nabla u^n_k}{\nabla\eta_k}_{L^2(\Omega)} \,,
    \\
    \label{NUM:4*}
    0 &= - \skp{g^{n+1}_k}{q_k}_{L^2(\Omega)} 
    + \skp{\nabla h^{n+1}_k}{\nabla q_k}_{L^2(\Omega)} \,.
\end{alignat}
\end{subequations}
holds for all test functions $(w_k, \vartheta_k, \eta_k, q_k) \in (\mathcal{S}_k)^4$.

Note that \eqref{NUM*} can be solved in two steps. First, the variables $h^{n+1}_k$ and $g^{n+1}_k$ are computed from the linear subsystem \eqref{NUM:3*}--\eqref{NUM:4*}. Afterwards, the nonlinear Cahn--Hilliard subsystem \eqref{NUM:1*}--\eqref{NUM:2*} is solved using Newton's method.
In practice, both the linear subsystem \eqref{NUM:3*}--\eqref{NUM:4*} and the linear systems arising in the Newton iterations for \eqref{NUM:1*}--\eqref{NUM:2*} are solved using well-established preconditioned iterative methods, see, e.g., \cite{elman_silvester_wathen_2014}. In fact, a preconditioned \texttt{MINRES} solver is employed for \eqref{NUM:3*}--\eqref{NUM:4*}, which forms a regular and symmetric saddle point system for $h^{n+1}_k$ and $g^{n+1}_k$ (when \eqref{NUM:4*} is multiplied by the constant $\kappa>0$). For the Cahn--Hilliard subsystem \eqref{NUM:1*}--\eqref{NUM:2*}, strategies from \cite{bosch_2016} can be applied. The overall system \eqref{NUM*} is implemented using the finite element library \texttt{FEniCS} \cite{fenics_book_2012}, which provides access to the linear algebra backend \texttt{PETSc} \cite{petsc-user-ref_2021}. The implementation builds upon the code developed in \cite{trautwein_2021}. 

We will present several numerical results in Section~\ref{sec:numerics}. Before that, we discuss several important properties of the discrete system \eqref{NUM*}.

\subsection{Well-posedness and properties}

In the following, we study properties of the discrete system \eqref{NUM*}. As a first step, we prove that \eqref{NUM*} admits a unique solution for any time step size $\tau > 0$, without requiring a smallness condition.

\begin{lemma}
Let $n\in\{0,\ldots,N-1\}$ be arbitrary and let $(u_k^n,h_k^n)\in(\mathcal{S}_k)^2$ be prescribed.
Then there exists a unique solution $(u^{n+1}_k,\mu^{n+1}_k,g^{n+1}_k,g^{n+1}_k)\in (\mathcal{S}_k)^4$ to the discrete system \eqref{NUM*}.
% For every $n\in\{0,\ldots,N-1\}$, there exists a unique solution to the discrete system \eqref{NUM*}.
\end{lemma}

\begin{proof}
To prove the existence of a unique solution to the subsystem \eqref{NUM:3*}--\eqref{NUM:4*}, it suffices to show that the homogeneous problem
\begin{align*}
    0 &= \frac{1}{\tau} 
    \skp{h^{n+1}_k}{\eta_k}_{L^2(\Omega)} 
    + \skp{G \nabla h^{n+1}_k}{\nabla \eta_k}_{L^2(\Omega)} 
    + \kappa \skp{\nabla g^{n+1}_k}{\nabla \eta_k}_{L^2(\Omega)} \,,
    \\
    0 &= - \skp{g^{n+1}_k}{q_k}_{L^2(\Omega)} 
    + \skp{\nabla h^{n+1}_k}{\nabla q_k}_{L^2(\Omega)} \,
\end{align*}
has only the trivial solution.
Therefore, we set $\eta_k=h^{n+1}_k$ and $q_k = -\kappa g^{n+1}_k$. Summing the resulting two equations leads to
\begin{align*}
    0 = \frac{1}{\tau} \norm{h^{n+1}_k}_{L^2(\Omega)}^2
    + \bignorm{\sqrt{G} \nabla h^{n+1}_k}_{L^2(\Omega)}^2
    + \kappa \norm{g^{n+1}_k}_{L^2(\Omega)}^2 \,,
\end{align*}
which implies $h^{n+1}_k = g^{n+1}_k = 0$.

Next, we study the existence and uniqueness of $u^{n+1}_k, \mu^{n+1}_k \in \mathcal{S}_k$ solving \eqref{NUM:1*}--\eqref{NUM:2*}.
Since $h^{n+1}_k$ is already computed, we can define the functional
\begin{align*}
    &J: \Big\{ u_k \in \mathcal{S}_k \;\Big|\; \scp{u_k}{1}_{L^2(\Omega)} = \scp{u^n_k}{1}_{L^2(\Omega)} \Big\} \to \R,
    \\[1ex]
    &\begin{aligned}
    J(u_k) 
    &\coloneq \frac{1}{2\tau} \norm{u_k - u^n_k}_{H_{(0)}^{-1}(\Omega)}^2
    + \frac\eps2 \norm{\nabla u_k}_{L^2(\Omega)}^2
    + \eps^{-1} \skp{W_1(u_k)}{1}_{L^2(\Omega)}
    + \eps^{-1} \skp{W'_2(u^n_k)}{u_k}_{L^2(\Omega)}
    \\
    &\qquad
    - \skp{L \nabla h^{n+1}_k}{\nabla u_k}_{L^2(\Omega)} \,.
    \end{aligned}
\end{align*}
Obviously all summands in the definition of $J$ are convex. Moreover, invoking Poincar\'e's inequality, we deduce that the second summand is even strictly convex. This implies that $J_n$ is strictly convex.
As $J$ is also continuous, it follows that there exists a unique minimizer $u_k\in \mathcal{S}_k$ and a unique Lagrange multiplier $\lambda_k\in \mathbb{R}$, such that $u_k$ satisfies the Euler--Lagrange equation 
%(Standard arguments from optimization in finite dimensions, since $J$ is strictly convex and there is one affine-linear equality constraint)
\begin{align*}
    0 &= \eps \skp{\nabla u_k}{\nabla \vartheta_k}_{L^2(\Omega)}
    + \eps^{-1} \skp{W_1'(u_k) + W_2'(u^n_k)}{\vartheta_k}_{L^2(\Omega)}
    - \skp{L \nabla h^{n+1}_k}{\nabla \vartheta_k}_{L^2(\Omega)}
    \\
    &\quad
    + \skp{ \frac{1}{\tau} (-\Delta_k)^{-1}(u_k - u^n_k)}{ \vartheta_k}_{L^2(\Omega)}
    - \lambda_k \skp{1}{\vartheta_k}_{L^2(\Omega)} 
\end{align*}
for all test functions $\vartheta_k \in \mathcal{S}_k$.
Here, $(-\Delta_k)^{-1} \colon \, \mathcal{S}_k\cap H^1_{(0)}(\Omega) \to \mathcal{S}_k \cap H^1_{(0)}(\Omega)$ is the solution operator to the discrete Neumann--Poisson problem on the finite element space. 
By setting
\begin{align*}
    \mu_k = \lambda_k - \frac{1}{\tau} (-\Delta_k)^{-1}(u_k - u^n_k)
\end{align*}
and noting that $\nabla \lambda_k = 0$, we can split the previous equation into the two equations
\begin{align*}
    0 &= \frac{1}{\tau} 
    \skp{u_k - u^n_k}{w_k}_{L^2(\Omega)} 
    + \skp{\nabla \mu_k}{\nabla w_k}_{L^2(\Omega)} \,,
    \\
    0 &= - \skp{\mu_k}{\vartheta_k}_{L^2(\Omega)}
    + \eps \skp{\nabla u_k}{\nabla \vartheta_k}_{L^2(\Omega)}
    + \eps^{-1} \skp{W_1'(u_k) + W_2'(u^n_k)}{\vartheta_k}_{L^2(\Omega)}
    \\
    &\quad
    - \skp{L \nabla h^{n+1}_k}{\nabla \vartheta_k}_{L^2(\Omega)} \,,
\end{align*}
which hold for all $w_k,\vartheta_k \in \mathcal{S}_k$.
Hence, defining $u^{n+1}_k \coloneqq u_k$ and $\mu^{n+1}_k\coloneqq \mu_k$, we found the unique solution $(u^{n+1}_k,\mu^{n+1}_k,g^{n+1}_k,g^{n+1}_k) \in (\mathcal{S}_k)^4$ of the discrete system \eqref{NUM*}.
\end{proof}

Next, we present a mass conservation property for $u_k^{n+1}$ and $h_k^{n+1}$.
\begin{lemma}
For any $n\in\{0,\ldots,N-1\}$, the unique solution to \eqref{NUM*} satisfies
\begin{align}
    \label{NUM:mass*}
    \skp{u^{n+1}_k}{1}_{L^2(\Omega)} = \skp{u^{0}_k}{1}_{L^2(\Omega)},
    \qquad
    \skp{h^{n+1}_k}{1}_{L^2(\Omega)} = \skp{h^{0}_k}{1}_{L^2(\Omega)}.
\end{align}
\end{lemma}
\begin{proof}
This follows by choosing $w_k = 1$ in \eqref{NUM:1*} and $\eta_k=1$ in \eqref{NUM:3*}.
\end{proof}

Another important property is a discrete energy inequality. As before, no smallness assumption for $\tau>0$ is needed.

\begin{lemma}
The unique solution to \eqref{NUM*} satisfies
\begin{align}
    \nonumber
    \label{NUM:energy_global*}
    & \max_{m\in \{0,\ldots,N-1\}}
    \left( 
    \frac{\eps}{2} \norm{\nabla u^{m+1}_k}_{L^2(\Omega)}^2
    + \eps^{-1} \skp{W(u^{m+1}_k)}{1}_{L^2(\Omega)} 
    - \skp{L \nabla h^{m+1}_k}{\nabla u^{m+1}_k}_{L^2(\Omega)} \right)
    \\ \nonumber
    &\quad
    + \max_{m\in \{0,\ldots,N-1\}}
    \left( 
    \frac12 \norm{\sqrt{G} \nabla h^{m+1}_k}_{L^2(\Omega)}
    + \frac{\kappa}{2} \norm{g^{m+1}_k}_{L^2(\Omega)}^2
    \right)
    \\ \nonumber
    &\quad
    + \sum_{n=0}^{N-1} \left( \frac{\eps}{2}  \norm{\nabla u^{n+1}_k - \nabla u^n_k}_{L^2(\Omega)}^2
    + \frac12 \norm{\sqrt{G} \nabla (h^{n+1}_k - h^n_k)}_{L^2(\Omega)}
    + \frac{\kappa}{2} \norm{g^{n+1}_k - g^n_k}_{L^2(\Omega)}^2 
    \right)
    \\ \nonumber
    &\quad
    + \tau \sum_{n=0}^{N-1} \left( \norm{\mu^{n+1}_k}_{L^2(\Omega)}^2
    + \norm{\frac{1}{\tau} (h^{n+1}_k - h^n_k) }_{L^2(\Omega)}^2 \right)
    \\
    &\leq \nonumber
    \frac{\eps}{2} \norm{\nabla u^{0}_k}_{L^2(\Omega)}^2
    + \eps^{-1} \skp{W(u^{0}_k)}{1}_{L^2(\Omega)} 
    - \skp{L \nabla h^{0}_k}{\nabla u^{0}_k}_{L^2(\Omega)}
    + \frac12 \norm{\sqrt{G} \nabla h^{0}_k}_{L^2(\Omega)}
    \\
    &\quad
    + \frac{\kappa}{2} \norm{g^{0}_k}_{L^2(\Omega)}^2  \, .
\end{align}
\end{lemma}
\begin{proof}
By choosing $w_k = \mu^{n+1}_k$ in \eqref{NUM:1*}, $\vartheta=\frac{1}{\tau} (u^{n+1}_k - u^n_k)$ in \eqref{NUM:2*}, $\eta_k = \frac{1}{\tau} (h^{n+1}_k - h^n_k)$ in \eqref{NUM:3*}, and using \eqref{NUM:4*}, the convex-concave splitting of $W = W_1 + W_2$, and the elementary identity 
\begin{align*}
    2 b (b-a) = b^2 - a^2 + (b-a)^2 \quad \text{for all $a,b \in \mathbb{R}$} \, ,
\end{align*}
we obtain
\begin{align}
    \label{NUM:energy*} \nonumber
    &\frac{\eps}{2} \norm{\nabla u^{n+1}_k}_{L^2(\Omega)}^2
    + \frac{\eps}{2} \norm{\nabla u^{n+1}_k - \nabla u^n_k}_{L^2(\Omega)}^2
    + \eps^{-1} \skp{W(u^{n+1}_k)}{1}_{L^2(\Omega)} 
    - \skp{L \nabla h^{n+1}_k}{\nabla u^{n+1}_k}_{L^2(\Omega)}
    \\ \nonumber
    &\quad 
    + \frac12 \norm{\sqrt{G} \nabla h^{n+1}_k}_{L^2(\Omega)}
    + \frac12 \norm{\sqrt{G} \nabla (h^{n+1}_k - h^n_k)}_{L^2(\Omega)}
    + \frac{\kappa}{2} \norm{g^{n+1}_k}_{L^2(\Omega)}^2
    + \frac{\kappa}{2} \norm{g^{n+1}_k - g^n_k}_{L^2(\Omega)}^2
    \\ \nonumber
    &\quad
    + \tau \norm{\mu^{n+1}_k}_{L^2(\Omega)}^2
    + \tau \norm{\frac{1}{\tau} (h^{n+1}_k - h^n_k) }_{L^2(\Omega)}^2
    \\
    &\leq 
    \frac{\eps}{2} \norm{\nabla u^{n}_k}_{L^2(\Omega)}^2
    + \eps^{-1} \skp{W(u^{n}_k)}{1}_{L^2(\Omega)} 
    - \skp{L \nabla h^{n}_k}{\nabla u^{n}_k}_{L^2(\Omega)}
    + \frac12 \norm{\sqrt{G} \nabla h^{n}_k}_{L^2(\Omega)}
    + \frac{\kappa}{2} \norm{g^{n}_k}_{L^2(\Omega)}^2  \, .
\end{align}
Then, the claim follows from summation over $n \in \{0,\ldots,m\}$, $m\in\{0,\ldots,N-1\}$, in \eqref{NUM:energy*}, and by taking the maximum over $m\in\{0,\ldots,N-1\}$.
\end{proof}

% \subsection{Computational aspects}
% \begin{itemize}
% \item describe initial setting
% \item how to solve the discrete system in practice
% \item adaptive time stepping strategy
% %\item
% \end{itemize}

\subsection{Numerical experiments}
\label{sec:numerics}

In the following, we present two numerical experiments based on \eqref{NUM*} to study the behavior of the model and, in particular, the influence of the physical parameters. The computations always use a time step size of $\tau = 10^{-4}$ and a uniform Friedrichs--Keller triangulation with $160 \times 160$ vertices.

We employ a double-well potential that coincides with the logarithmic Flory--Huggins potential \ref{ASS:POT:SING} on the interval $[-1 + \delta, 1 - \delta]$ with $\delta = 0.02$, and extend it to $\mathbb{R}$ using a second-order Taylor expansion of its convex part. Specifically, we define
\begin{align*}
    W(s) = 
    \begin{cases}
    W_{\mathrm{log},1}(\delta-1) + (s+1-\delta) W'_{\mathrm{log},1}(\delta-1)
    + \frac12 (s+1-\delta)^2 W''_{\mathrm{log},1}(\delta-1)
    - \frac{5}{2} s^2
    & s \leq \delta-1,
    \\
    W_{\mathrm{log}}(s) & s\in(\delta-1,1-\delta),
    \\
    W_{\mathrm{log},1}(1-\delta) + (s-1+\delta) W'_{\mathrm{log},1}(1-\delta)
    + \frac12 (s-1+\delta)^2 W''_{\mathrm{log},1}(1-\delta)
    - \frac{5}{2} s^2
    & s \geq  1-\delta ,
    \end{cases}
\end{align*}
where $W_{\mathrm{log}}(s) = W_{\mathrm{log},1}(s) - \frac{5}{2} s^2$, with its convex part being defined as 
\[
W_{\mathrm{log},1}(s) = 2 \big( (1+s)\ln(1+s) + (1-s)\ln(1-s) \big)
\]
for all $s\in[-1,1]$. 
This corresponds to the setting in \ref{ASS:POT:SING} with parameters $\mu_0 = 0$, $\theta = 4$, and $\theta_c = 5$. 
The potential attains its global minima at $s\approx \pm 0.71$.
%For this potential, we apply a convex–concave splitting.

\begin{table}
\centering
\begin{tabular}{c||c|c|c|c|c}%{c||c|c|c|c|c|c|c}
     Parameter  & period & coupling strength & surface tension & bending rigidity & interface par.
     \\
     & & $\Lambda = |L|$ & $\sigma = |G| $ %& $a_4/\varepsilon$ 
     & $\kappa$ & $\varepsilon$
     \\[1ex] \hline 
     Physical unit &
     m & J/m & J/$\mathrm{m}^2$ %& J/$\mathrm{m}^2$ 
     & J & J
     \\[1ex] 
     Literature & 
     $10^{-6}$  & $4.9 \cdot 10^{-12}$ & $\left[5\cdot 10^{-6}, 10^{-4}\right]$ %& $10^{-5}$  
     & $10^{-19}$ & $5\cdot 10^{-19}$
     \\[1ex]
     Simulation 1 &
     $10^{-6}$  & $\left[10^{-11}, 10^{-10}\right]$  & $5\cdot 10^{-5}$  %& $10^{-5}$ J/$\mathrm{m}^2$ 
     & $\left[ 3.75\cdot10^{-19}, 10^{-18}\right]$ & $5\cdot 10^{-19}$ 
     \\[1ex]
     Simulation 2 &
     $10^{-6}$  & $[10^{-11}, 3\cdot10^{-11}]$  & $\left[5\cdot 10^{-7}, 5\cdot 10^{-4}\right]$  %& $10^{-5}$ J/$\mathrm{m}^2$ 
     & $10^{-19}$ & $5\cdot 10^{-19}$ 
\end{tabular}
\caption{Physical parameters (modulation period, coupling strength, surface tension, bending rigidity, and interface parameter) with their physical units, the estimated values from \cite{Komura:Langmuir:2006, RaedlerFSS1995_parameters, rozovsky2005formation, Seifert1995_parameters}, and the parameters from 
\eqref{NUM:parameters} and \eqref{NUM:parameters2} converted to physical units, see \eqref{NUM:parameters_converted}. 
% For the conversion, we use the pattern period $10^{-6}\,\mathrm{m}$ as a characteristic length and a characteristic energy of $5\cdot10^{-17}\,\mathrm{J}$, corresponding to the ratio between the physical and simulated values of the interface parameter $\varepsilon$.
Here, ``Simulation 1'' and `Simulation 2'' correspond to our numerical computations based on the parameter choices \eqref{NUM:parameters} and \eqref{NUM:parameters2}, respectively.}
\label{tab:physical_parameters}
\end{table}
%% Eref = 5*10^{-19}J/\varepsilon
%% xref = 10^{-6}m
%% Lphys = L * Eref/xref
%% sigmaphys = sigma * Eref/(xref^2)
%% kphys = k * Eref

In Table~\ref{tab:physical_parameters}, we compare estimated physical values from the literature \cite{Komura:Langmuir:2006, RaedlerFSS1995_parameters, rozovsky2005formation, Seifert1995_parameters} with the parameter choices of our two parameter settings: \eqref{NUM:parameters} and \eqref{NUM:parameters2}. In all cases, we consider isotropic choices for the surface tension $G = \sigma \mathrm{I}$ and coupling strength $L = \Lambda \mathrm{I}$, respectively. Here, $\mathrm{I}\in\mathbb{R}^{2\times2}$ is the identity matrix.
To convert our non-dimensional parameters from \eqref{NUM:parameters} and \eqref{NUM:parameters2} to physical units, we introduce a characteristic energy scale $E_{c} = 5\cdot10^{-19}\, \mathrm{J}/\varepsilon$. With this choice, an interface parameter $\varepsilon=0.01$ corresponds to the physical value $5\cdot 10^{-19}\, \mathrm{J}$, as reported in \cite{Komura:Langmuir:2006, RaedlerFSS1995_parameters, Seifert1995_parameters}. 
Following \cite{Komura:Langmuir:2006, rozovsky2005formation}, we take the pattern modulation period as the characteristic length $x_{c} = 10^{-6}\, \mathrm{m}$. Using these scales, we then compute the physical values of the parameters $\Lambda$, $\sigma$ and $\kappa$ as follows:
\begin{align} \label{NUM:parameters_converted}
    \Lambda_\textit{phys} = \frac{\Lambda E_{c}}{x_{c}},
    \quad
    \sigma_\textit{phys} = \frac{\sigma E_{c}}{x_{c}^2}, 
    \quad
    \kappa_\textit{phys} = \kappa E_{c} .
\end{align}

According to \cite{Komura:Langmuir:2006}, if the condition $\Lambda^2 > \sigma\varepsilon$ is satisfied, a homogeneous initial state may become unstable with respect to long-wave fluctuations. 
Depending on the values of $\Lambda$ and $\sigma$, three characteristic types of patterns can appear. By increasing $\sigma$ (or alternatively, by decreasing $\Lambda$), one typically observes a transition from striped structures to ordered dots, and finally to disordered dotted domains. Stripes tend to form when $\sigma$ is small (or $\Lambda$ is large). This behavior agrees with experimental observations in mixed lipid bilayers reported in \cite{rozovsky2005formation}. A more detailed classification of these regimes is shown in \cite[Figure~3]{Komura:Langmuir:2006}.

In addition, the mean value $\mean{u_0}$ of the order parameter also influences the type of pattern that develops (see \cite[Figure~2]{Komura:Langmuir:2006}). Striped phases are typically observed when $\mean{u_0}$ is close to zero, while disordered dotted structures appear when $\mean{u_0}$ approaches the physical bounds $\pm1$. For intermediate values of $\mean{u_0}$, ordered dotted patterns are more likely. Motivated by these findings, we examine two parameter studies. First, in \eqref{NUM:parameters}, we vary $\Lambda$ and $\kappa$ for two different mean values $\mean{u_0}$. Afterwards, in \eqref{NUM:parameters2}, we vary $\sigma$ and $\Lambda$ while again considering two different choices of $\mean{u_0}$. We note that the coupling condition $\Lambda^2 > \sigma\varepsilon$ is fulfilled in all cases except for $\sigma=10$, $\Lambda=0.2$ in the second parameter setting \eqref{NUM:parameters2}, where the ranges $\sigma\in[0.01,10]$ and $\Lambda\in[0.2,0.6]$ are studied.

In the first setting, the general idea is to study the influence of the bending rigidity $\kappa$ and the (isotropic) coupling strength $\Lambda$. Overall, the initial data are set to $h_0 = 0$, while $u_0$ is given by a random perturbation of amplitude $0.2$ around a prescribed mean value $\mean{u_0} \in \{0.1, 0.3\}$. Also, we fix the interface parameter, the (isotropic) surface tension, the bending rigidity, and the coupling strength to be
\begin{gather}\label{NUM:parameters}
    \varepsilon = 0.01, \quad 
    \sigma = 1, \quad
    \kappa \in [0.75 \varepsilon, 2\varepsilon], \quad
    \Lambda \in [0.2, 2]. 
\end{gather}

%%% needed for \includegraphics[height=0.98\ht\mycolumn]{...}
\newsavebox{\mycolumn}
\sbox{\mycolumn}{%
\begin{minipage}{.33\textwidth}
\includegraphics[width=\linewidth]{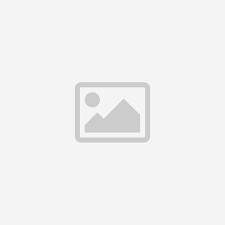}\\[-0em]
\includegraphics[width=\linewidth]{images/dummy_figure.png}
\end{minipage}%
}

\begin{figure}[htb!]
\centering
\begin{minipage}{.86\textwidth}
\setlength{\tabcolsep}{2pt} % default ~6pt
\includegraphics[width=.16\textwidth,trim={0.9cm 1.6cm 0.9cm 1.6cm}, clip]{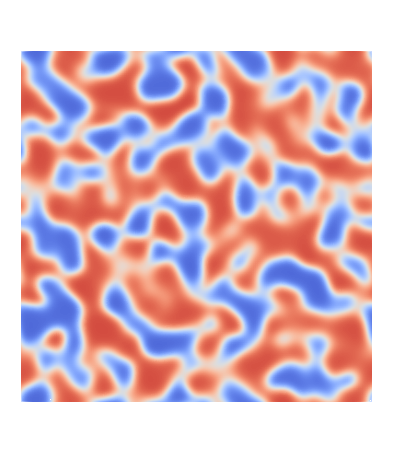} 
\includegraphics[width=.16\textwidth,trim={0.9cm 1.6cm 0.9cm 1.6cm}, clip]{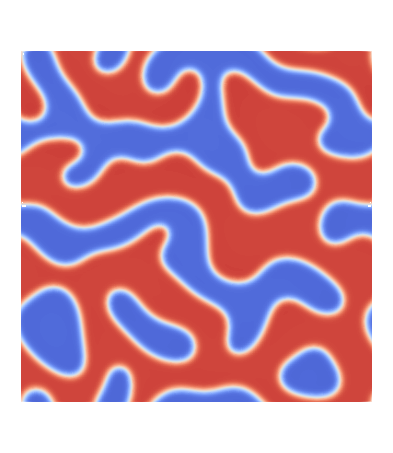} 
\includegraphics[width=.16\textwidth,trim={0.9cm 1.6cm 0.9cm 1.6cm}, clip]{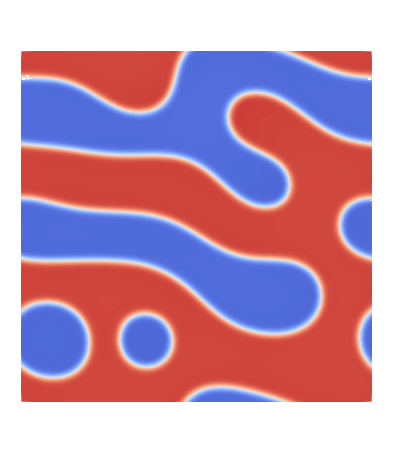}
\includegraphics[width=.16\textwidth,trim={0.9cm 1.6cm 0.9cm 1.6cm}, clip]{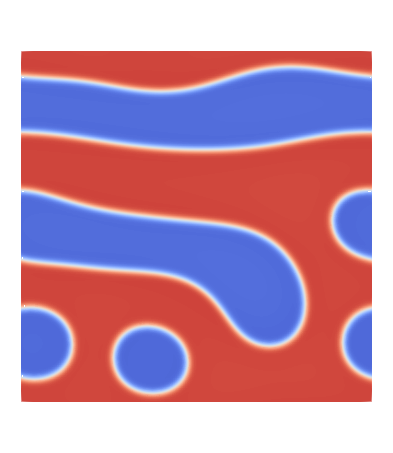}
\includegraphics[width=.16\textwidth,trim={0.9cm 1.6cm 0.9cm 1.6cm}, clip]{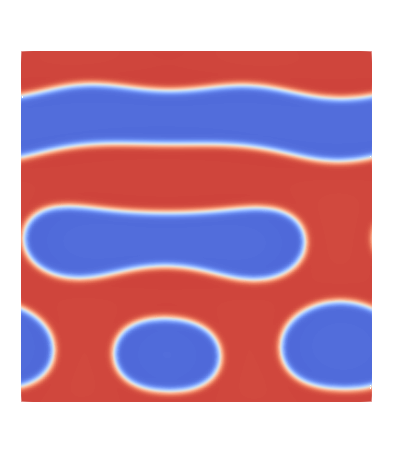}
\\
\includegraphics[width=.16\textwidth,trim={0.2cm 2cm 0.8cm 2cm}, clip]{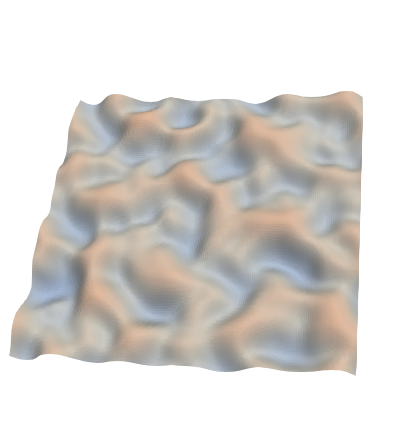} 
\includegraphics[width=.16\textwidth,trim={0.2cm 2cm 0.8cm 2cm}, clip]{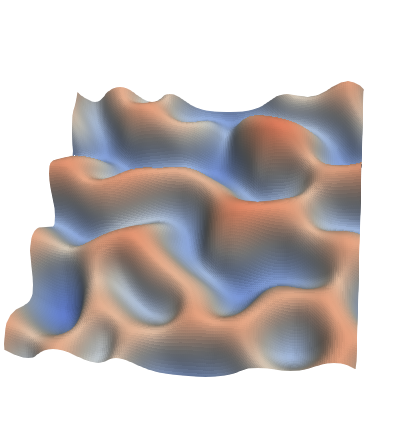} 
\includegraphics[width=.16\textwidth,trim={0.2cm 2cm 0.8cm 2cm}, clip]{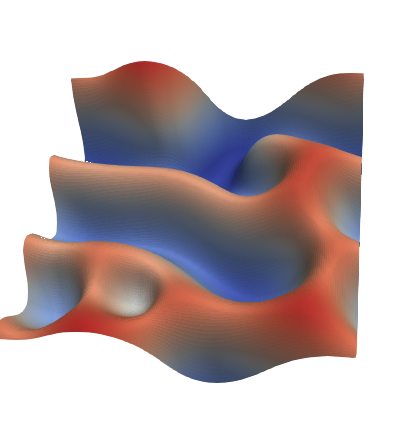}
\includegraphics[width=.16\textwidth,trim={0.2cm 2cm 0.8cm 2cm}, clip]{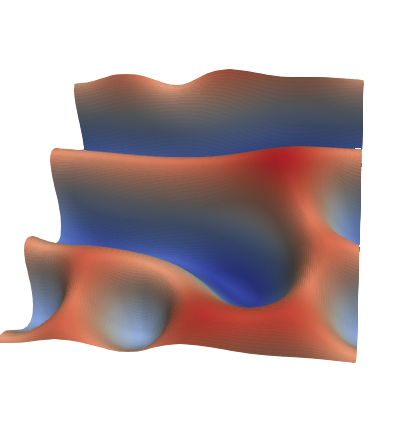}
\includegraphics[width=.16\textwidth,trim={0.2cm 2cm 0.8cm 2cm}, clip]{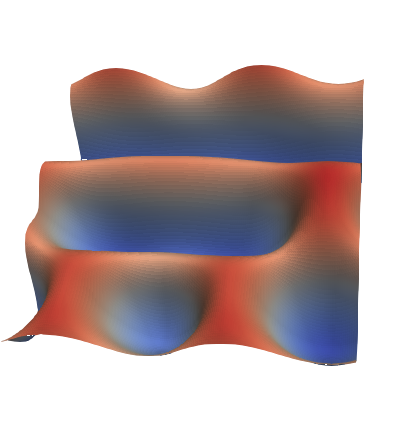}
\end{minipage}%
\hspace{-0.15\textwidth}  %% adapt gap
\begin{minipage}{.06\textwidth}
\mbox{}\vspace{0.1cm}

\includegraphics[height=0.45\ht\mycolumn]{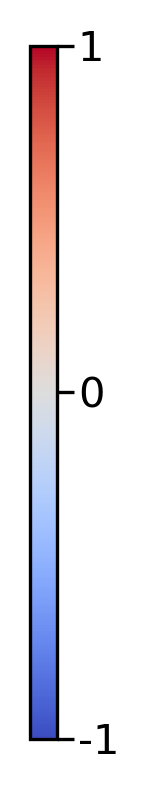}
\\[-0.1cm]
% \vspace{-0.2cm}
\includegraphics[height=0.40\ht\mycolumn]{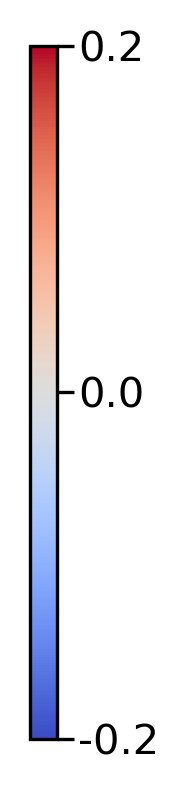}
\end{minipage}%
\caption{Time evolution of the order parameter $u$ (first row) and the membrane height $h$ (second row) for $\sigma=1$, $\kappa=\varepsilon=0.01$, $\Lambda=0.6$, and $\mean{u_0}=0.1$ at times $t\in\{0.002, 0.01, 0.1, 1, 10\}$.}
\label{fig:test1_long}
\end{figure}

In Figure~\ref{fig:test1_long}, we present a long-time evolution of the order parameter $u$ and the height profile $h$ of the membrane until time $t=10$. In this test, the bending rigidity is fixed as $\kappa=\varepsilon$ and the coupling strength is set to $\Lambda=0.6$, while the other parameters are chosen as in \eqref{NUM:parameters}.
Starting from an initial mixture with small random perturbations around the mean value $\mean{u_0}=0.1$, the system undergoes phase separation. The order parameter $u$ forms snake-like patterns that gradually merge and thicken over time, eventually forming stripe-like patterns including few dots. The membrane height $h$ follows the spatial pattern of $u$, forming valleys where $u\approx -1$ and hills where $u\approx +1$.

% \begin{figure}[htb!]
%   \centering
%   \resizebox{!}{0.20\textheight}{%
%     \begin{minipage}{\textwidth}
%       \centering
%       \includegraphics[width=\textwidth]{images/dummy_figure.png}
%     \end{minipage}
%   }
%   \caption{Caption}
%   \label{fig:placeholder}   
% \end{figure}

\begin{figure}[htb!]
\centering
\begin{minipage}{.86\textwidth}
\setlength{\tabcolsep}{2pt} % default ~6pt
\begin{tabular}{cM{.16\textwidth}M{.16\textwidth}M{.16\textwidth}M{.16\textwidth}M{.16\textwidth}M{.16\textwidth}}
& $\Lambda=0.2$& $\Lambda=0.4$ & $\Lambda=0.6$  & $\Lambda=0.8$ & $\Lambda=1.0$ & $\Lambda=2.0$    
\\ 
\rotatebox{90}{\hspace{-0.6cm} $\kappa=0.75\varepsilon$} 
&\includegraphics[width=.16\textwidth,trim={0.9cm 1.6cm 0.9cm 1.6cm}, clip]{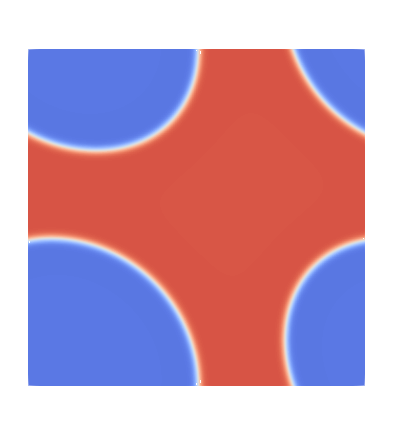} 
&\includegraphics[width=.16\textwidth,trim={0.9cm 1.6cm 0.9cm 1.6cm}, clip]{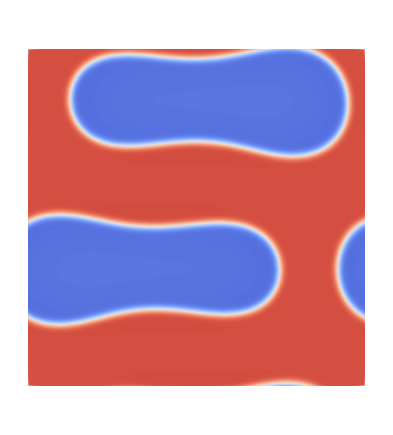} 
&\includegraphics[width=.16\textwidth,trim={0.9cm 1.6cm 0.9cm 1.6cm}, clip]{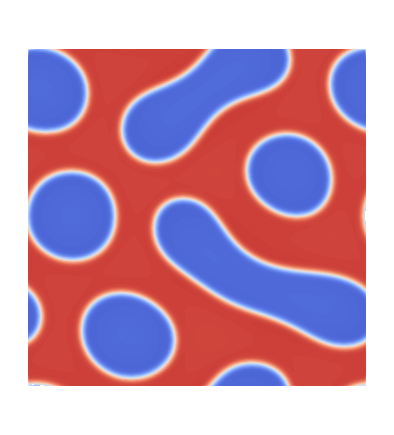} 
& \includegraphics[width=.16\textwidth,trim={0.9cm 1.6cm 0.9cm 1.6cm}, clip]{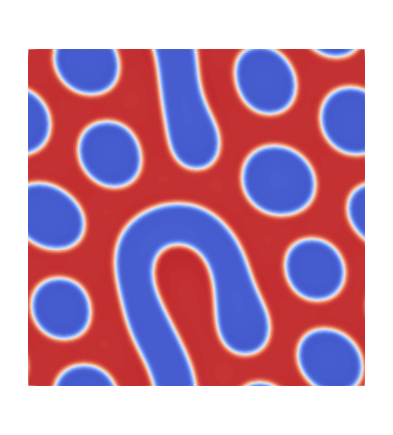} 
& \includegraphics[width=.16\textwidth,trim={0.9cm 1.6cm 0.9cm 1.6cm}, clip]{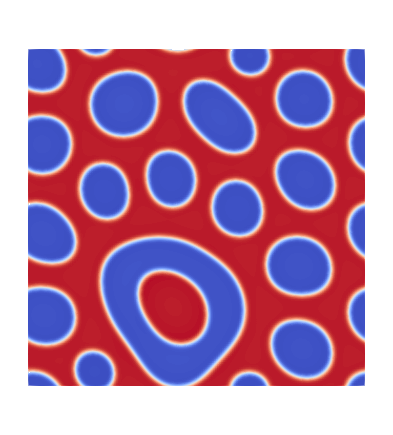} 
& \includegraphics[width=.16\textwidth,trim={0.9cm 1.6cm 0.9cm 1.6cm}, clip]{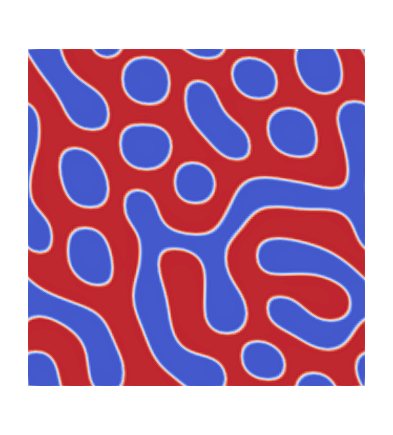} 
\\[-0em]
\rotatebox{90}{\hspace{-0.6cm} $\kappa=\varepsilon$} 
&\includegraphics[width=.16\textwidth,trim={0.9cm 1.6cm 0.9cm 1.6cm}, clip]{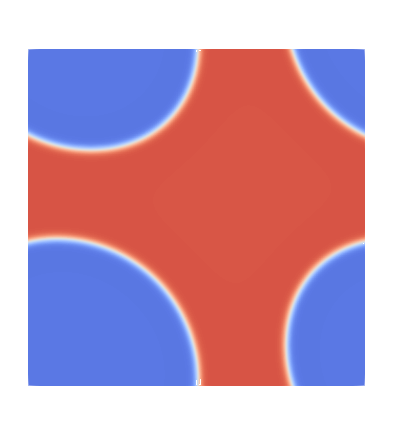} 
&\includegraphics[width=.16\textwidth,trim={0.9cm 1.6cm 0.9cm 1.6cm}, clip]{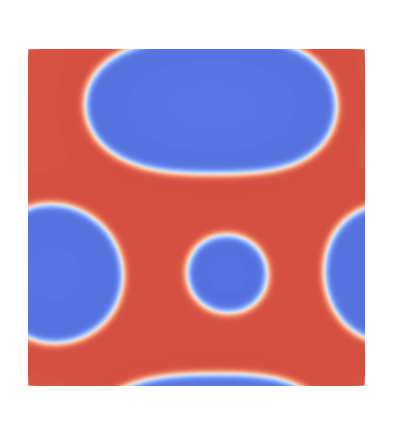} 
&\includegraphics[width=.16\textwidth,trim={0.9cm 1.6cm 0.9cm 1.6cm}, clip]{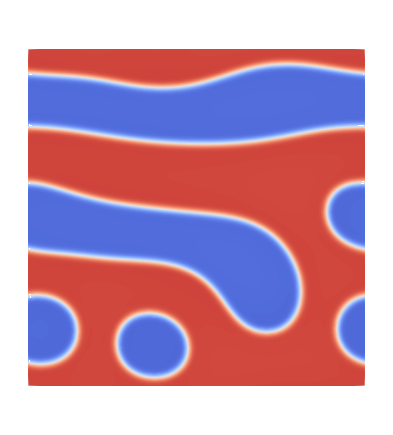} 
& \includegraphics[width=.16\textwidth,trim={0.9cm 1.6cm 0.9cm 1.6cm}, clip]{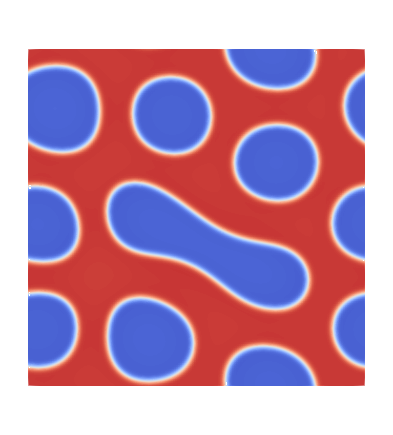} 
& \includegraphics[width=.16\textwidth,trim={0.9cm 1.6cm 0.9cm 1.6cm}, clip]{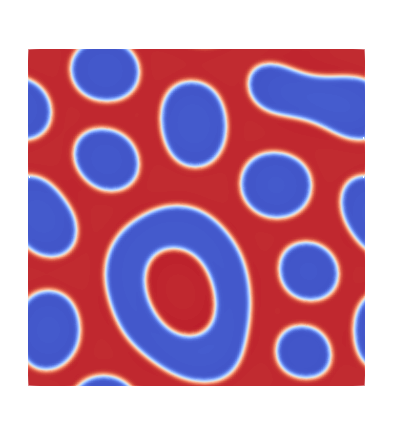} 
& \includegraphics[width=.16\textwidth,trim={0.9cm 1.6cm 0.9cm 1.6cm}, clip]{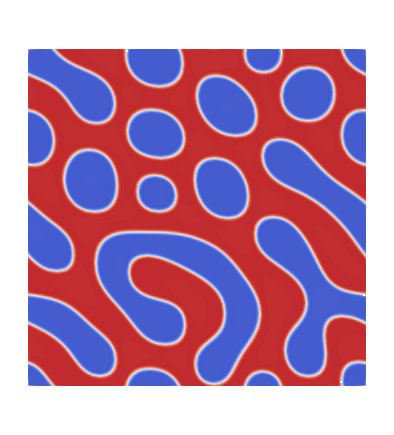} 
\\[-0em]
\rotatebox{90}{\hspace{-0.6cm} $\kappa=1.5\varepsilon$}
&\includegraphics[width=.16\textwidth,trim={0.9cm 1.6cm 0.9cm 1.6cm}, clip]
{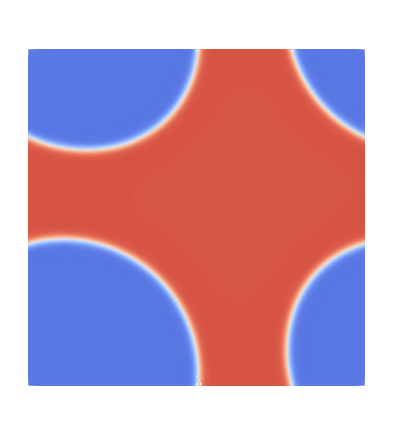}  
&\includegraphics[width=.16\textwidth,trim={0.9cm 1.6cm 0.9cm 1.6cm}, clip]{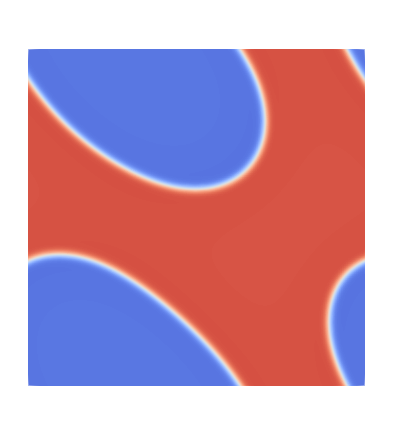} 
&\includegraphics[width=.16\textwidth,trim={0.9cm 1.6cm 0.9cm 1.6cm}, clip]{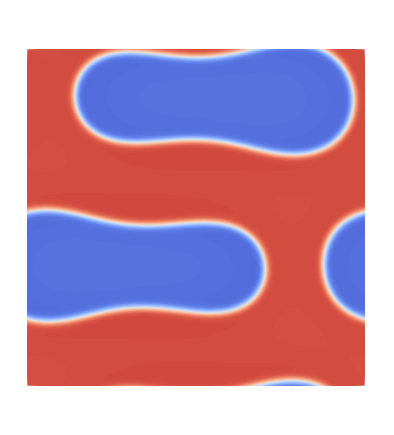} 
& \includegraphics[width=.16\textwidth,trim={0.9cm 1.6cm 0.9cm 1.6cm}, clip]{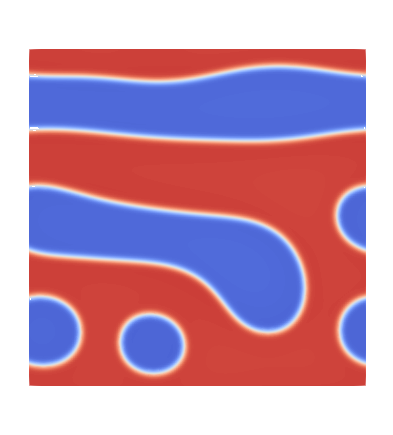} 
& \includegraphics[width=.16\textwidth,trim={0.9cm 1.6cm 0.9cm 1.6cm}, clip]{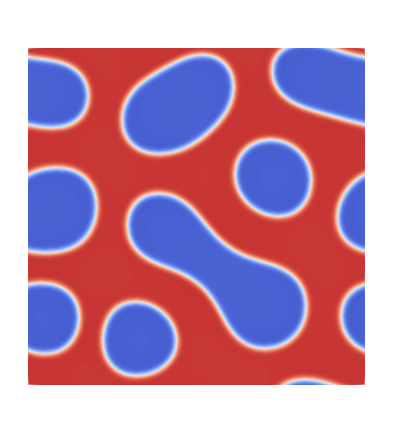} 
& \includegraphics[width=.16\textwidth,trim={0.9cm 1.6cm 0.9cm 1.6cm}, clip]{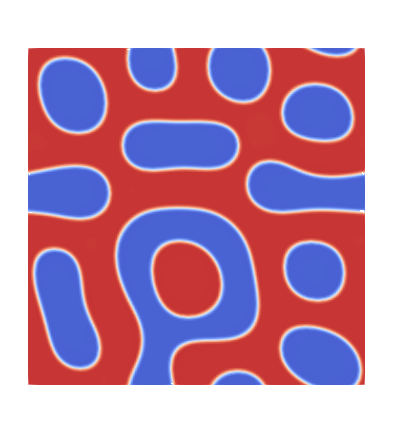} 
\\[-0em]
\rotatebox{90}{\hspace{-0.6cm} $\kappa=2\varepsilon$} 
&\includegraphics[width=.16\textwidth,trim={0.9cm 1.6cm 0.9cm 1.6cm}, clip]
{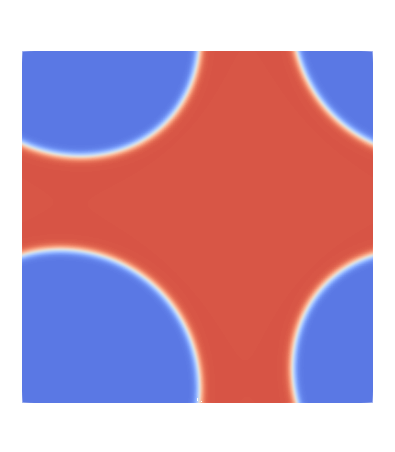} 
&\includegraphics[width=.16\textwidth,trim={0.9cm 1.6cm 0.9cm 1.6cm}, clip]{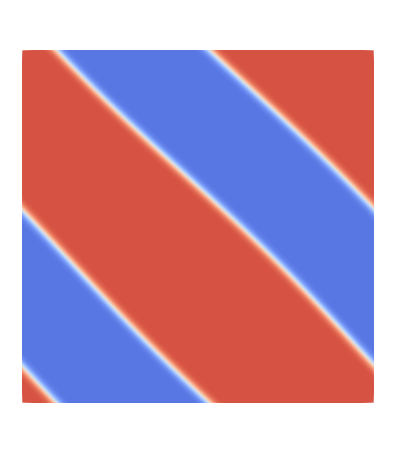} 
&\includegraphics[width=.16\textwidth,trim={0.9cm 1.6cm 0.9cm 1.6cm}, clip]{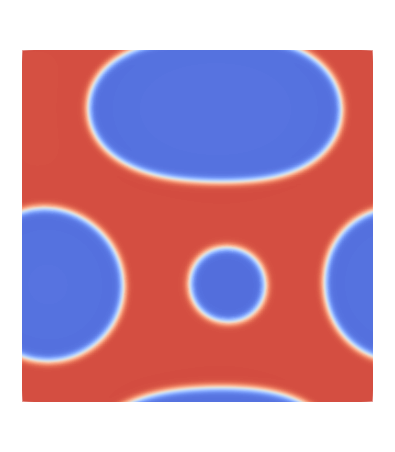} 
& \includegraphics[width=.16\textwidth,trim={0.9cm 1.6cm 0.9cm 1.6cm}, clip]{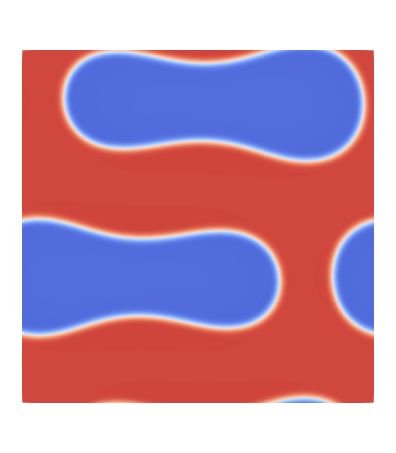} 
& \includegraphics[width=.16\textwidth,trim={0.9cm 1.6cm 0.9cm 1.6cm}, clip]{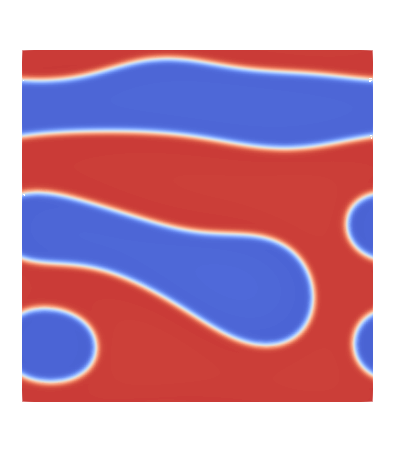} 
& \includegraphics[width=.16\textwidth,trim={0.9cm 1.6cm 0.9cm 1.6cm}, clip]{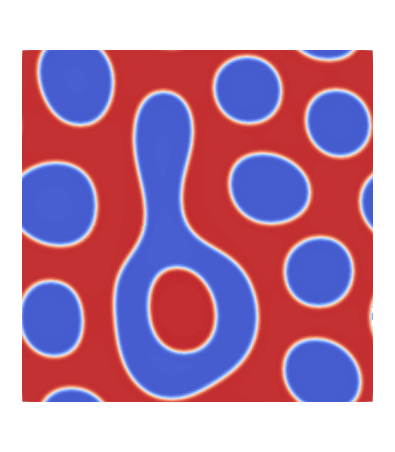} 
\end{tabular}
\end{minipage}%
\hspace{0.035\textwidth}  %% adapt gap
\begin{minipage}{.06\textwidth}
\mbox{}\vspace{.3cm}

\includegraphics[height=0.88\ht\mycolumn]{images/colorbar.png}
\end{minipage}%
\caption{Order parameter $u$ at time $t=1$ for different values of $\kappa \in [0.75\varepsilon, 2\varepsilon]$ and $\Lambda\in[0.2, 2]$. The other parameters are fixed as $\sigma = 1$, $\varepsilon=0.01$, and $\mean{u_0} = 0.1$.
With increasing $\Lambda$, the pattern transitions from dotted 
(small $\Lambda$) to striped (intermediate $\Lambda$) and unstructured 
(large $\Lambda$) phases. Increasing $\kappa$ shifts this transition toward larger $\Lambda$.
}\vspace{-1em}
\label{fig:test1}
\end{figure}

We now study the influence of the coupling strength $\Lambda\in[0.2, 2]$ and the bending rigidity $\kappa\in[0.75\varepsilon, 2\varepsilon]$ on the pattern formation. The mean value of the order parameter is fixed as $\mean{u_0} = 0.1$ and the parameters are specified in \eqref{NUM:parameters}.
The corresponding order parameters $u$ for different values of $\Lambda$ and $\kappa$ at time $t=1$ are shown in Figure~\ref{fig:test1}.
To discuss the effect of varying $\Lambda$, we fix $\kappa$ (i.e., consider one row in Figure~\ref{fig:test1}). 
We observe a transition from a purely dotted phase for small $\Lambda$ to a rather striped pattern for intermediate $\Lambda$, and finally to an mixed pattern for large $\Lambda$, containing dots and snake-like structures. Increasing $\kappa$ shifts this transition toward larger values of $\Lambda$.
Moreover, we observe the following general trends. Very large values of $\Lambda$ cause the order parameter $u$ to exceed the physically relevant interval $[-1,1]$. Moreover, if $\kappa$ is chosen too small compared to $\varepsilon$, the resulting diffuse interface becomes very thin, so that the spatial resolution of $160\times160$ vertices is too coarse to resolve it accurately.

\begin{figure}[htb!]
\centering
\begin{minipage}{.86\textwidth}
\setlength{\tabcolsep}{2pt} % default ~6pt
\begin{tabular}{cM{.16\textwidth}M{.16\textwidth}M{.16\textwidth}M{.16\textwidth}M{.16\textwidth}M{.16\textwidth}}
& $\Lambda=0.2$& $\Lambda=0.4$ & $\Lambda=0.6$  & $\Lambda=0.8$ & $\Lambda=1.0$ & $\Lambda=2.0$    
\\ 
\rotatebox{90}{\hspace{-0.6cm} $\kappa=0.75\varepsilon$} 
% & \includegraphics[width=.16\textwidth]{images/dummy_figure.png}
&\includegraphics[width=.16\textwidth,trim={0.9cm 1.6cm 0.9cm 1.6cm}, clip]{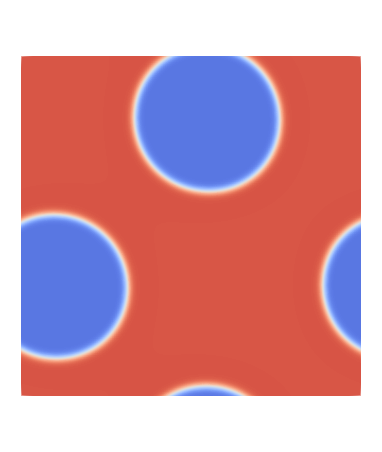} 
&\includegraphics[width=.16\textwidth,trim={0.9cm 1.6cm 0.9cm 1.6cm}, clip]{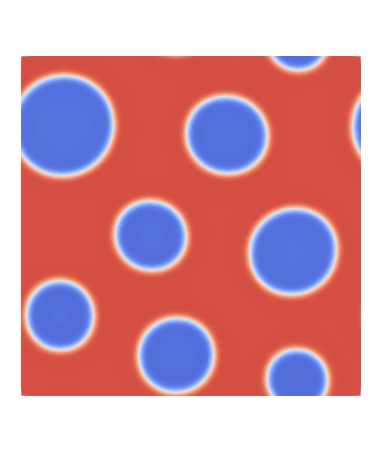} 
&\includegraphics[width=.16\textwidth,trim={0.9cm 1.6cm 0.9cm 1.6cm}, clip]{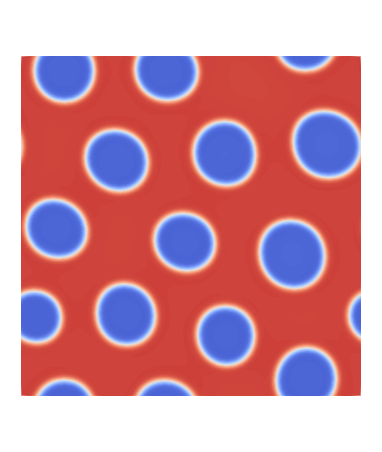} 
& \includegraphics[width=.16\textwidth,trim={0.9cm 1.6cm 0.9cm 1.6cm}, clip]{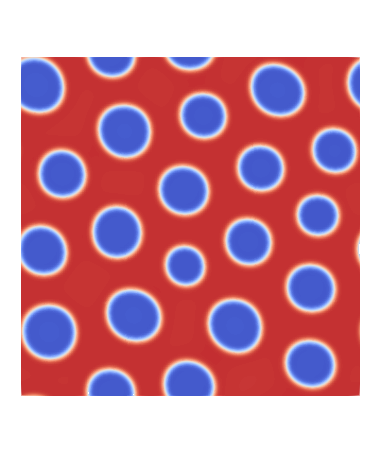} 
& \includegraphics[width=.16\textwidth,trim={0.9cm 1.6cm 0.9cm 1.6cm}, clip]{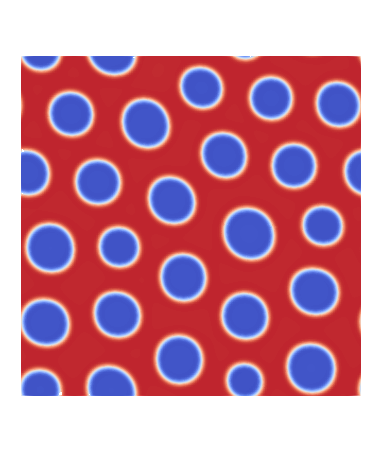} 
& \includegraphics[width=.16\textwidth,trim={0.9cm 1.6cm 0.9cm 1.6cm}, clip]{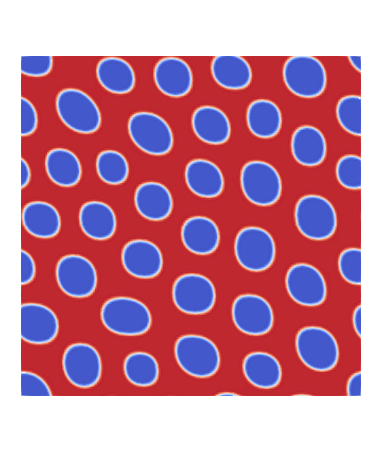} 
\\[-0em]
\rotatebox{90}{\hspace{-0.6cm} $\kappa=\varepsilon$} 
&\includegraphics[width=.16\textwidth,trim={0.9cm 1.6cm 0.9cm 1.6cm}, clip]{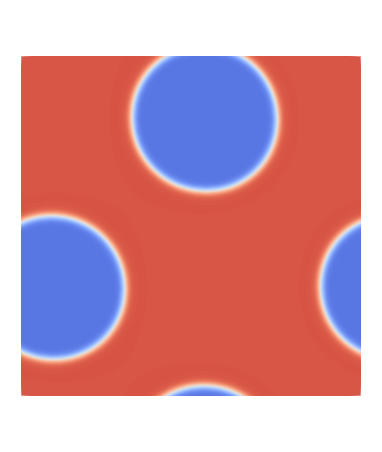} 
&\includegraphics[width=.16\textwidth,trim={0.9cm 1.6cm 0.9cm 1.6cm}, clip]{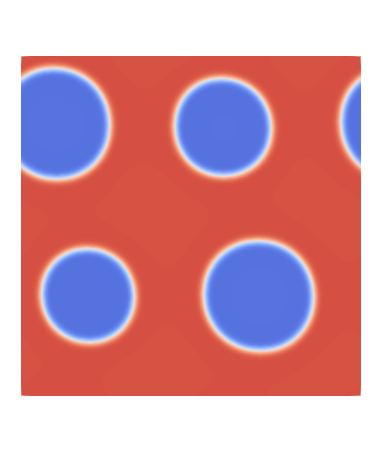} 
&\includegraphics[width=.16\textwidth,trim={0.9cm 1.6cm 0.9cm 1.6cm}, clip]{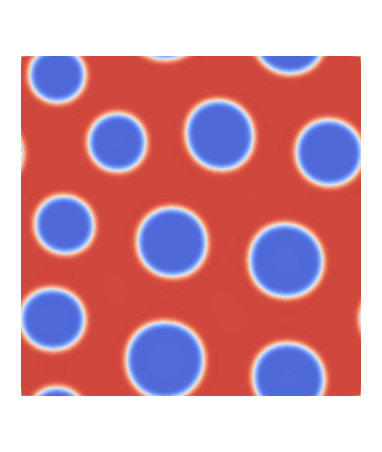} 
& \includegraphics[width=.16\textwidth,trim={0.9cm 1.6cm 0.9cm 1.6cm}, clip]{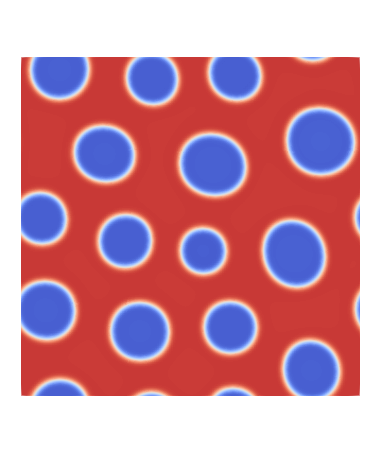} 
& \includegraphics[width=.16\textwidth,trim={0.9cm 1.6cm 0.9cm 1.6cm}, clip]{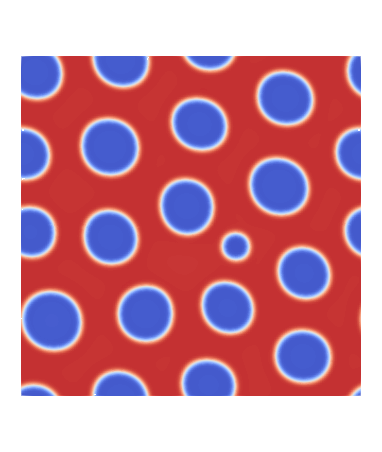} 
& \includegraphics[width=.16\textwidth,trim={0.9cm 1.6cm 0.9cm 1.6cm}, clip]{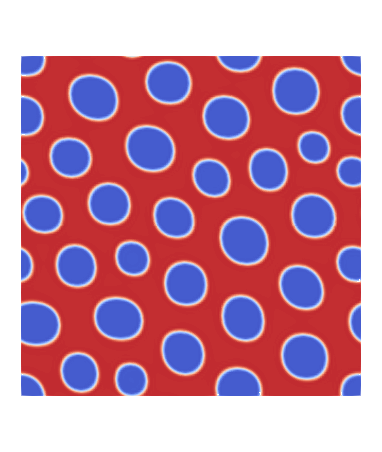} 
\\[-0em]
\rotatebox{90}{\hspace{-0.6cm} $\kappa=1.5\varepsilon$}
&\includegraphics[width=.16\textwidth,trim={0.9cm 1.6cm 0.9cm 1.6cm}, clip]{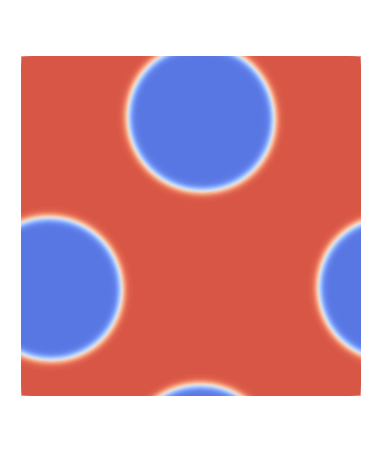} 
&\includegraphics[width=.16\textwidth,trim={0.9cm 1.6cm 0.9cm 1.6cm}, clip]{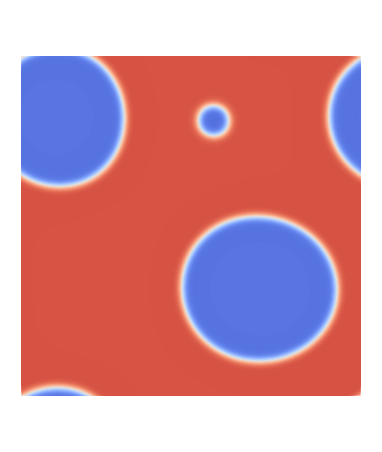} 
&\includegraphics[width=.16\textwidth,trim={0.9cm 1.6cm 0.9cm 1.6cm}, clip]{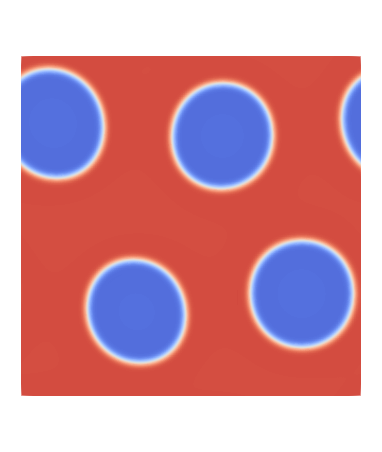} 
& \includegraphics[width=.16\textwidth,trim={0.9cm 1.6cm 0.9cm 1.6cm}, clip]{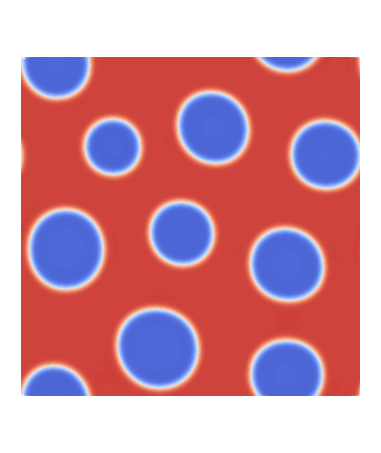} 
& \includegraphics[width=.16\textwidth,trim={0.9cm 1.6cm 0.9cm 1.6cm}, clip]{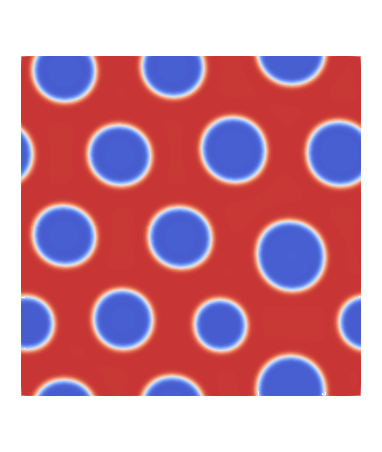} 
& \includegraphics[width=.16\textwidth,trim={0.9cm 1.6cm 0.9cm 1.6cm}, clip]{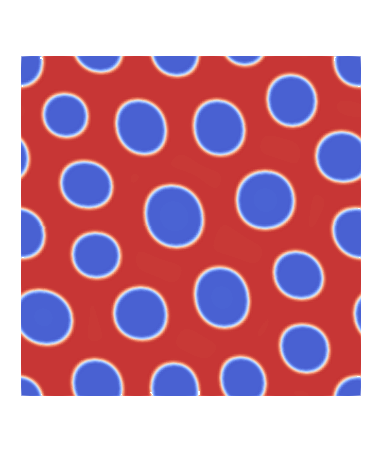} 
\\[-0em]
\rotatebox{90}{\hspace{-0.6cm} $\kappa=2\varepsilon$} 
&\includegraphics[width=.16\textwidth,trim={0.9cm 1.6cm 0.9cm 1.6cm}, clip]{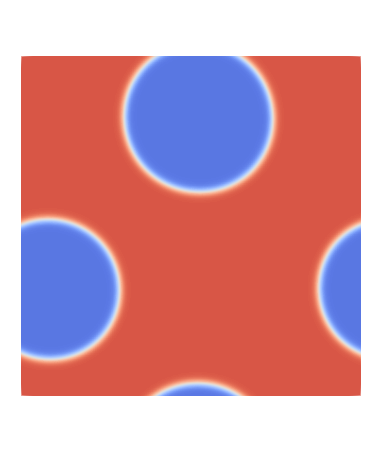} 
&\includegraphics[width=.16\textwidth,trim={0.9cm 1.6cm 0.9cm 1.6cm}, clip]{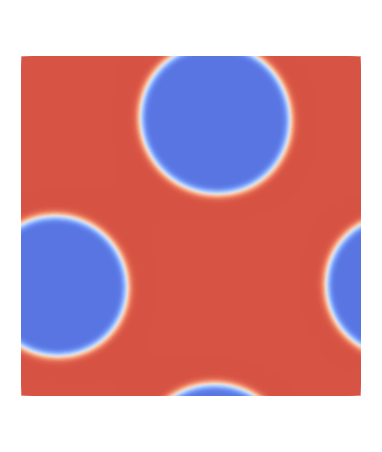} 
&\includegraphics[width=.16\textwidth,trim={0.9cm 1.6cm 0.9cm 1.6cm}, clip]{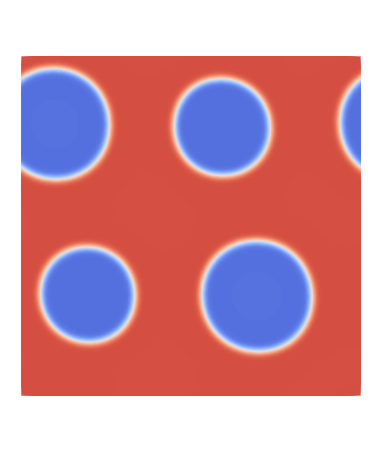} 
& \includegraphics[width=.16\textwidth,trim={0.9cm 1.6cm 0.9cm 1.6cm}, clip]{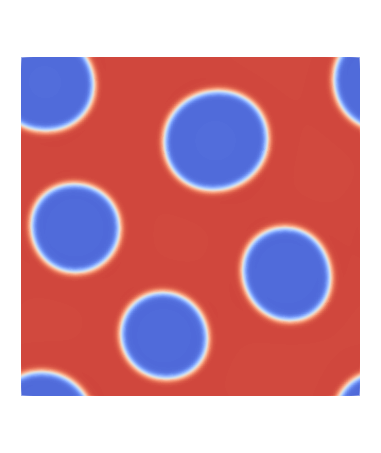} 
& \includegraphics[width=.16\textwidth,trim={0.9cm 1.6cm 0.9cm 1.6cm}, clip]{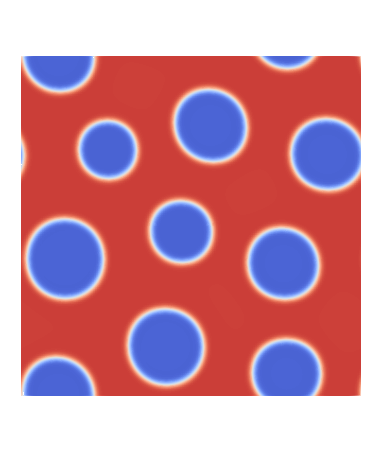} 
& \includegraphics[width=.16\textwidth,trim={0.9cm 1.6cm 0.9cm 1.6cm}, clip]{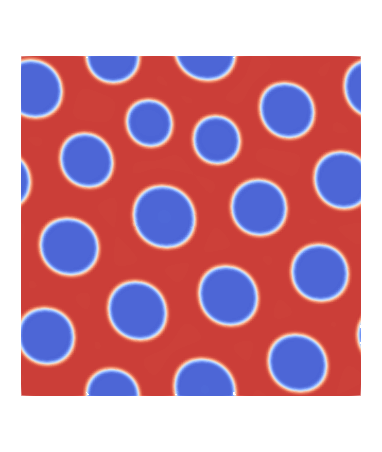} 
\end{tabular}
\end{minipage}%
\hspace{0.035\textwidth}  %% adapt gap
\begin{minipage}{.06\textwidth}
\mbox{}\vspace{.3cm}

\includegraphics[height=0.88\ht\mycolumn]{images/colorbar.png}
\end{minipage}%
\caption{Order parameter $u$ at time $t=1$ for different values of $\kappa \in [0.75\varepsilon, 2\varepsilon]$ and $\Lambda\in[0.2, 2]$. The other parameters are fixed as $\sigma = 1$, $\varepsilon=0.01$, and $\mean{u_0} = 0.3$. 
In contrast to Figure~\ref{fig:test1}, only dotted patterns appear. For fixed $\kappa$, the dots become more numerous and smaller with increasing $\Lambda$. For larger $\kappa$, this transition occurs at higher $\Lambda$.
}\vspace{-1em}
\label{fig:test1b}
\end{figure}

We now repeat the test with a larger mean value for the order parameter, $\mean{u_0} = 0.3$, while keeping the remaining model parameters as in \eqref{NUM:parameters}. The corresponding order parameters $u$ for different values of $\Lambda\in[0.2, 2]$ and $\kappa\in[0.75\varepsilon, 2\varepsilon]$ at time $t=1$ are shown in Figure~\ref{fig:test1b}.
In contrast to Figure~\ref{fig:test1}, only dotted patterns are observed. When $\kappa$ is fixed (i.e., for one row in Figure~\ref{fig:test1b}), the number of dots increases with larger $\Lambda$, while their size decreases. For larger value of $\kappa$, this transition shifts toward higher values of $\Lambda$.

%
%
%   numerical experiment 2
%
%
In the second experiment, we study the influence of the surface tension $\sigma$ and the coupling strength $\Lambda$. We consider the following parameter setting:
\begin{gather}\label{NUM:parameters2}
    \varepsilon = 0.01, \quad 
    \sigma \in[0.01, 10], \quad
    \kappa = 0.2\varepsilon, \quad
    \Lambda \in [0.2, 0.6].
\end{gather}
As before, the initial data are set as $h_0 = 0$, while $u_0$ is given by a random perturbation of amplitude $0.2$ around a prescribed mean value $\mean{u_0} \in \{0.1, 0.3\}$.

% We now fix the mean value $\mean{u_0} = 0.1$ and vary the surface tension $\sigma \in \{0.2, 0.4, \ldots, 2\}$ and the coupling strength $\Lambda \in \{0.2, 0.4, 0.6, 0.8, 1, 2\}$, while the bending rigidity stays fixed $\kappa=0.2\varepsilon$.
% The corresponding order parameters $u$ at time $t=1$ are shown in Figure~\ref{fig:test2}.
% To discuss the effect of varying $\sigma$, we fix $\Lambda$ (i.e., consider one row in Figure~\ref{fig:test2}). We observe ...

\begin{figure}[htb!]
\centering
\begin{minipage}{.86\textwidth}
\setlength{\tabcolsep}{2pt} % default ~6pt
\begin{tabular}{cM{.16\textwidth}M{.16\textwidth}M{.16\textwidth}M{.16\textwidth}M{.16\textwidth}M{.16\textwidth}}
& $\sigma=0.01$ & $\sigma=0.1$  & $\sigma=1.0$ & $\sigma=2.0$ & $\sigma=4.0$ & $\sigma=10.0$   
\\ 
\rotatebox{90}{\hspace{-0.6cm} $\Lambda=0.2$} 
& \includegraphics[width=.16\textwidth,trim={0.9cm 1.6cm 0.9cm 1.6cm}, clip]{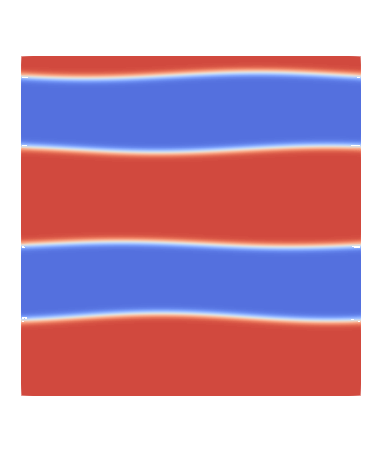} 
& \includegraphics[width=.16\textwidth,trim={0.9cm 1.6cm 0.9cm 1.6cm}, clip]{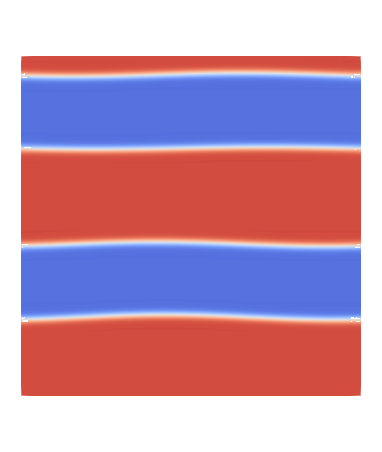} 
& \includegraphics[width=.16\textwidth,trim={0.9cm 1.6cm 0.9cm 1.6cm}, clip]{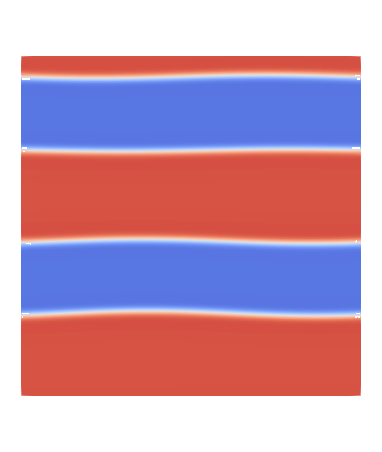} 
& \includegraphics[width=.16\textwidth,trim={0.9cm 1.6cm 0.9cm 1.6cm}, clip]{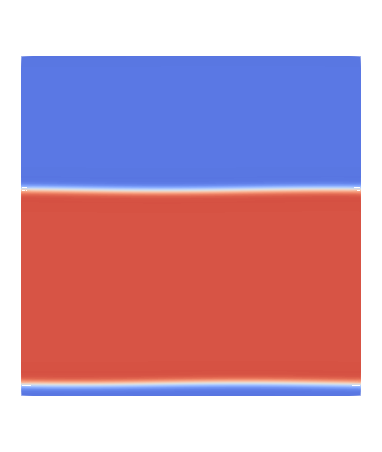} 
& \includegraphics[width=.16\textwidth,trim={0.9cm 1.6cm 0.9cm 1.6cm}, clip]{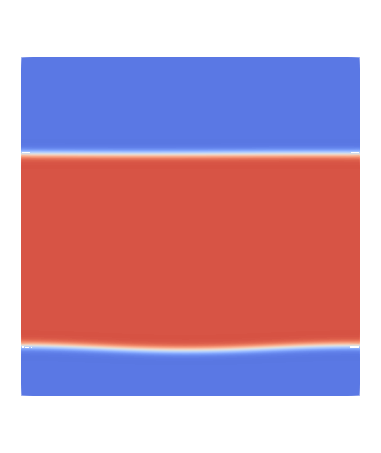} 
& \includegraphics[width=.16\textwidth,trim={0.9cm 1.6cm 0.9cm 1.6cm}, clip]{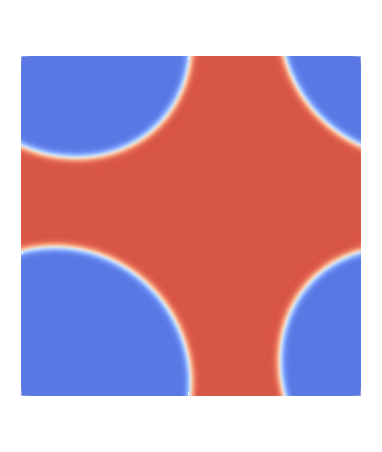} 
\\[-0em]
\rotatebox{90}{\hspace{-0.6cm} $\Lambda=0.3$} 
& \includegraphics[width=.16\textwidth,trim={0.9cm 1.6cm 0.9cm 1.6cm}, clip]{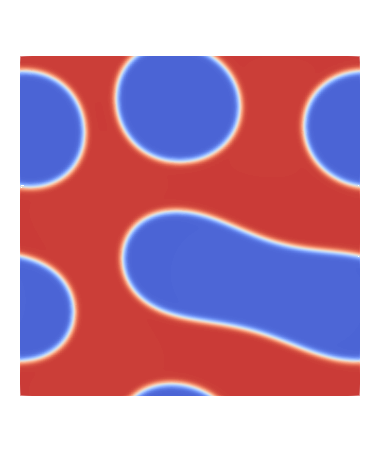} 
& \includegraphics[width=.16\textwidth,trim={0.9cm 1.6cm 0.9cm 1.6cm}, clip]{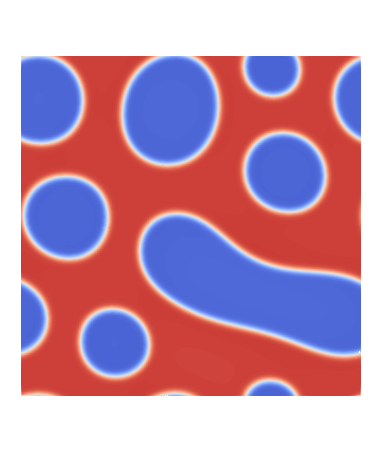} 
& \includegraphics[width=.16\textwidth,trim={0.9cm 1.6cm 0.9cm 1.6cm}, clip]{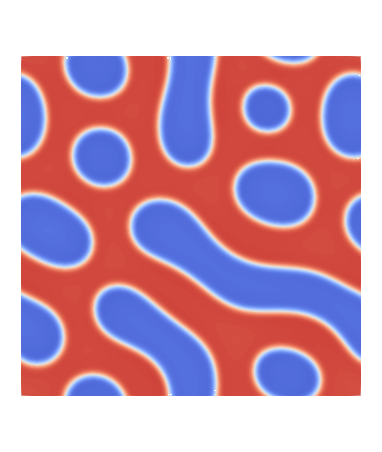} 
& \includegraphics[width=.16\textwidth,trim={0.9cm 1.6cm 0.9cm 1.6cm}, clip]{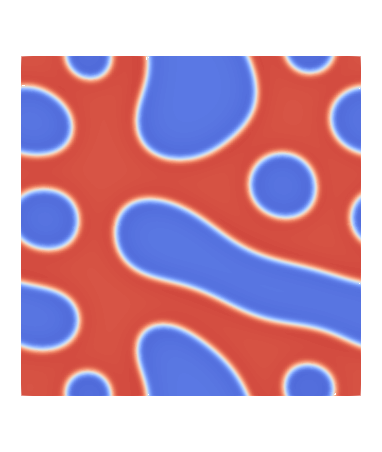} 
& \includegraphics[width=.16\textwidth,trim={0.9cm 1.6cm 0.9cm 1.6cm}, clip]{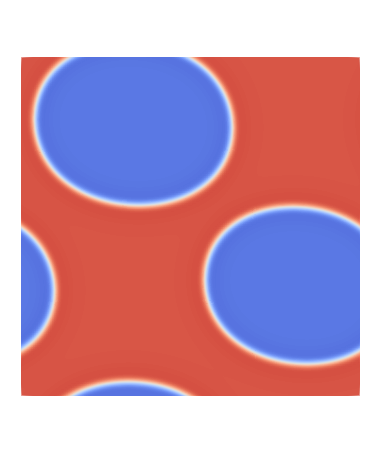} 
& \includegraphics[width=.16\textwidth,trim={0.9cm 1.6cm 0.9cm 1.6cm}, clip]{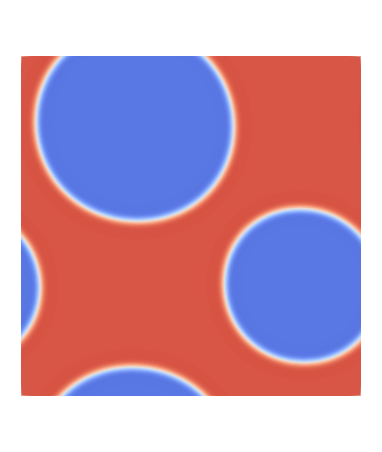} 
\\[-0em]
\rotatebox{90}{\hspace{-0.6cm} $\Lambda=0.4$} 
& \includegraphics[width=.16\textwidth,trim={0.9cm 1.6cm 0.9cm 1.6cm}, clip]{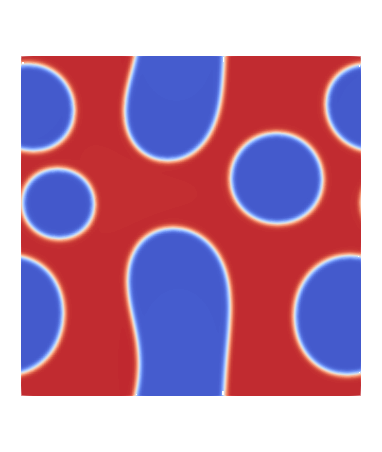} 
& \includegraphics[width=.16\textwidth,trim={0.9cm 1.6cm 0.9cm 1.6cm}, clip]{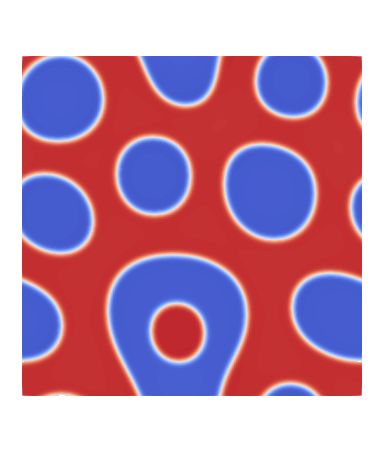} 
& \includegraphics[width=.16\textwidth,trim={0.9cm 1.6cm 0.9cm 1.6cm}, clip]{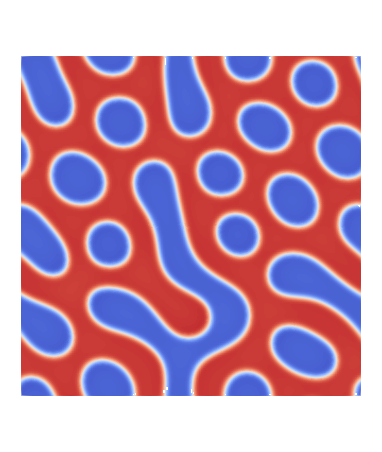} 
& \includegraphics[width=.16\textwidth,trim={0.9cm 1.6cm 0.9cm 1.6cm}, clip]{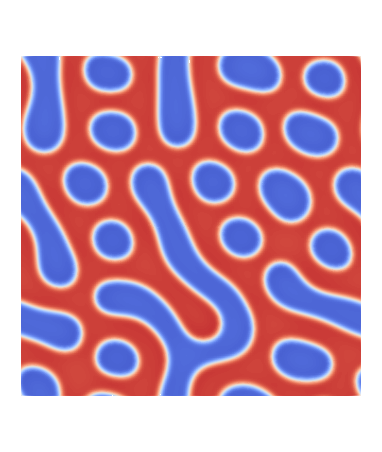} 
& \includegraphics[width=.16\textwidth,trim={0.9cm 1.6cm 0.9cm 1.6cm}, clip]{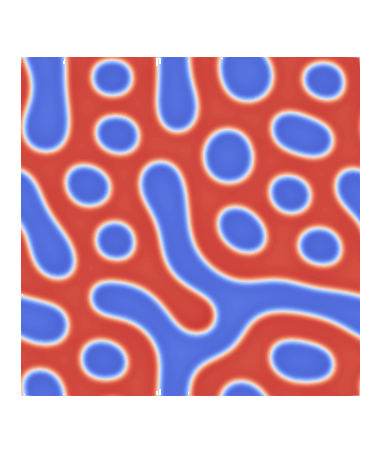} 
& \includegraphics[width=.16\textwidth,trim={0.9cm 1.6cm 0.9cm 1.6cm}, clip]{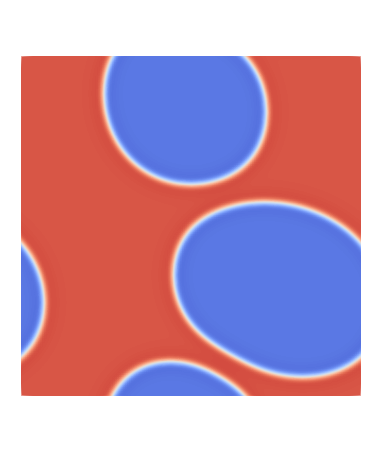} 
\\[-0em]
\rotatebox{90}{\hspace{-0.6cm} $\Lambda=0.6$} 
& \includegraphics[width=.16\textwidth,trim={0.9cm 1.6cm 0.9cm 1.6cm}, clip]{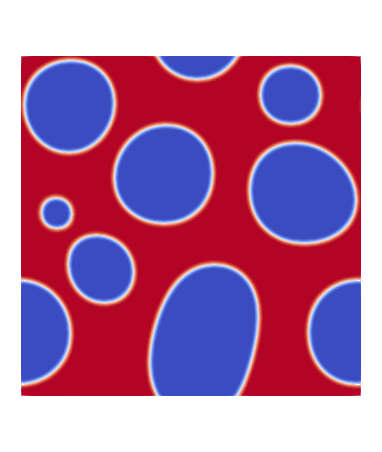} 
& \includegraphics[width=.16\textwidth,trim={0.9cm 1.6cm 0.9cm 1.6cm}, clip]{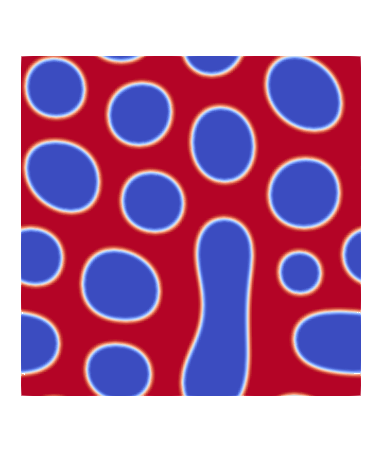} 
& \includegraphics[width=.16\textwidth,trim={0.9cm 1.6cm 0.9cm 1.6cm}, clip]{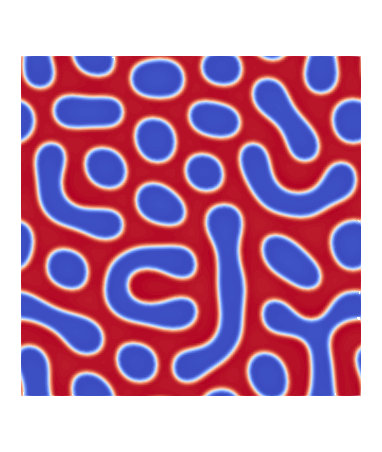} 
& \includegraphics[width=.16\textwidth,trim={0.9cm 1.6cm 0.9cm 1.6cm}, clip]{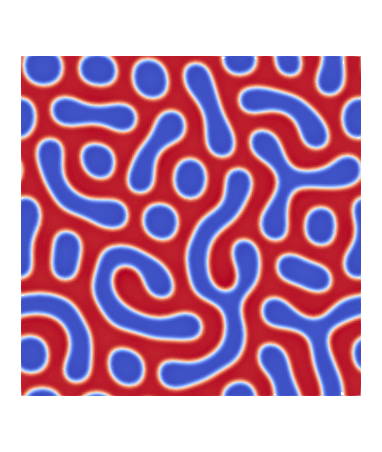} 
& \includegraphics[width=.16\textwidth,trim={0.9cm 1.6cm 0.9cm 1.6cm}, clip]{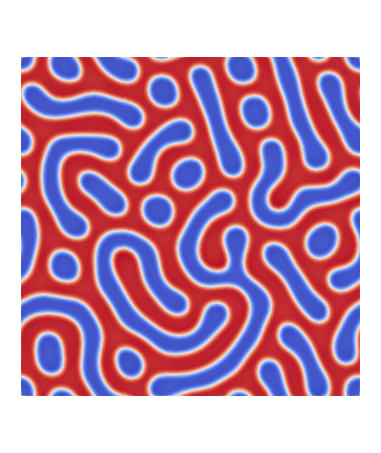} 
& \includegraphics[width=.16\textwidth,trim={0.9cm 1.6cm 0.9cm 1.6cm}, clip]{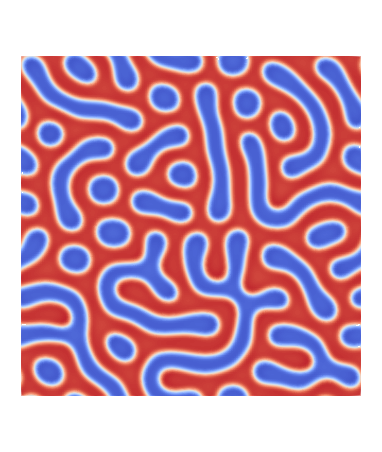} 
\end{tabular}
\end{minipage}%
\hspace{0.035\textwidth}  %% adapt gap
\begin{minipage}{.06\textwidth}
\mbox{}\vspace{.3cm}

\includegraphics[height=0.88\ht\mycolumn]{images/colorbar.png}
\end{minipage}%
\caption{Order parameter $u$ at time $t=1$ for different values of $\Lambda \in [0.2,0.6]$ and $\sigma\in[0.01, 10]$. The other parameters are fixed as $\kappa = 0.2\varepsilon$, $\varepsilon=0.01$, and $\mean{u_0} = 0.1$. Depending on $\Lambda$ and $\sigma$, the system exhibits striped, snake-like, or dotted patterns.}\vspace{-1em}
\label{fig:test2}
\end{figure}

First, we fix the mean value $\mean{u_0} = 0.1$ and vary the coupling strength $\Lambda \in [0.2, 0.6]$ and the surface tension $\sigma \in [0.01, 10]$.
The corresponding order parameter $u$ at time $t = 1$ is shown in Figure~\ref{fig:test2}.
To analyze the influence of $\sigma$, we fix $\Lambda$ (i.e., consider one row in Figure~\ref{fig:test2}).
For $\Lambda = 0.2$, a striped pattern appears for $\sigma \leq 1$. As $\sigma$ increases to $2$ and $4$, the number of stripes decreases while their width grows. For the largest value $\sigma = 10$, the pattern transitions from stripes to dots.
For higher values of $\Lambda$, purely striped structures do not occur. Instead, small $\sigma$ produces a rather dotted pattern that gradually develops into interconnected, snake-like structures as $\sigma$ increases. For large values of $\sigma$, the snake-like domains thicken and eventually merge into larger, dot-like regions, where the dots appear to be rather disordered.
For larger $\Lambda$, this transition occurs at higher $\sigma$ values. In particular, for $\Lambda = 0.6$, the dotted phase becomes dominant for $\sigma \geq 20$, although these cases are not shown in Figure~\ref{fig:test2}.

\begin{figure}[htb!]
\centering
\begin{minipage}{.86\textwidth}
\setlength{\tabcolsep}{2pt} % default ~6pt
\begin{tabular}{cM{.16\textwidth}M{.16\textwidth}M{.16\textwidth}M{.16\textwidth}M{.16\textwidth}M{.16\textwidth}}
& $\sigma=0.01$ & $\sigma=0.1$  & $\sigma=1.0$ & $\sigma=2.0$ & $\sigma=4.0$ & $\sigma=10.0$   
\\ 
\rotatebox{90}{\hspace{-0.6cm} $\Lambda=0.2$} 
% & \includegraphics[width=.16\textwidth]{images/dummy_figure.png}
& \includegraphics[width=.16\textwidth,trim={0.9cm 1.6cm 0.9cm 1.6cm}, clip]{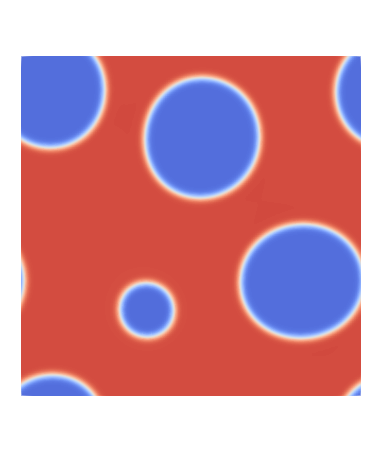} 
& \includegraphics[width=.16\textwidth,trim={0.9cm 1.6cm 0.9cm 1.6cm}, clip]{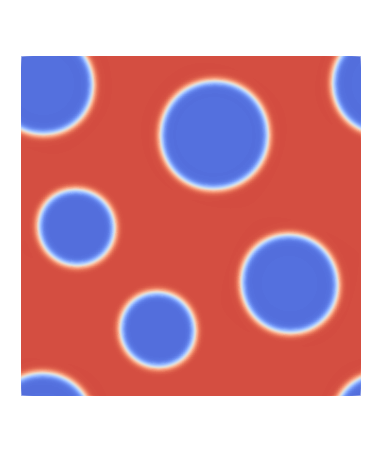} 
& \includegraphics[width=.16\textwidth,trim={0.9cm 1.6cm 0.9cm 1.6cm}, clip]{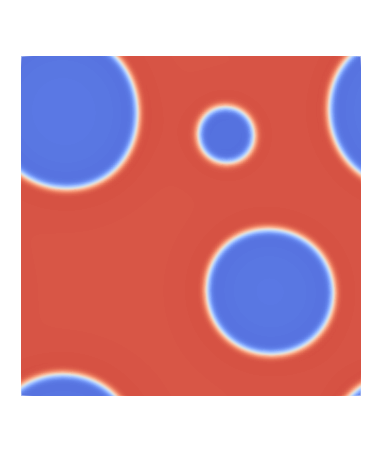} 
& \includegraphics[width=.16\textwidth,trim={0.9cm 1.6cm 0.9cm 1.6cm}, clip]{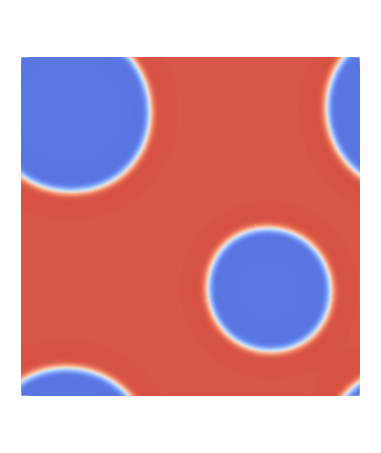} 
& \includegraphics[width=.16\textwidth,trim={0.9cm 1.6cm 0.9cm 1.6cm}, clip]{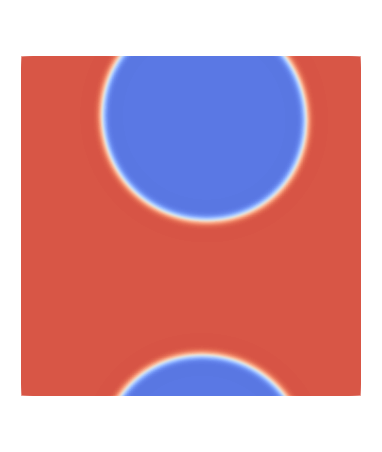} 
& \includegraphics[width=.16\textwidth,trim={0.9cm 1.6cm 0.9cm 1.6cm}, clip]{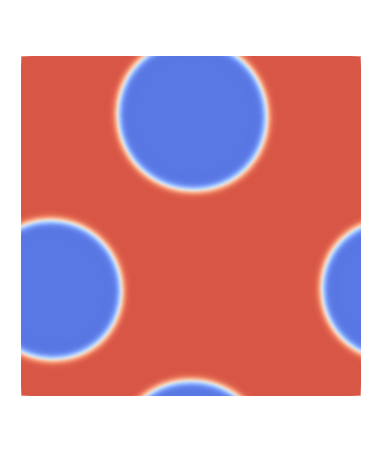} 
\\[-0em]
\rotatebox{90}{\hspace{-0.6cm} $\Lambda=0.3$} 
& \includegraphics[width=.16\textwidth,trim={0.9cm 1.6cm 0.9cm 1.6cm}, clip]{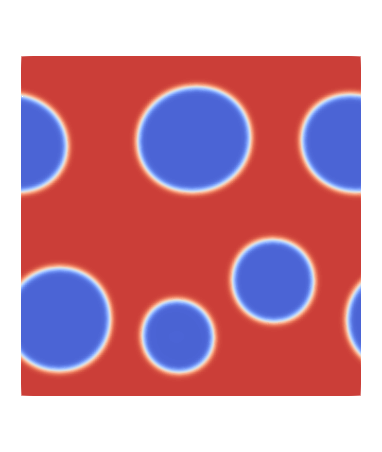} 
& \includegraphics[width=.16\textwidth,trim={0.9cm 1.6cm 0.9cm 1.6cm}, clip]{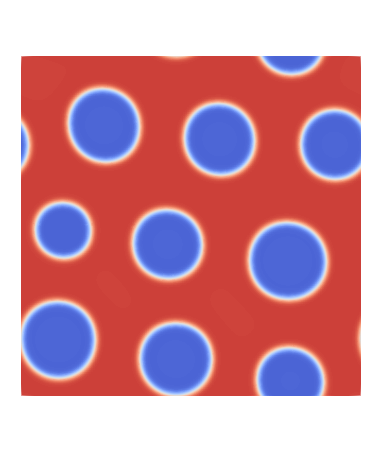} 
& \includegraphics[width=.16\textwidth,trim={0.9cm 1.6cm 0.9cm 1.6cm}, clip]{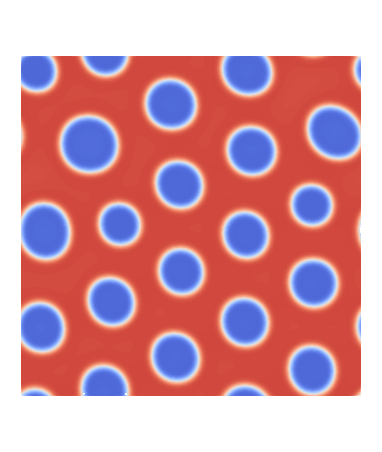} 
& \includegraphics[width=.16\textwidth,trim={0.9cm 1.6cm 0.9cm 1.6cm}, clip]{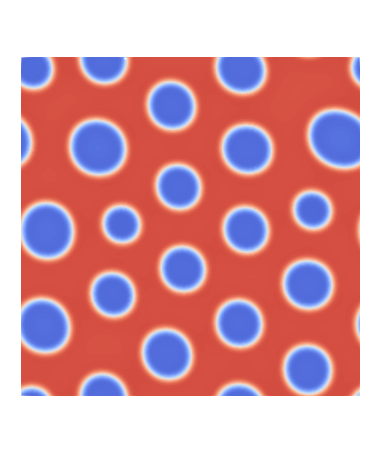} 
& \includegraphics[width=.16\textwidth,trim={0.9cm 1.6cm 0.9cm 1.6cm}, clip]{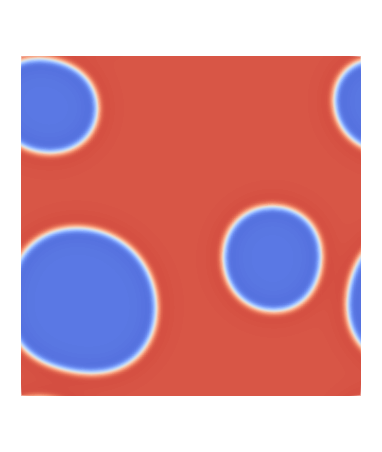} 
& \includegraphics[width=.16\textwidth,trim={0.9cm 1.6cm 0.9cm 1.6cm}, clip]{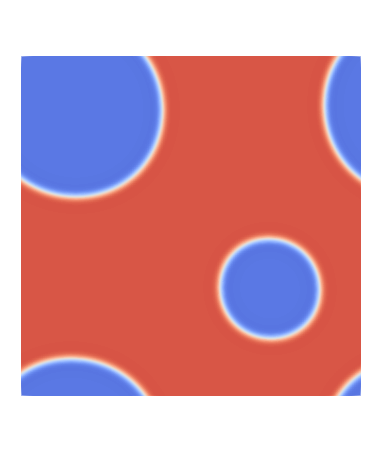} 
\\[-0em]
\rotatebox{90}{\hspace{-0.6cm} $\Lambda=0.4$} 
& \includegraphics[width=.16\textwidth,trim={0.9cm 1.6cm 0.9cm 1.6cm}, clip]{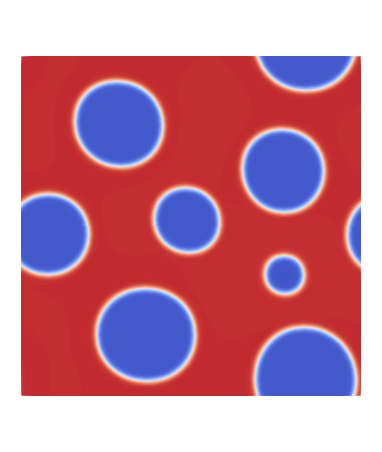} 
& \includegraphics[width=.16\textwidth,trim={0.9cm 1.6cm 0.9cm 1.6cm}, clip]{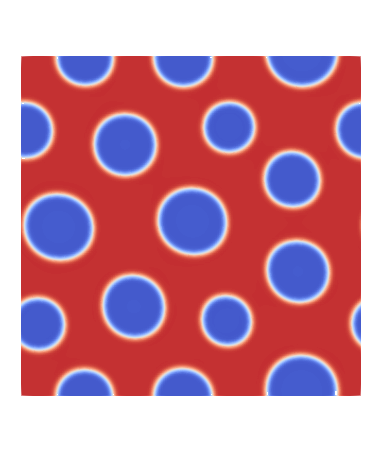} 
& \includegraphics[width=.16\textwidth,trim={0.9cm 1.6cm 0.9cm 1.6cm}, clip]{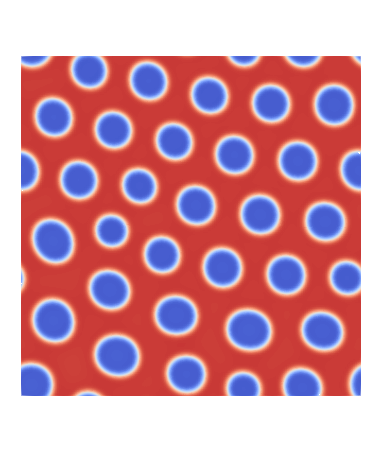} 
& \includegraphics[width=.16\textwidth,trim={0.9cm 1.6cm 0.9cm 1.6cm}, clip]{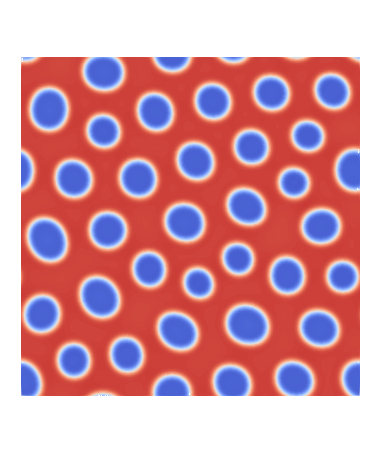} 
& \includegraphics[width=.16\textwidth,trim={0.9cm 1.6cm 0.9cm 1.6cm}, clip]{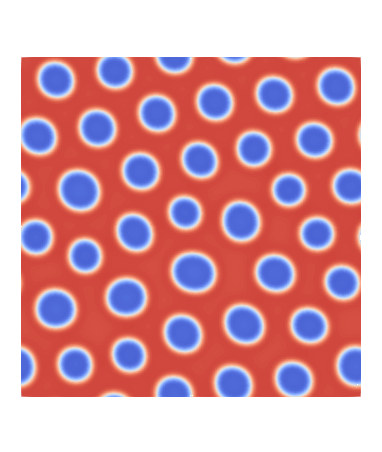} 
& \includegraphics[width=.16\textwidth,trim={0.9cm 1.6cm 0.9cm 1.6cm}, clip]{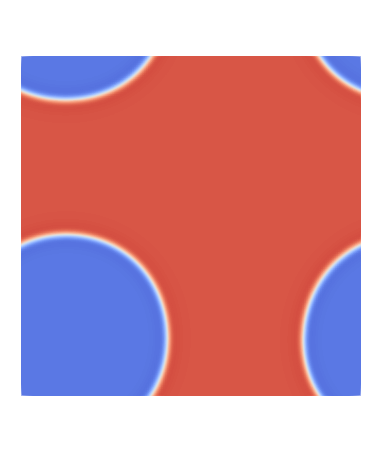} 
\\[-0em]
\rotatebox{90}{\hspace{-0.6cm} $\Lambda=0.6$} 
& \includegraphics[width=.16\textwidth,trim={0.9cm 1.6cm 0.9cm 1.6cm}, clip]{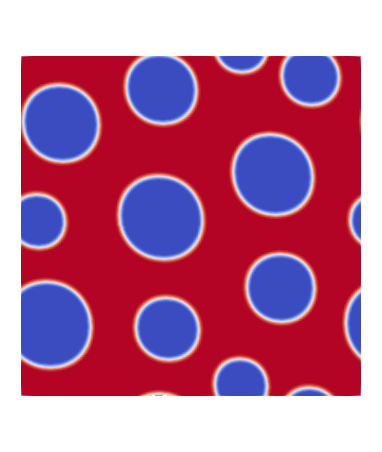} 
& \includegraphics[width=.16\textwidth,trim={0.9cm 1.6cm 0.9cm 1.6cm}, clip]{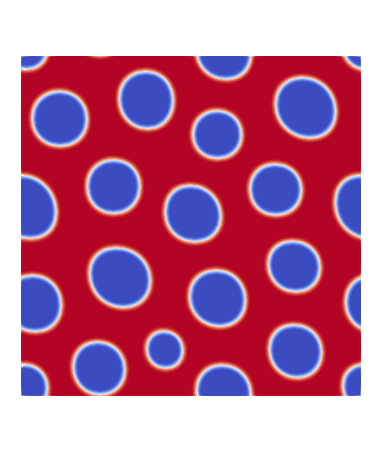} 
& \includegraphics[width=.16\textwidth,trim={0.9cm 1.6cm 0.9cm 1.6cm}, clip]{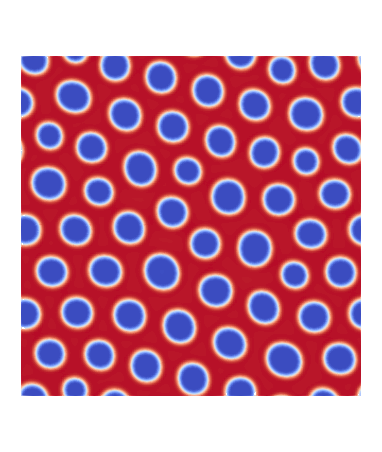} 
& \includegraphics[width=.16\textwidth,trim={0.9cm 1.6cm 0.9cm 1.6cm}, clip]{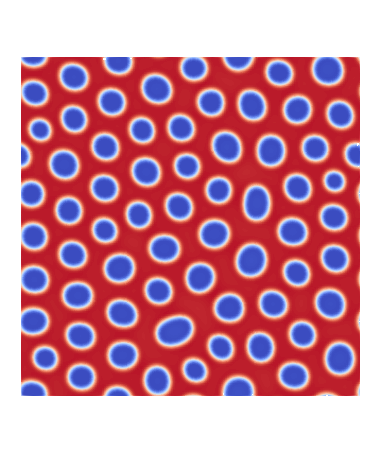} 
& \includegraphics[width=.16\textwidth,trim={0.9cm 1.6cm 0.9cm 1.6cm}, clip]{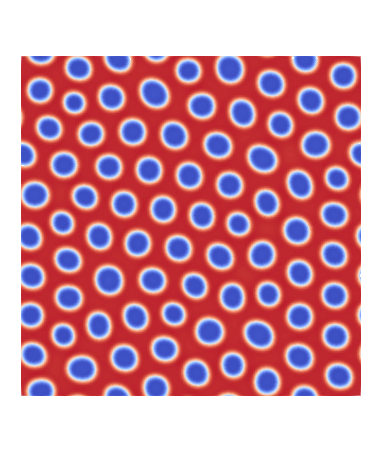} 
& \includegraphics[width=.16\textwidth,trim={0.9cm 1.6cm 0.9cm 1.6cm}, clip]{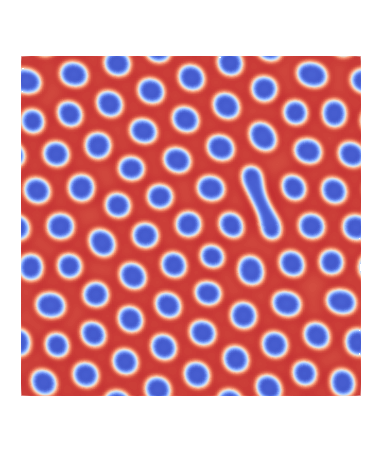} 
\end{tabular}
\end{minipage}%
\hspace{0.035\textwidth}  %% adapt gap
\begin{minipage}{.06\textwidth}
\mbox{}\vspace{.3cm}

\includegraphics[height=0.88\ht\mycolumn]{images/colorbar.png}
\end{minipage}%
\caption{Order parameter $u$ at time $t=1$ for different values of $\Lambda \in [0.2,0.6]$ and $\sigma\in[0.01, 10]$. The other parameters are fixed as $\kappa = 0.2\varepsilon$, $\varepsilon=0.01$, and $\mean{u_0} = 0.3$. In contrast to Figure~\ref{fig:test2}, only dotted patterns are observed. 
For fixed $\sigma$, the dots become more numerous and smaller with increasing $\Lambda$, while varying $\sigma$ leads to a non-monotone change in dot size and number.
%For fixed $\sigma$, the number of dots increases and their size decreases with increasing $\Lambda$. Fixing $\Lambda$, increasing $\sigma$ first reduces the number of dots and enlarges them, but above a critical value of $\sigma$, the number of dots grows again while their size shrinks. For larger $\Lambda$, this transition shifts to higher values of $\sigma$.
}\vspace{-1em}
\label{fig:test2b}
\end{figure}

Similar to before, we repeat the test with a larger mean value of the order parameter, $\mean{u_0} = 0.3$, while keeping the remaining model parameters as in \eqref{NUM:parameters2}.
The corresponding order parameters $u$ at time $t=1$ for different values of $\Lambda \in [0.2, 0.6]$ and $\sigma \in [0.01, 10]$ are shown in Figure~\ref{fig:test2b}.
In contrast to Figure~\ref{fig:test2}, only dotted patterns are observed now.
When $\sigma$ is fixed (i.e., for one column in Figure~\ref{fig:test2b}), the number of dots increases with increasing $\Lambda$, while their size decreases.
To illustrate the influence of $\sigma$, we fix $\Lambda = 0.3$. As $\sigma$ increases from $0.01$ to $2$, the number of dots decreases while their size grows, while they become more and more ordered. Beyond a critical range of $\sigma \in [2, 4]$, the trend reverses: the number of dots increases again, their size becomes smaller.
For larger $\Lambda$, this critical transition appears to shift toward higher $\sigma$ values.
Even though in Figure~\ref{fig:test2} we observed stripes and snake-like patterns for smaller mean values (i.e., $\mean{u_0}=0.1$), a similar critical value of $\sigma$ seems to exist there as well, above which only large dot-like structures remain.

% \newpage

\section*{Acknowledgment}
P.K. and D.T. are partially supported by the Graduiertenkolleg 2339 IntComSin of the Deutsche For\-schungs\-ge\-mein\-schaft (DFG, German Research Foundation) -- Project-ID 321821685.
D.T.~also acknowledges the support from the Swedish Research Council (grant no.~2021-06594) during his stay at Institut Mittag-Leffler, Djursholm, Sweden, in 2025.

\section*{Declaration}
\textbf{Data availability.} Data will be made available on reasonable request.

\bibliographystyle{abbrv}
\bibliography{lipidrafts}

\end{document}